\documentclass[11pt,reqno]{amsart}
\pdfoutput = 1
\usepackage{amsaddr}
\usepackage[utf8]{inputenc}
\usepackage{amsmath,amssymb,amsfonts,epsfig,mathtools,float,latexsym,color}  
\DeclareMathAlphabet{\mathcalold}{OMS}{cmsy}{m}{n}
\usepackage{textcomp}
\usepackage{geometry}
\usepackage{natbib}

\newtheorem{theorem}{Theorem}

\numberwithin{equation}{section}
\numberwithin{table}{section}
\numberwithin{figure}{section}
\allowdisplaybreaks

\newcommand{\qbbar}{\bar{\bar{q}}}

\makeatletter
\renewcommand\@biblabel[1]{}
\makeatother

\title[Froth flotation with drainage]{A degenerating convection-diffusion system \\  modelling froth flotation with drainage}

\author[B\"urger]{Raimund B\"urger}
\address[A1]{\vspace{-0.3cm}CI${}^{\mathrm{2}}$MA and Departamento de Ingenier\'{\i}a Matem\'{a}tica, Facultad de Ciencias F\'{i}sicas y Matem\'{a}ticas, Universidad de Concepci\'{o}n, Casilla 160-C, Concepci\'{o}n, Chile}
\author[Diehl]{\vspace{-0.5cm}Stefan Diehl}
\address[A2]{\vspace{-0.3cm}Centre for Mathematical Sciences, Lund University, P.O.\ Box 118, S-221 00 Lund, Sweden}
\author[Mart\'{\i}]{\vspace{-0.5cm}M.\ Carmen Mart\'{\i}}
\address[A3]{\vspace{-0.3cm}Departament de Matem\`atiques, Universitat de Val\`encia,  Avda. Vicent Andr\'es Estell\'es s/n, Burjassot, Val\`encia, Spain}
\author[V\'asquez]{\vspace{-0.5cm}Yolanda~V\'{a}squez$^*$}
\address[A1]{\vspace{-0.3cm}CI${}^{\mathrm{2}}$MA and Departamento de Ingenier\'{\i}a Matem\'{a}tica, Facultad de Ciencias F\'{i}sicas y Matem\'{a}ticas, Universidad de Concepci\'{o}n, Casilla 160-C, Concepci\'{o}n, Chile\\
\vspace{0.3cm}
 {\small \rm $^*$Corresponding author: \textsf{yvasquez@ing-mat.udec.cl}}}

 \usepackage[font=small,labelfont=bf]{caption}

\begin{document}

\begin{abstract}
{Froth flotation is a common unit operation used  in   mineral processing. It serves
 to separate valuable mineral particles from worthless gangue particles in finely 
  ground ores. The  valuable mineral particles are hydrophobic and 
  attach to bubbles of air injected into the pulp. This creates bubble-particle aggregates that  rise to the top of the flotation column where they accumulate to a froth or foam layer that is removed through a launder for further processing.
   At the same time, the hydrophilic gangue particles settle and are removed continuously.  
  The drainage of liquid due to capillarity is essential  for  the formation of a stable froth layer.
This effect is included into a previously formulated 
  hyperbolic system of partial differential equations that models the volume fractions of floating aggregates and settling hydrophilic solids [R.\ B\"{u}rger, S.\ Diehl and M.C.\ Mart\'i, {\it IMA J.\ Appl.\ Math.} {\bf 84} (2019) 930--973]. 
The construction of desired steady-state solutions with a froth layer is detailed and feasibility conditions on the feed volume fractions and the volumetric flows of feed, underflow and wash water are visualized in so-called operating charts. A monotone numerical scheme is derived and employed to simulate  the  dynamic behaviour of a flotation column. It is also proven that, under a suitable Courant-Friedrichs-Lewy (CFL) condition, the approximate volume fractions are bounded between zero and one when the initial data are.}
\\[5pt]
\\[5pt]
\textbf{2000 Math Subject Classification:} 35L65, 35P05, 35R05.
\end{abstract}

\maketitle


\section{Introduction}

\subsection{Scope}
\begin{figure}[t]
\begin{center}
\includegraphics[width=0.85\textwidth]{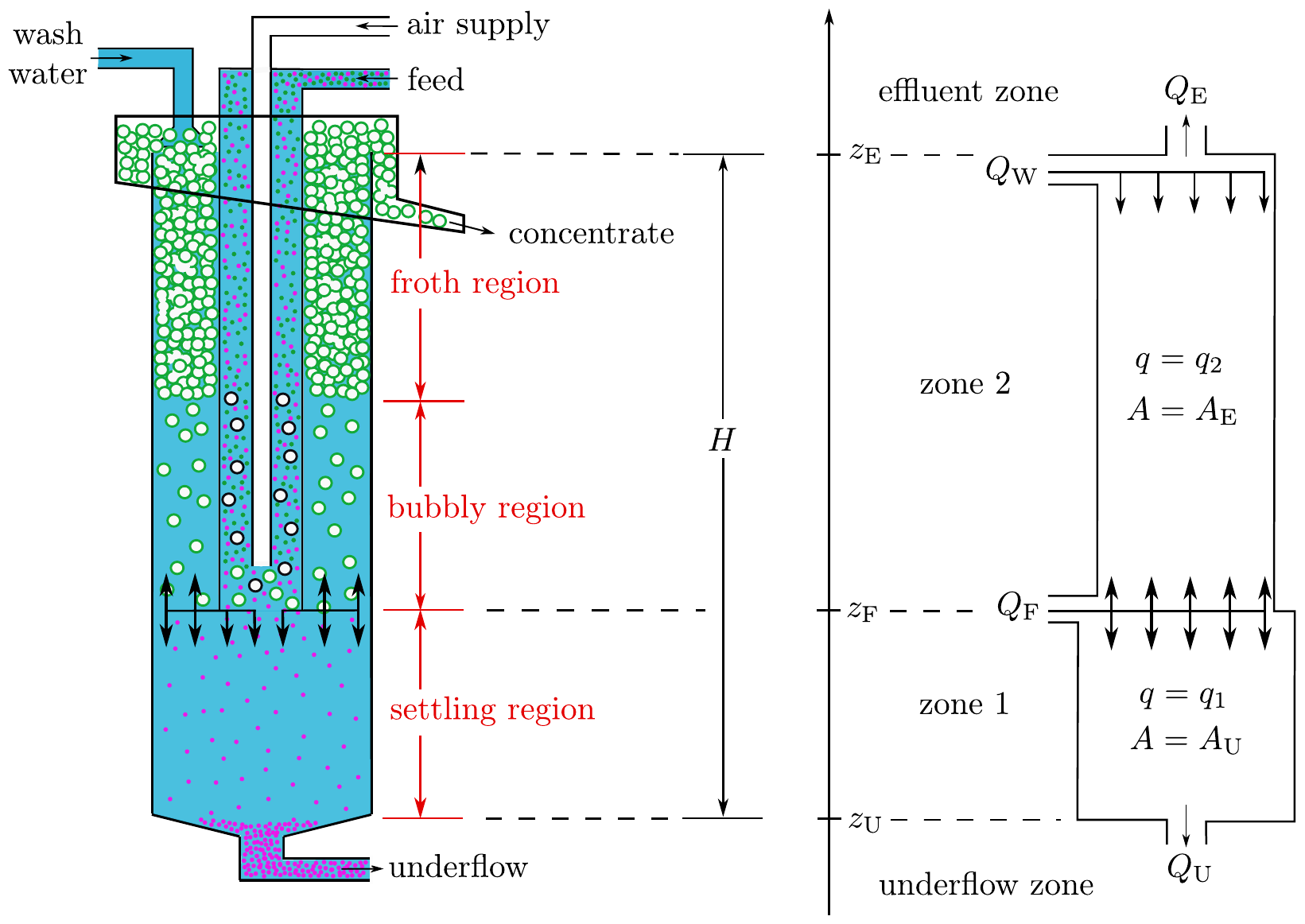}
\caption{Left: Schematic of a flotation column; cf.\ the Reflux Flotation Cell by  \cite{Dickinson2014}.
Right: The corresponding one-dimensional conceptual model with a non-constant cross-sectional area $A(z)$. Wash water is sprinkled at the effluent level $z=z_{\mathrm{E}}$ and a mixture of aggregates and feed slurry is fed at $z=z_{\mathrm{F}}$, where $z_{\mathrm{U}}<z_{\mathrm{F}}<z_{\mathrm{E}}$ divide the real line into the zones inside the column and the underflow and effluent zones.}
\label{fig:Column}
\end{center}
\end{figure}

Flotation is a separation process where air bubbles are used to attract hydrophobic particles or droplets from a mixture of solids in water.
The process is used in mineral processing, where valuable mineral particles are separated out from crushed ore, and in wastewater treatment to remove floating solids, residual chemicals,  and droplets of oil and fat.
The process is often applied in a column to which both a mixture of particles (or droplets) and air bubbles are injected.
The effluent at the top should consist of a concentrate of hydrophobic particles that are attached to the bubbles, while the hydrophilic particles settle to the bottom, where they are removed (see Figure~\ref{fig:Column}). 
A layer of froth at the top is preferred since the effluent then consists of a minimum amount of water and the froth works as a filter enhancing the separation process.
In our previous models of column froth flotation (with or without simultaneous sedimentation of hydrophilic particles) \citep{SDNHM1,SDIMAflot2019,SDmineng3_flot_sed,SDwst_DAF_CFF}, a particular  constitutive assumption on the bubble velocity  leads to a hyperbolic system of partial differential equations (PDEs) that  models the layer of froth with a constant horizontal average volume fraction of bubbles~$\phi$, or equivalently, the volume fraction $\varepsilon=1-\phi$ of liquid (or suspension with hydrophilic particles)   that fills  the interstices outside the bubbles.
It is however known that~$\varepsilon$ varies with the height in the froth because of capillarity and drainage of liquid; see~\cite{Neethling2018} and references therein.

It is the purpose of this work to extend the previous hyperbolic model to one that includes capillarity.
 To this end we partly generalize the well-known drainage equation to hold for all bubble volume fractions, and partly generalize our previous model of column froth flotation with simultaneous sedimentation.
The latter is a nonlinear system of PDEs where the unknowns are the volume fractions of aggregates (bubbles/droplets loaded with hydrophobic particles) and solid hydrophilic particles.
A  numerical scheme for the new governing PDE system is presented. 
We show that the approximate volume fractions stay between zero and one if a suitable Courant-Friedrichs-Lewy (CFL)  condition is used.
Furthermore, we construct desired steady-state solutions and provide algebraic equations and inequalities 
 that establish the dependence of steady states on the input and control variables.
Such dependences are conveniently visualized in so-called operating charts that  constitute a graphical tool for controlling the process. The particular importance of steady states comes from the application under study; namely they describe 
 the ability of the model to capture steady operation of the flotation device without the necessity of permanent 
  control actions. 

\subsection{Some preliminaries} 

Froth is assumed to form when the volume fraction of bubbles~$\phi$ is above a critical value~$\phi_\mathrm{c}=1-\varepsilon_\mathrm{c}$ when the bubbles are in contact with each other.
Then capillarity forces are involved, which means that the governing PDE is parabolic, whereas it is hyperbolic in regions without froth. The present 
 derivation is based on the traditional one by~\cite{Leonard1965laminar}, \cite{Goldfarb1988}, and \cite{Verbist1996foam}, leading to the drainage equation for low liquid content~$\varepsilon$.
We then combine results by~\cite{Neethling2002foam} and~\cite{Stevenson2004} to obtain a constitutive relationship between the relative fluid-gas velocity~$u$ and the liquid volume fraction~$\varepsilon\leq\varepsilon_\mathrm{c}$ when capillarity forces are present.
With a compatibility condition at~$\varepsilon_\mathrm{c}$, we obtain a constitutive relationship of the relative fluid-bubble velocity~$u$ as a function of~$\varepsilon\in[0,1]$, which for~$\varepsilon>\varepsilon_\mathrm{c}$ is the common  \cite{Richardson1954} power-law expression for separated bubbles~\citep{Galvin2014II}.
The resulting generalized drainage PDE is (in a closed vessel)
\begin{equation}\label{eq:gen_drainage}
\partial_t \varepsilon 
-\partial_z \bigl(\varepsilon\tilde{v}_\mathrm{f}(\varepsilon) \bigr) 
= \partial_z^2\tilde{D}(\varepsilon),
\end{equation}
where $t$~is time, $z$~is height, $\tilde{v}_\mathrm{f}(\varepsilon)$~is a nonlinear fluid-velocity function,  and 
$\tilde{D}(\varepsilon)$ an integrated diffusion function modelling capillarity, which is zero for~$\varepsilon>\varepsilon_\mathrm{c}$. 

Equation~\eqref{eq:gen_drainage} can alternatively be written in terms of the volume fraction of bubbles~$\phi$.
We assume that~$\phi$ denotes the volume fraction of aggregates, by which we mean bubbles that are fully loaded with hydrophobic particles.
Under  a common constitutive assumption for the settling of hydrophilic particles within the liquid outside the bubbles,  the following system of PDEs  models the combined flotation-drainage-sedimentation process in a vertical column with a feed inlet of air-slurry mixture at the height~$z=z_\mathrm{F}$ with the volumetric flow~$Q_\mathrm{F}(t)$ (see Figure~\ref{fig:Column}):
\begin{align}\label{eq:gov} \begin{split} 
& A(z) \partial_t \begin{pmatrix}
\phi \\ 
\psi\end{pmatrix} +  
\partial_z \left(A(z)\begin{pmatrix}
J(\phi,z,t) \\
-\tilde{F}(\psi,\phi,z,t)
\end{pmatrix}\right) 
 \\ & = 
 \partial_z \left(A(z)\gamma(z) \partial_z  D(\phi) \begin{pmatrix}
1 \\
- {\psi}/{(1-\phi )}
\end{pmatrix}\right)
+
Q_\mathrm{F}(t)\begin{pmatrix}
\phi_\mathrm{F}(t)\\
\psi_\mathrm{F}(t)\end{pmatrix}\delta_{z_\mathrm{F}}.
\end{split} 
\end{align} 
Here, $\psi$ is the volume fraction of solids, $A(z)$ the cross-sectional area of the tank, and $J$ and $\tilde{F}$ are convective flux functions that depend discontinuously on~$z$ at the locations of the feed and wash water inlets and the outlets at the top and bottom.  The system \eqref{eq:gov} is valid for~$t>0$ and all $z\in\mathbb{R}$ where the characteristic function~$\gamma(z)=1$ indicates the interior  of the tank and $\gamma(z)=0$ outside, and $\delta$ is the delta function.
Outside the tank, the mixture is assumed to follow the outlet streams; consequently, boundary conditions are not needed; conservation of mass determines the outlet volume fractions in a natural way.

Similarly to the role of~$\tilde{D}$ in \eqref{eq:gen_drainage}, 
the nonlinear function~$D$ models the capillarity present when bubbles are in contact. 
 Precisely, with a function $d(\phi)$ (specified later) we define 
 \begin{align}  \label{Dphidef} 
  D( \phi) := \int_0^{\phi}  d(s) \, \mathrm{d} s. 
  \end{align} 
   The function~$d$  is assumed to satisfy 
   \begin{align} \label{by3.2}  
    d( \phi) = D' (\phi) = \begin{cases}  
     0 & \text{for $0 \leq \phi \leq \phi_{\mathrm{c}}$,} \\
      >0 & \text{for $\phi_{\mathrm{c}} < \phi \leq 1$.} \end{cases} 
\end{align} 
Consequently, at each point $(z, t)$ where $\phi(z,t) \leq \phi_{\mathrm{c}}$, there  holds $D( \phi(z, t)) =0$, 
 and therefore  \eqref{eq:gov} degenerates at such points into a first-order system of conservation laws 
   of hyperbolic type (as was shown in \cite{SDIMAflot2019}). Since this degeneration occurs for 
    $0 \leq \phi \leq \phi_{\mathrm{c}}$ and $0 \leq \psi 
      \leq 1 - \phi$,  
        that is, on a set of positive two-dimensional measure, \eqref{eq:gov} is called {\em strongly} degenerate. 
         While it is clear that the first PDE in \eqref{eq:gov}  is parabolic for $\phi_{\mathrm{c}} < \phi \leq 1$ 
          and this PDE, as well as \eqref{eq:gen_drainage}, are scalar strongly degenerate parabolic equations, 
the same cannot be said about the system.
We observe namely that with $A=\gamma=1$, the diffusion term on the right-hand side can be written as 
            \begin{align*} 
             \partial_z \biggl( d(\phi) \partial_z \phi \begin{pmatrix} 1  \\ - \psi/(1-\phi) \end{pmatrix} 
              \biggr)  = \partial_z \biggl( \boldsymbol{B} ( \phi, \psi) \begin{pmatrix} 
               \partial_z \phi \\ \partial_z \psi \end{pmatrix} \biggr) , \quad 
                \boldsymbol{B} ( \phi, \psi) := d( \phi) \begin{bmatrix} 1 & 0 \\ - \psi/(1-\phi) & 0 \end{bmatrix}. 
              \end{align*} 
              Since at least one of the eigenvalues of $\boldsymbol{B} ( \phi, \psi)$ is always zero, we observe that even 
               when  $d( \phi) >0$, the system \eqref{eq:gov} is not {\em strictly}  parabolic. 

\subsection{Related work}
Modelling flotation and developing strategies to control this process are  research areas 
  that have generated 
 many contributions; see~\cite{Cruz1997} and  the review  by~\cite{Quintanilla2021} and references therein. 
 The development of  control strategies requires  dynamic models along with a
   categorization of steady-state (stationary) solutions of such models. 
Since the volume fractions depend on both time and space, the resulting governing equations are PDEs.
With the aim of developing controllers, \cite{Tian2018,Tian2018II} and~\cite{Azhin2021mineng,Azhin2021} use hyperbolic systems of PDEs for the froth or pulp regions coupled to ODEs for the lower part of the column.
They include the attachment and detachment processes; however, the phases seem to have constant velocities, which is not in agreement with the established drift-flux theory~\citep{Wallis1969,Rietema1982,Brennen2005book}.
Nonlinear dependence of the phase velocities on the volume fractions give rise to discontinuities in the concentration profiles, which is confirmed experimentally~\citep{Cruz1997}.
 \cite{Azhin2021mineng,Azhin2021} show continuous steady-state profiles for their model.

\cite{Narsimhan2010} shows realistic conceptual transient solutions of a bubble-liquid suspension which is homogeneous initially. 
The rising bubbles form a layer of foam at the top which can undergo compressibility due to gravity and capillarity.
The governing PDE has similarities with~\eqref{eq:gen_drainage}; however, it is not utilized in all the different regions of solution.
Separate equations are derived for the foam region and boundary assumptions between regions have to be imposed.
The purpose of our previous and the present contributions is to let a single equation such as \eqref{eq:gen_drainage} govern the bubble-liquid behaviour under any dynamic situation.

Phenomenological models for two-phase systems with bubbles rising (or, analogously, particles settling) in a liquid, are derived from physical laws of conservation of mass and momentum~\citep{Bascur1991,Bustos1999, Burger&W&C2000,Brennen2005book}.
Under  certain simplifying assumptions on the stress tensor and partial pressure of the bubbles/solids, one can obtain 
first- or second-order PDEs involving one or two constitutive (material specific) functions, respectively. 

The resulting first-order PDE modelling such a separation process in a one-dimensional column of rising bubbles is a scalar conservation law $\phi_t+(\phi\tilde{v}(\phi))_z$=0 with a  drift-flux velocity function~$\tilde{v}(\phi)$.
This is in agreement with the drift-flux theory by~\cite{Wallis1969}.
With additional bulk flows due to the inlets and outlets of the column that theory has mostly been used for steady-state investigations of flotation columns~\citep{Vandenberghe2005,Stevenson2008,Dickinson2014, Galvin2014II, Galvin2014III}.
Models of and numerical schemes for column froth flotation with the drift-flux assumption and possibly simultaneous sedimentation  have been presented by the authors \citep{SDNHM1,SDIMAflot2019,SDmineng3_flot_sed,SDwst_DAF_CFF}.

The analogy of the drift-flux theory  for sedimentation is the established solids-flux theory~\citep{Kynch1952, Ekama1997, SDjem1, SDwatres2,  Ekama2004, LaMotta2007}.
With an additional constitutive assumption on sediment compressibility, the model becomes a second-order degenerate parabolic PDE~\citep{Burger&W&C2000}.
Sedimentation in a clarifier-thickener unit  is mathematically similar to the column-flotation case.
A full PDE model of such a vessel necessarily contains source terms and spatial discontinuities at both inlets and outlets.
Steady-state analyses, numerical schemes, dynamic simulations and control of such models can be found in~\cite{Burger&K&R&T2004b, Burger&K&T2005a} and \cite{SDsiam2,SDsiam3,SDjem5}.
Because of the discontinuous coefficients and degenerate diffusion term of the PDE, so-called entropy conditions are needed to guarantee a unique physically relevant solution~\citep{Burger&K&T2005a, Evje2000, SDjhde1}.
Those results will be utilized in the present work.

The first-order PDE of the flotation process does not include capillarity in the foam.
Such effects have been studied intensively by~\cite{Neethling2002foam, Neethling2003} and~\cite{Neethling2018}; see more references in~\cite{Quintanilla2021}.
Solids motion in froth can be found in~\cite{Neethling2002_solids_froth}.
It is the aim of the present contribution to extend our previous hyperbolic PDE model to include capillarity.

\subsection{Outline of the paper} 

In Section \ref{sec:drainage}, we consider a simplified two-phase bubble-fluid system in a closed vessel and derive a generalized drainage equation governing the flotation of the bubbles with formation of froth and drainage of liquid from it.
In Section \ref{sec:model}, we extend the equation derived to the process of column flotation with sedimentation of solid particles and with froth drainage at the top. The treatment of the feed inlets and the definition of the flux density functions in each zone (see Figure~\ref{fig:Column}) are detailed.  
Section \ref{sec:steadystates} is devoted to the construction of steady-state solutions having a froth layer at the top of the tank and  bubble-free   underflow. Necessary conditions for those so-called desired steady-states to appear, in terms of inequalities involving the volumetric flows $Q_\mathrm{U}$, $Q_\mathrm{F}$ and $Q_\mathrm{W}$ and the incoming volume fractions of aggregates $\phi_\mathrm{F}$ and solids $\psi_\mathrm{F}$, are derived  in Sections~\ref{sec:constructionSS} and \ref{sec:desiredSS}.
In Section \ref{sec:num_method}, the numerical scheme for simulation of the process is introduced. 
It is proven that  under a CFL condition, the approximate volume fractions of aggregates and solids remain between zero and one provided that the initial data do.
The proofs are outlined  in Appendix~A.
Some simulations are provided in Section~\ref{sec:numexp}.
They show fill-up of a flotation column and froth formation, illustrating the response of the system to changes of operating conditions. Finally, conclusions are drawn in Section~\ref{sec:conclusions}.

\section{A generalized drainage equation in a closed tank} \label{sec:drainage}

The two-phase system has bubbles of volume fraction~$\phi$ and phase velocity~$v$, and fluid of volume fraction~$\varepsilon=1-\phi$ and phase velocity~$v_\mathrm{f}$, where $0\leq\phi,\varepsilon\leq 1$.
When the bubbles are mono-sized and separated from each other (i.e., there is no froth), a common expression for their velocity 
  in a closed container without any bulk flow is~\citep{Pal1989, Vandenberghe2005}
\begin{equation*}
v(\phi)= v_\mathrm{term}(1-\phi)^{n_{\mathrm{b}}}
\qquad\text{(separated bubbles)},
\end{equation*}
where $v_\mathrm{term}$ is the velocity of a single bubble far away from others ($\phi\approx 0$) and~$n_{\mathrm{b}}$ a dimensionless parameter (similar to the Richardson-Zaki exponent within the analogous expression for the sedimentation of  mono-sized and equal-density particles in a liquid, see Section~\ref{subsec:threephase}).
We thus let velocities be positive in the upward direction of the $z$-axis.
The relative velocity of fluid to bubbles is $u:=v_\mathrm{f}-v$.
In a closed container, the volume-average velocity is zero; hence, $0=\phi v + \varepsilon v_\mathrm{f}$, and we get
\begin{equation}\label{eq:power}
u=-\frac{(1-\varepsilon)v}{\varepsilon}-v
=-\frac{v}{\varepsilon} = -v_\mathrm{term}\varepsilon^{n_{\mathrm{b}}-1}
\qquad\text{(separated bubbles)},
\end{equation}
which is negative because the fluid flows downwards.
We also obtain  the  identities  
$v_\mathrm{f}=(1-\varepsilon)u$ and $v=-(1-\phi)u$.

If $\phi$ exceeds 
 a certain critical volume fraction~$\phi_\mathrm{c}=1-\varepsilon_\mathrm{c}$, the bubbles touch each other and a foam is formed.
The larger~$\phi>\phi_\mathrm{c}$, or smaller~$\varepsilon<\varepsilon_\mathrm{c}$, the more deformed are the bubbles.
Randomly packed rigid spheres leave a volume fraction of $\varepsilon_\mathrm{c}=1-0.64=0.36$; cf.\ 
\cite[Table~1]{Brito2012advantages}.
For froth, we assume the value $\varepsilon_{\mathrm{c}} = 0.26$ \cite[Eq.~(21)]{Neethling2003} and \cite{Narsimhan2010}.

We  discuss  below  the most difficult intermediate fluid volume fractions 
 when $\varepsilon$~is smaller than, but close to~$\varepsilon_\mathrm{c}$.
We consider, however, first a layer of foam with a very low volume fraction of liquid~$\varepsilon$ and 
  recall the derivation of the drainage equation \citep{Leonard1965laminar,Goldfarb1988,Verbist1996foam}.
 In this case the deformed bubbles are separated by very thin lamellae, which are separated by channels, so called \emph{Plateau borders}, which are connected at vertices, or nodes, so that a network is formed.
It is assumed that almost all the liquid is contained in the Plateau borders, whose cross section is the plane   region bounded by three externally tangential circles all of radius~$r$.
This deformed triangular-shaped region has the area
\begin{equation}\label{eq:calA}
\mathcalold{A}=C^2 r^2,\quad\text{with}\quad
C:=\bigl(\sqrt{3}-\pi/2 \bigr)^{1/2}.
\end{equation}
If the radius~$r$ changes along the Plateau border, this is related to a pressure difference according to the Young-Laplace law:
\begin{equation*}
p_\mathrm{f}=p_\mathrm{b}-\frac{\gamma_\mathrm{w}}{r},
\end{equation*}
where $p_\mathrm{f}$ and $p_\mathrm{b}$ are the fluid and bubble pressure, respectively, and $\gamma_\mathrm{w}$ is the surface tension of water.
The bubble pressure~$p_\mathrm{b}$ is assumed to be constant.

There are three forces acting per volume fraction of the Plateau border:
\begin{alignat*}2
&\text{gravity:}&\quad&\rho_\mathrm{f}\boldsymbol{g},\\
&\text{dissipation:}&\quad&-\frac{C_\mathrm{PB}\mu}{\mathcalold{A}} \boldsymbol{u} 
= -\frac{C_\mathrm{PB}\mu}{C^2 r^2} \boldsymbol{u},\\
&\text{capillarity:}&\quad&-\nabla p_\mathrm{f} = -\frac{\gamma_\mathrm{w}}{r^2}\nabla r.
\end{alignat*}
Here, $\boldsymbol{g}$ is the gravity acceleration vector, $\boldsymbol{u}$~the fluid-bubble relative velocity, $\mu$~the fluid viscosity, and $C_\mathrm{PB}$the dimensionless Plateau border drag coefficient, which can be inferred to be 49.3 from the numerical calculations by~\cite{Leonard1965laminar}. 
The value $C_\mathrm{PB}=50$ is often used in the literature.
The sum of the three forces is zero if one neglects inertial forces.
Along a Plateau border tilted an angle~$\theta$ from the vertical $z$-axis, we place a $z_\theta$-axis with the coordinate relation $z=z_\theta\cos \theta$.
The force balance along the $z_\theta$-axis is
\begin{equation*}
-\rho_\mathrm{f}g\cos\theta 
-\frac{C_\mathrm{PB}\mu}{C^2 r^2}
u_\theta 
-\frac{\gamma_\mathrm{w}}{r^2} \partial_{z_\theta} r  =0,
\end{equation*}
where the relative fluid-gas velocity in the channel is $u_\theta={u_\mathrm{PB}}/{\cos\theta}$ and $u_\mathrm{PB}$~is its vertical contribution from one Plateau border of angle~$\theta$, which thus is
\begin{equation*}
u_\mathrm{PB} = -\frac{C^2r^2}{C_\mathrm{PB}\mu}\left(
\rho_\mathrm{f}g
+ \frac{\gamma_\mathrm{w}}{r^2} \partial_z  r 
\right)\cos^2\theta.
\end{equation*}
Under the assumption of randomly distributed Plateau borders with respect to the angle $0\leq\theta\leq\pi$, the likelihood that a Plateau border has an angle in the interval $(\theta,\theta +\mathrm{d}\theta)$ is the area $2\pi\sin\theta\,\mathrm{d}\theta$ of the circular strip of the unit sphere divided by its total area $4\pi$.
Since
\begin{equation*}
\langle\cos^2 \rangle
:=\int_{0}^{\pi}\cos^2\theta\,\frac{2\pi\sin\theta}{4\pi}\mathrm{d}\theta
=\frac{1}{3},
\end{equation*}
the relative vertical velocity~$u$ is defined as the average vertical relative fluid-gas velocity for many Plateau borders:
\begin{equation}\label{eq:uPBaver}
u=\langle u_\mathrm{PB}\rangle = -\frac{C^2r^2\rho_\mathrm{f}g}{3C_\mathrm{PB}\mu}\left(
1
+ \frac{\gamma_\mathrm{w}}{r^2\rho_\mathrm{f}g} \partial_z r 
\right).
\end{equation}
This velocity can be expressed in~$\mathcalold{A}$ by~\eqref{eq:calA}, and substituting the resulting expression into the conservation law $\partial_t\mathcalold{A}+\partial_z(\mathcalold{A}(1-\varepsilon)u)=0$ and setting $\varepsilon=0$ (recall that  $v_\mathrm{f}=(1-\varepsilon)u$) one obtains the classical drainage equation for low liquid content.

We want an equation for the volume fraction~$\varepsilon$, which is equal to~$\mathcalold{A}$ times the length of Plateau borders per unit volume; cf.~\cite{Neethling2018}.
Then the length~$L$ and number of such channels should be estimated.
Since we also want an equation for all $0\leq\varepsilon\leq\varepsilon_\mathrm{c}$, the estimation of such numbers becomes difficult since the Plateau borders are only narrow channels for small~$\varepsilon$, their lengths are not well defined and the volume and dissipation effect in the nodes varies.
 \cite{Koehler2000} presented a relationship between~$\varepsilon$, $L$ and $r$, valid for at least $\varepsilon$ up to 0.1.
They derived a generalized foam drainage equation which covers the two limiting cases of channel- and node-dominated models, respectively.
To remove the variable~$L$, \cite{Neethling2002foam} made the common assumption that for small~$\varepsilon$, bubbles can be assumed to have the form of a tetrakaidecahedron (Kelvin cell) and used the equation $4\pi r_\mathrm{b}^3/{3}=(1-\varepsilon)2^{2/7}L^3$, where~$r_\mathrm{b}$ is the bubble radius.
Thereby, they obtained the algebraic equation
\begin{equation}\label{eq:Neeth}
\varepsilon= 0.3316\left(
\frac{r}{r_\mathrm{b}}\right)^2(1-\varepsilon)^{2/3}
+ 0.5402\left(
\frac{r}{r_\mathrm{b}}\right)^3(1-\varepsilon),
\end{equation}
which is implicit in all its variables.
Containing these three variables, they derived a PDE valid for $0\leq\varepsilon\leq\varepsilon_\mathrm{c}$ by considering dissipation both from the Plateau borders and the nodes.
Assuming~$r_\mathrm{b}$ is constant, their PDE and algebraic equation defines the unknowns~$\varepsilon$ and~$r$.

\cite{Stevenson2006} demonstrated that the effective relative fluid-gas velocity~$u$  could be very well approximated by a power law of the type~\eqref{eq:power}, at least for fluid volume fraction up to $\varepsilon\approx 0.2$.
In particular, \cite{Stevenson2004} approximated Equation~\eqref{eq:Neeth} by 
\begin{equation*}
\frac{r}{r_\mathrm{b}} = m\varepsilon^{n_\mathrm{S}},
\quad\text{with\quad $m=1.28$, $n_\mathrm{S}=0.46$.} 
\end{equation*}
This equation can be substituted into~\eqref{eq:uPBaver} to give
\begin{equation}\label{eq:drainage}
u = -v_\mathrm{drain}\varepsilon^{2n_\mathrm{S}}\bigl(
1+d_\mathrm{cap}\varepsilon^{-(1+n_\mathrm{S})}
 \partial_z \varepsilon \bigr) \quad \text{for 
$0\leq\varepsilon<\varepsilon_\mathrm{c}$,} 
\end{equation}
where the drainage velocity~$v_\mathrm{drain}$ (with respect to gravity and dissipation) and the dimensionless capillarity-to-gravity parameter~$d_\mathrm{cap}$ are given by
\begin{equation*}
v_\mathrm{drain}
:= \frac{m^2 C^2r_\mathrm{b}^2\rho_\mathrm{f}g}{3C_\mathrm{PB}\mu},
\qquad
d_\mathrm{cap}:=\frac{n_\mathrm{S}\gamma_\mathrm{w}}{mr_\mathrm{b}\rho_\mathrm{f}g}.
\end{equation*}
The derivative term in~\eqref{eq:drainage} models the capillarity that is not present for separated bubbles; see~\eqref{eq:power}.
Hence, we suggest the relative fluid-gas velocity
\begin{equation*}
u:=-\begin{cases}
v_\mathrm{drain}\varepsilon^{2n_\mathrm{S}} \bigl(
1+d_\mathrm{cap}\varepsilon^{-(1+n_\mathrm{S})}
 \partial_z \varepsilon \bigr) 
&\text{for $0\leq\varepsilon<\varepsilon_\mathrm{c}$,} \\
v_{\mathrm{term}}\varepsilon^{n_{\mathrm{b}} -1} 
&\text{for $\varepsilon_\mathrm{c}\leq\varepsilon\leq 1$} 
\end{cases}
\end{equation*}
with the compatibility condition (continuity across $\varepsilon = \varepsilon_\mathrm{c}$)
\begin{equation}\label{eq:compat}
v_\mathrm{drain}\varepsilon_\mathrm{c}^{2n_\mathrm{S}}
= v_{\mathrm{term}}\varepsilon_\mathrm{c}^{n_{\mathrm{b}} -1}
\quad\Leftrightarrow\quad
\frac{v_\mathrm{drain}}{v_{\mathrm{term}}}
=\varepsilon_\mathrm{c}^{n-1-2n_\mathrm{S}}.
\end{equation}
Values for~$n_{\mathrm{b}}$ in the literature range from~2 to~3.2  \citep{Dickinson2014, Galvin2014II, Pal1989, Vandenberghe2005}.

Recalling once again that $v_\mathrm{f}=(1-\varepsilon)u$, we now define the velocity  function 
\begin{align*}
\tilde{v}_\mathrm{f}(\varepsilon) &:= 
\begin{cases}
v_\mathrm{drain}(1-\varepsilon)\varepsilon^{2n_\mathrm{S}}
& \text{for $0\leq\varepsilon<\varepsilon_\mathrm{c}$,} \\
v_{\mathrm{term}}(1-\varepsilon)\varepsilon^{n_{\mathrm{b}} -1} 
&\text{for $\varepsilon_\mathrm{c}\leq\varepsilon\leq 1$} 
\end{cases}
\end{align*}
and the diffusion function 
\begin{align*}
d_\mathrm{f}(\varepsilon) :=\begin{cases}
v_\mathrm{drain}d_\mathrm{cap}(1-\varepsilon)\varepsilon^{n_\mathrm{S}} 
&\text{for $0\leq\varepsilon<\varepsilon_\mathrm{c}$,} \\
0
&\text{for $\varepsilon_\mathrm{c}\leq\varepsilon\leq 1$,} 
\end{cases}
\end{align*}
so that the liquid flux (in a closed vessel) becomes  
\begin{equation} \label{lflux} 
\varepsilon v_\mathrm{f} = \varepsilon(1-\varepsilon)u =-\varepsilon\tilde{v}_\mathrm{f}(\varepsilon) -d_\mathrm{f}(\varepsilon) \partial_z \varepsilon 
=-\varepsilon\tilde{v}_\mathrm{f}(\varepsilon)
- \partial_z \tilde{D}(\varepsilon),
\end{equation}
where the integrated diffusion function is
\begin{equation*}
\tilde{D}(\varepsilon)
:=\int_{0}^{\varepsilon}d_\mathrm{f}(\xi)\,\mathrm{d}\xi
=\begin{cases}
v_\mathrm{drain}d_\mathrm{cap}
\biggl(\dfrac{\varepsilon^{n_\mathrm{S}+1}}{n_\mathrm{S}+1}
- \dfrac{\varepsilon^{n_\mathrm{S}+2}}{n_\mathrm{S}+2}
\biggr) 
&\text{for $0\leq\varepsilon<\varepsilon_\mathrm{c}$,} \\[10pt]
v_\mathrm{drain}d_\mathrm{cap}
\biggl(\dfrac{\varepsilon_\mathrm{c}^{n_\mathrm{S}+1}}{n_\mathrm{S}+1}
- \dfrac{\varepsilon_\mathrm{c}^{n_\mathrm{S}+2}}{n_\mathrm{S}+2}
\biggr) 
& \text{for $\varepsilon_\mathrm{c}\leq\varepsilon\leq 1$.} 
\end{cases}
\end{equation*}
(Notice that $\tilde{D}( \varepsilon)$ is constant, and therefore $\tilde{D}' ( \varepsilon) =0$, 
 for $\varepsilon_{\mathrm{c}} \leq \varepsilon \leq 1$.)
Inserting the expression \eqref{lflux} into the conservation law for the fluid phase $\partial_t\varepsilon+\partial_z(\varepsilon v_\mathrm{f})=0$, we obtain the generalized Equation~\eqref{eq:gen_drainage} modelling both rising bubbles and drainage of froth in a closed container.

\section{A model of flotation including froth drainage} \label{sec:model}

\subsection{Assumption on the tank and mixture}

We use a one-dimensional setup of the Reflux Flotation Cell by~\cite{Dickinson2014}; see Figure~\ref{fig:Column}.
A mixture of slurry and aggregates is fed at the height $z=z_\mathrm{F}$ at the volumetric flow~$Q_\mathrm{F}>0$ and wash water is injected at the top effluent level $z=z_\mathrm{E}$ at $Q_\mathrm{W}\geq 0$.
At $z=z_\mathrm{U}$, a volumetric flow~$Q_\mathrm{U}\geq 0$ is taken out.
In the one-dimensional model on the real line, there are four zones, two inside the vessel plus the underflow and effluent zones.
The resulting effluent volumetric overflow~$Q_{\mathrm{E}}:=Q_{\mathrm{W}} + Q_{\mathrm{F}} - Q_{\mathrm{U}}$ is assumed to be positive so that the mixture is conserved and the vessel is always completely filled.
In comparison to the  previous treatments \citep{SDIMAflot2019, SDmineng3_flot_sed}, here we do not separate the wash water inlet and the effluent level, i.e., $z_\mathrm{W}=z_\mathrm{E}$.
The cross-sectional area is assumed to satisfy
\begin{align*}
A(z)=\begin{cases}
A_\mathrm{E}&\text{for $z \geq z_{\mathrm{F}}$,} \\
A_\mathrm{U} &\text{for $z < z_{\mathrm{F}}$.}
\end{cases} 
\end{align*}

Particles trapped in the froth region influence the drainage of fluid \citep{Ata2012,Haffner2015}, but for simplicity 
 we   nevertheless  assume  that the volume fraction of aggregates (bubbles with attached hydrophobic particles) can be determined as a function of height and time by a single equation.
Thus,  the suspension in the interstices outside the bubbles is assumed to behave independently  of the volume fraction of (hydrophilic) particles.
Such particles may however settle within the suspension, which undergoes bulk transport.
 From now on we denote by~$\phi=1-\varepsilon$  the volume fraction of aggregates.
As a first approximation in a closed vessel, $\phi$~can be obtained by solving~\eqref{eq:gen_drainage} for~$\varepsilon$ and setting $\phi=1-\varepsilon$, but we proceed to derive an explicit equation for~$\phi$  since that will be extended to the more complicated model of a  flotation column with in- and outlets.

\subsection{Equation for aggregates with froth drainage in a closed tank}

\begin{figure} 
\centering 
\includegraphics[width=0.45\textwidth]{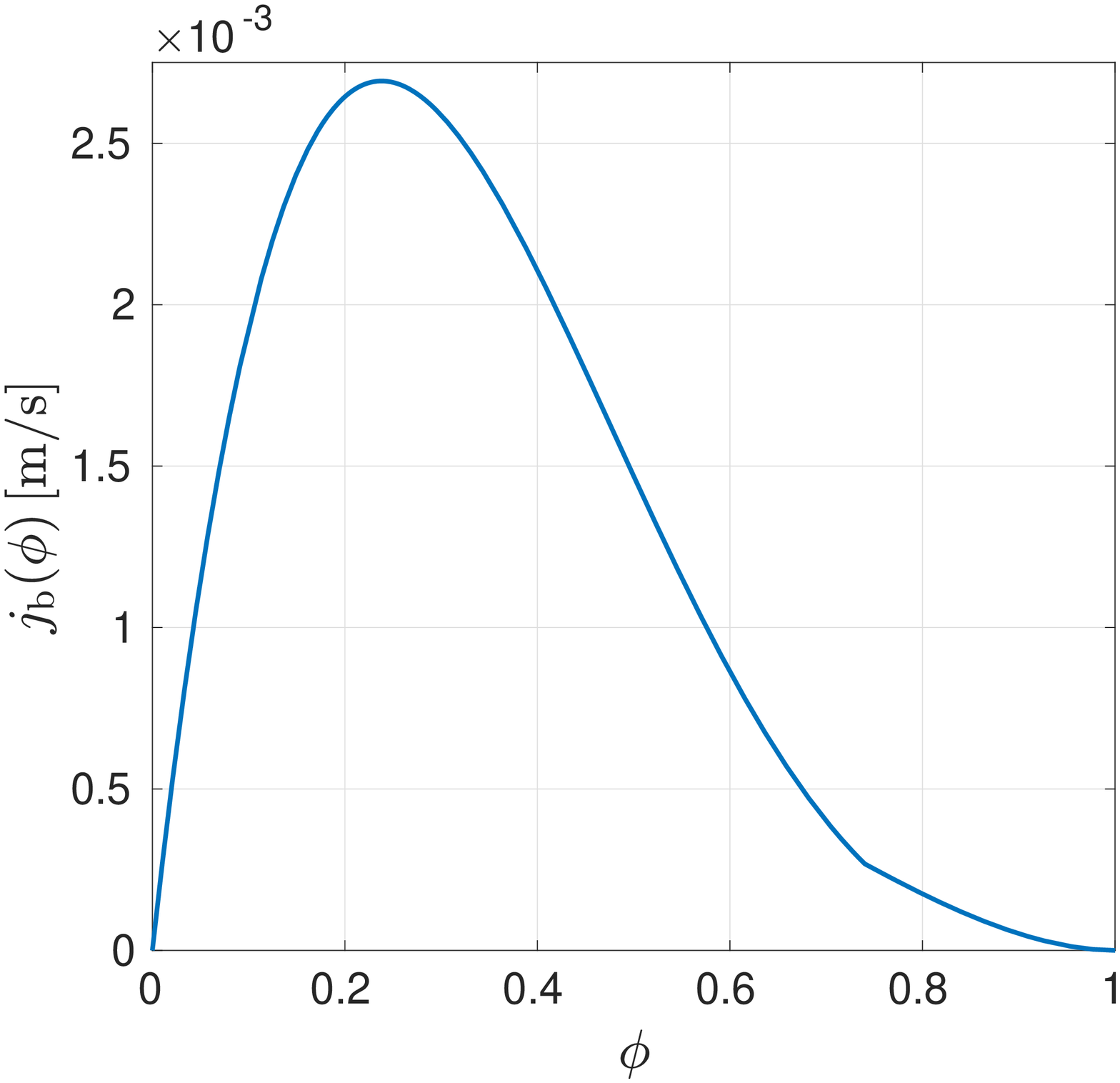}  \includegraphics[width=0.45\textwidth]{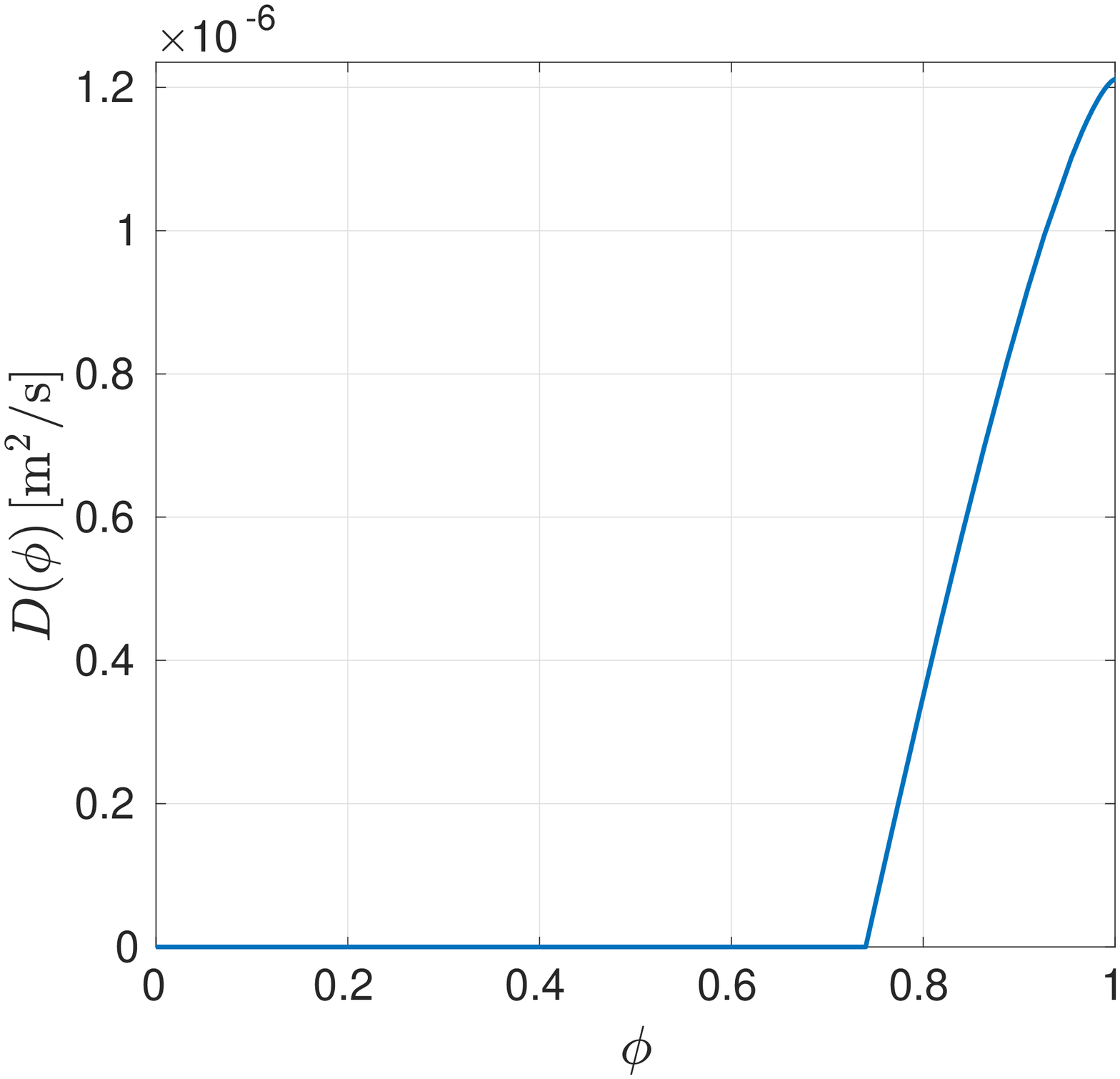} 
\caption{Left: function $j_\mathrm{b}(\phi)=\phi\tilde{v}(\phi)$. Right: diffusion function $D(\phi)$ modelling capillarity.
Note the behaviour of these functions at the critical concentration $\phi_\mathrm{c}=0.74$.} \label{fig:graphs}
\end{figure}

We recall the gas-phase velocity~$v=-(1-\phi)u$ and the compatibility condition~\eqref{eq:compat}, and define
\begin{align}
\tilde{v}(\phi) &:=
\begin{cases}
v_{\mathrm{term}}(1-\phi)^{n_{\mathrm{b}}} 
&\text{for $0\leq\phi\leq\phi_\mathrm{c}$,} \\
v_\mathrm{drain}(1-\phi)^{2n_\mathrm{S}+1}
=v_{\mathrm{term}}
\dfrac{(1-\phi)^{2n_\mathrm{S}+1}}{(1-\phi_\mathrm{c})^{2n_\mathrm{S}+1-n_{\mathrm{b}}}} 
&\text{for $\phi_\mathrm{c}<\phi\leq 1$},
\end{cases}\label{eq:vtilde}\\
{d}(\phi) &:=\begin{cases}
0 
& \text{for $0\leq\phi\leq\phi_\mathrm{c}$,} \\
v_\mathrm{drain}d_\mathrm{cap}\phi(1-\phi)^{n_\mathrm{S}}
=v_{\mathrm{term}}d_\mathrm{cap}
\dfrac{\phi(1-\phi)^{n_\mathrm{S}}}{(1-\phi_\mathrm{c})^{2n_\mathrm{S}+1-n_{\mathrm{b}}}}
& \text{for $\phi_\mathrm{c}<\phi\leq 1$.} 
\end{cases}  \label{dphidef} 
\end{align}
With the batch-drift flux function~$j_\mathrm{b}(\phi):=\phi\tilde{v}(\phi)$, where $\tilde{v}(\phi)$ is given by~\eqref{eq:vtilde},  we can write the aggregate-phase flux (in a closed container) as
\begin{equation*}
\phi v = -\phi(1-\phi)u 
= \phi\tilde{v}(\phi) + {d}(\phi) \partial_z (1-\phi) 
= \phi\tilde{v}(\phi) - {d}(\phi) \partial_z \phi 
= j_\mathrm{b}(\phi) - \partial_z {D}(\phi),
\end{equation*}
where $D(\phi)$  is defined by \eqref{Dphidef}. In light of \eqref{dphidef} we obtain 
\begin{equation}\label{eq:D}
{D}(\phi) 
=\begin{cases}
0 
&\text{for $0\leq\phi\leq\phi_\mathrm{c}$,} \\
v_\mathrm{drain}d_\mathrm{cap}\dfrac{\omega(\phi_\mathrm{c})-\omega(\phi)}
{(n_\mathrm{S}+1)(n_\mathrm{S}+2)} 
&\text{for $\phi_\mathrm{c}<\phi\leq 1$,}   
\end{cases}
\end{equation}
where 
$\omega(\phi):=(1-\phi)^{n_\mathrm{S}+1}((n_\mathrm{S}+1)\phi+1)$ and we reconfirm the 
 property \eqref{by3.2}.  
The conservation law $\partial_t\phi+\partial_z(\phi v)=0$ now yields the following equation
 for the volume fraction $\phi= \phi(z,t) \in [0,1]$ of aggregates in a closed vessel: 
\begin{equation}\label{eq:batch_phi}
\partial_t \phi 
+ \partial_z  j_\mathrm{b}(\phi) 
= \partial_z^2{D}(\phi). 
\end{equation}
The graphs of the constitutive functions $j_\mathrm{b}(\phi)$ and $D(\phi)$ are drawn in Figure~\ref{fig:graphs}.

\subsection{Three phases and constitutive assumptions} \label{subsec:threephase} 

The three phases and their volume fractions are the fluid~$\phi_\mathrm{f}$,  the solids~$\psi$, and the aggregates~$\phi$, 
 where $\phi_\mathrm{f} + \psi + \phi=1$. 
By suspension we mean the fluid and solid phases.
The volume fraction of solids {\em within the suspension}  $\varphi$ is defined by
\begin{equation*}
\varphi :=\frac{\psi}{\psi + \phi_\mathrm{f}}=\frac{\psi}{1-\phi}.
\end{equation*}
The drift-flux and solids-flux theories utilize  constitutive functions for the aggregate upward batch flux~$j_\mathrm{b}(\phi)$ and the solids batch sedimentation flux~$f_\mathrm{b}(\varphi):=\varphi v_\mathrm{hs}(\varphi)$, where $v_\mathrm{hs}(\varphi)$ is the hindered-settling   function.  For simplicity, we employ the common  expression~\citep{Richardson1954}
\begin{equation}\label{eq:vhs}
v_\mathrm{hs}(\varphi) =v_\infty(1-\varphi)^{n_\mathrm{RZ}}, \quad \text{where $n_\mathrm{RZ}>1$.} 
\end{equation}

Applying the conservation of mass to each of the three phases, introducing the volume-average velocity, or bulk velocity, of the mixture~$q$ and the relative velocities of both the aggregate-suspension and the solid-fluid, \cite{SDIMAflot2019} derived the PDE model~\eqref{eq:gov} without the capillarity function~$D(\phi)$.
In particular, the volumetric flows in and out of the flotation column define explicitly
\begin{equation}\label{eq:q}
q(z,t):=\begin{cases}
q_\mathrm{E}:=(-{Q_\mathrm{U}} +{Q_\mathrm{F}} +{Q_\mathrm{W}})/{A_\mathrm{E}} & \text{for $z\geq z_\mathrm{E}$}, \\
q_2:=(-{Q_\mathrm{U}} +{Q_\mathrm{F}})/{A_\mathrm{E}} & \text{for $z_\mathrm{F}\leq z<z_\mathrm{E}$,} \\
q_1=q_\mathrm{U}
:=-{Q_\mathrm{U}}/{A_\mathrm{U}} & \text{for $z<z_\mathrm{F}$}.
\end{cases}
\end{equation}
In the underflow  and effluent zones all phases are assumed to have the same velocity, i.e., they follow the bulk flow. Then the total convective fluxes for~$\phi$ and~$\varphi$ are given by
\begin{align*}
J(\phi,z,t)& =\begin{cases}
j_\mathrm{E}(\phi ,t):= q_\mathrm{E}(t)\phi & \text{for $z\geq z_\mathrm{E}$,} \\
j_2(\phi,t):=q_2(t)\phi +j_\mathrm{b}(\phi) & \text{for $z_\mathrm{F}\leq z<z_\mathrm{E}$,} \\
j_1(\phi,t):=q_1(t)\phi +j_\mathrm{b}(\phi) & \text{for $z_\mathrm{U}\leq z<z_\mathrm{F}$,} \\
j_\mathrm{U}(\phi,t):=q_1(t)\phi & \text{for $z<z_\mathrm{U}$,} 
\end{cases}\\
F(\varphi,\phi,z,t)& =\begin{cases}
f_\mathrm{E}(\varphi ,\phi ,t):=-(1-\phi) q_\mathrm{E}(t)\varphi & \text{for $z\geq z_\mathrm{E}$,}  \\
f_2(\varphi ,\phi,t) & \text{for $z_\mathrm{F}\leq z<z_\mathrm{E}$,} \\
f_1(\varphi ,\phi,t) & \text{for $z_\mathrm{U}\leq z<z_\mathrm{F}$,} \\
f_\mathrm{U}(\varphi,\phi ,t):=-(1-\phi )q_1(t)\varphi & \text{for $z<z_\mathrm{U}$} 
\end{cases}
\end{align*}
with the zone-settling flux functions (positive in the direction of sedimentation (decreasing $z$))  
\begin{align*}
f_k(\varphi,\phi,t):= &\,(1-\phi)f_\mathrm{b}(\varphi) + \big(j_\mathrm{b}(\phi) - (1-\phi)q_k(t)\big)\varphi\notag\\
=&\, (1-\phi )f_\mathrm{b}(\varphi ) + \big(j_k(\phi,t) - q_k(t)\big)\varphi,\quad k=1,2.
\end{align*}
With the capillarity function~$D(\phi)$, the batch flux~$j_\mathrm{b}(\phi)$ is extended to~$j_\mathrm{b}(\phi)-\partial_z D(\phi)$; cf.~\eqref{eq:batch_phi}.
Hence,  the total flux of the aggregates for any $z\in\mathbb{R}$ is
\begin{equation*}
\Phi(\phi,\partial_z\phi,z,t) := J(\phi,z,t)-\gamma(z) \partial_z  D(\phi),
\end{equation*}
where the characteristic function is
\begin{equation*}
\gamma(z) := 
\begin{cases}
1 & \text{for $z\in \left[ z_\mathrm{U},z_\mathrm{E}\right)$,}  \\
0 & \text{for $z\notin \left[ z_\mathrm{U},z_\mathrm{E}\right)$,} 
\end{cases}
\end{equation*}
and the total flux of the solids in the $z$-direction is ($F$ and $\tilde{F}$ are positive in the downwards direction of sedimentation, which is opposite to the $z$-direction)
\begin{equation}\label{eq:Ftotalflux}
\Psi(\psi,\partial_z\psi,\phi,z,t):= -\tilde{F}(\psi,\phi,z,t)
+\gamma(z)\frac{\psi}{1-\phi} \partial_z D(\phi),
\end{equation}
where
\begin{equation*}
\tilde{F}(\psi,\phi,z,t) := 
\begin{cases}
F\left(\dfrac{\psi}{1-\phi},\phi,z,t\right) & \text{if $0\leq\phi<1$,} \\
0 & \text{if $\phi=1$.} 
\end{cases}
\end{equation*}

The conservation law applied on the two phases with the total fluxes~$\Phi$ and $\Psi$ yields the governing system of equations~\eqref{eq:gov} in the case capillarity are included.
That system defines  solutions on the real line and next we define the outlet concentrations of the flotation column.

\subsection{Outlet concentrations}

Given the PDE solutions $\phi=\phi(z,t)$ and $\varphi=\varphi(z,t)$ of~\eqref{eq:gov}, we define the boundary concentrations at each in- or outlet by~$\phi_\mathrm{U}^{\pm}=\phi_\mathrm{U}^{\pm}(t):=\phi(z_\mathrm{U}^{\pm},t)$, etc.
Conservation of mass across $z=z_\mathrm{U}$ yields 
\begin{align}
j_1(\phi_\mathrm{U}^+,t) - \left. \partial_z  D(\phi) \right|_{z=z_\mathrm{U}^+}
&=j_\mathrm{U}(\phi_\mathrm{U}^-,t),
\label{eq:phiU}\\
f_1(\varphi_\mathrm{U}^+,t)
-\varphi_\mathrm{U}^+\left. \partial_z  D(\phi) \right|_{z=z_\mathrm{U}^+}
&=f_\mathrm{U}(\varphi_\mathrm{U}^-,t). 
\label{eq:varphiU}
\end{align}
The underflow concentrations of the flotation column are defined by~$\phi_\mathrm{U}(t):=\phi_\mathrm{U}^-(t)$ and $\varphi_\mathrm{U}(t):=\varphi_\mathrm{U}^-(t)$.
These concentrations can in fact be obtained from the solution inside the column ($z_\mathrm{U}<z<z_\mathrm{E}$) from~\eqref{eq:phiU} and \eqref{eq:varphiU} together with a uniqueness condition; see~\cite{SDjhde1}.

For the effluent level $z=z_\mathrm{E}$, the analogous situation holds:
\begin{align}
j_2(\phi_\mathrm{E}^-,t) 
- \left. \partial_z  D(\phi) \right|_{z=z_\mathrm{E}^-}
&= j_\mathrm{E}(\phi_\mathrm{E}^+,t),\label{eq:BCgasE}\\
f_2(\varphi_\mathrm{E}^-,\phi_\mathrm{E}^-,t) 
-\varphi_\mathrm{E}^- \left. \partial_z  D(\phi) \right|_{z=z_\mathrm{E}^-}
&= f_\mathrm{E}(\varphi_\mathrm{E}^+,\phi_\mathrm{E}^+,t),
\label{eq:BCsolidE}
\end{align}
In the one-dimensional PDE model~\eqref{eq:gov} without boundary conditions, the solution $\phi=\phi(z,t)$ (analogously for $\varphi$) in the interval $z>z_\mathrm{E}$ is governed by the linear transport PDE $\partial_t\phi+ (Q_\mathrm{E}/A_\mathrm{E})\partial_z\phi=0$ and the boundary value~$\phi_\mathrm{E}^+(t)$.
The effluent outlet concentrations are defined by $\phi_\mathrm{E}:=\phi_\mathrm{E}^+$ and $\varphi_\mathrm{E}:=\varphi_\mathrm{E}^+$.
In the concluding section, we discuss how bursting bubbles at the top can be incorporated in the model.

\section{Steady-state analysis}\label{sec:steadystates}

\subsection{Definition of a desired steady state}

In the case of no capillarity, \cite{SDIMAflot2019} provided  detailed constructions of all steady states, and  \cite{SDmineng3_flot_sed, SDwst_DAF_CFF} sorted out the most interesting steady states for the applications and how to control these by letting the volumetric flows satisfy certain nonlinear inequalities, which can be visualized in so-called operating charts.
We assume that $Q_\mathrm{F}$, $\phi_\mathrm{F}$, and$\psi_\mathrm{F}$ are given variables and that $Q_\mathrm{U}$ and $Q_\mathrm{W}$ are control variables.
The purpose here is to provide  an improved model of the froth region and we therefore focus on the steady states when a layer of froth in zone~2 is possible.
We consider only solutions where the froth layer does not fill the entire zone~2, so that there is at least a small region above the feed inlet with aggregate volume fraction below the critical one.
As mentioned before, it is  assumed  that the wash water is sprinkled at the top of the column, which is commonly done and gives fewer steady states to analyse.
A \emph{desired steady state} is  defined to be a stationary solution that has
\begin{equation}\label{eq:desired}
\begin{alignedat}2
&\text{no aggregates below the feed level}&\quad&\Rightarrow\quad\phi_\mathrm{U}=0,\\
&\text{no solids above the feed level}&\quad&\Rightarrow\quad\varphi_\mathrm{E}=0,\\
&\text{a froth layer  that does not fill the entire zone~2}&\quad&\Rightarrow\quad\phi(z_\mathrm{F}^+)<\phi_\mathrm{c}.
\end{alignedat}
\end{equation}
The reversed implications do not hold in the two first statements for the following reasons.
Since the bulk flow in zone~1 is directed downwards, there exist steady-state solutions with a standing layer of aggregates below the feed level.
Analogously, if the bulk flow in zone~2 is directed upwards, there may be a layer of standing solids when their settling velocity is balanced by the upward bulk velocity; see \cite{SDIMAflot2019}.

\subsection{Properties of the batch-flux density functions}

\begin{figure}[t] 
\centering\includegraphics[width=0.5\textwidth]{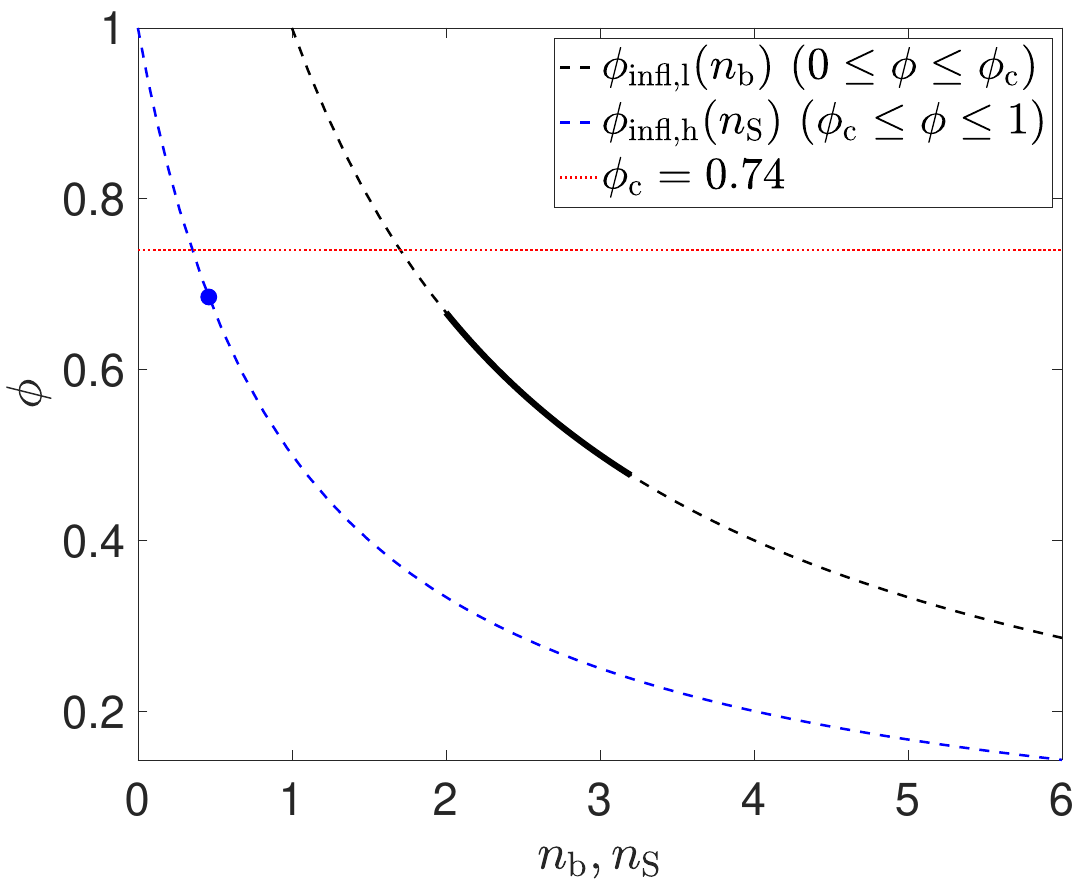}
\caption{Evolution of the inflection points of $j_\mathrm{bl}$ and $j_\mathrm{bh}$.
The literature values $2\leq n_{\mathrm{b}} \leq 3.2$ give an interval (solid black) of possible $\phi_\mathrm{infl,l}$ that lie entirely below~$\phi_\mathrm{c}=0.74$ (red line).
With $n_\mathrm{S}=0.46$, the inflection point (blue dot)~$\phi_\mathrm{infl,h}=1/(n_\mathrm{S}+1)\approx 0.685<\phi_\mathrm{c}$; hence, $j_\mathrm{bh}$ is strictly convex for $\phi\geq\phi_\mathrm{c}$.
\label{fig:inflection}}
\end{figure}

\begin{figure}[t] 
\centering\includegraphics[width=0.9\textwidth]{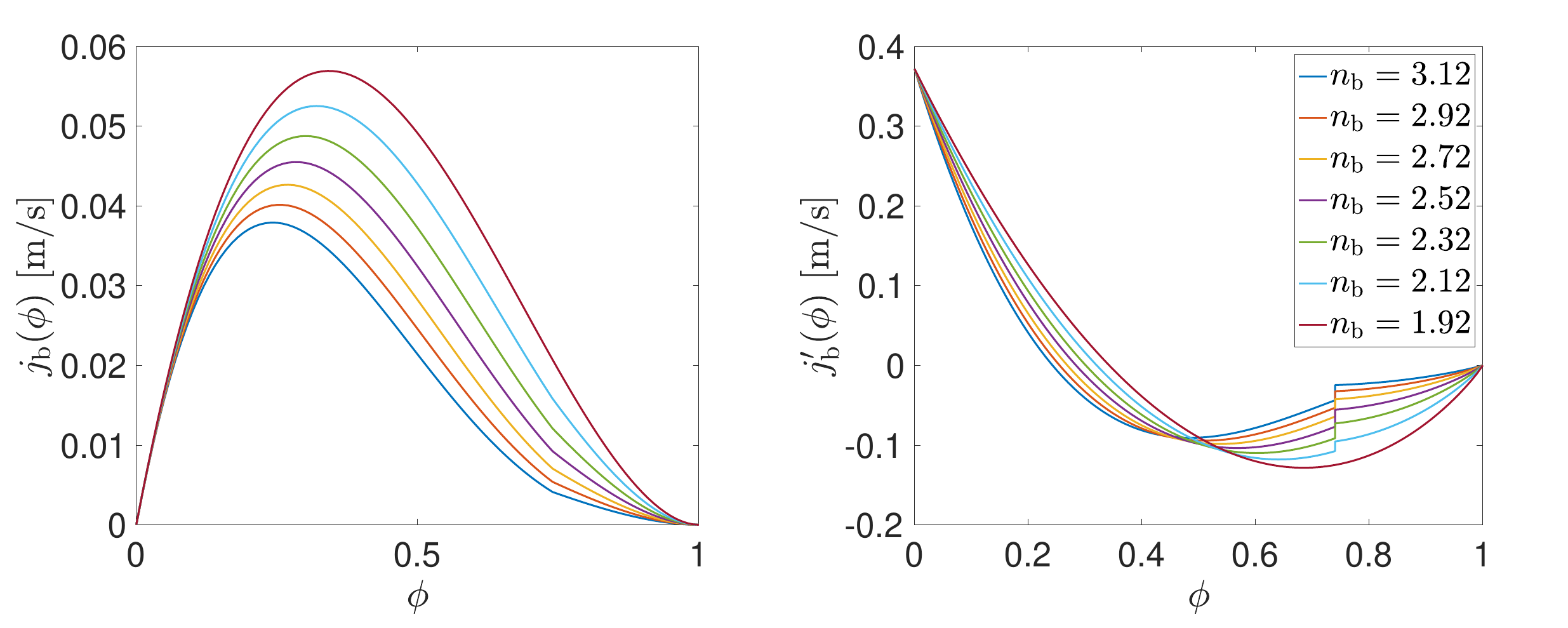}
\caption{Plots of $j_\mathrm{b}(\phi)$ (left) and $j'_\mathrm{b}(\phi)$ (right) for $n_\mathrm{S}=0.46$, $v_\mathrm{term}=0.3718$ and  various  values of $n_{\mathrm{b}}$ that satisfy~\eqref{eq:ncond}.
\label{fig:evol_jb2}}
\end{figure}

With $\tilde{v}$ given by~\eqref{eq:vtilde}, the continuous batch-drift flux function is
\begin{equation*}
j_\mathrm{b}(\phi) := \phi \tilde{v}(\phi) =
\begin{cases}
j_\mathrm{bl}(\phi):=\phi v_{\mathrm{term}}(1-\phi)^{n_{\mathrm{b}} } 
& \text{for $0\leq\phi\leq\phi_\mathrm{c}$,} \\
j_\mathrm{bh}(\phi):=\phi v_{\mathrm{term}} \dfrac{(1-\phi)^{2n_\mathrm{S}+1}}{(1-\phi_\mathrm{c})^{2n_\mathrm{S}+1-n_{\mathrm{b}}}} 
& \text{for $\phi_\mathrm{c}<\phi\leq 1$,} 
\end{cases}
\end{equation*}
where we have introduced the low $j_\mathrm{bl}$ and high $j_\mathrm{bh}$ parts of it.
Any function $u \mapsto u(1-u)^n$ has a unique inflection point at $u_\mathrm{infl}=2/(n+1)$.
Figure~\ref{fig:inflection} shows the inflection points
\begin{equation*}
\phi_\mathrm{infl,l}(n_{\mathrm{b}})=\frac{2}{n_{\mathrm{b}}+1},\qquad
\phi_\mathrm{infl,h}(n_\mathrm{S})=\frac{2}{2n_\mathrm{S}+1+1} =\frac{1}{n_\mathrm{S}+1} 
\end{equation*}
of $j_\mathrm{bl}$ and $j_\mathrm{bh}$, as functions of the exponents~$n$ and $n_\mathrm{S}$, respectively.
With the values~$\phi_\mathrm{c}=0.74$ and $n_\mathrm{S}=0.46$ suggested in the literature (see Section~\ref{sec:drainage}),  
and the interval $2\leq n_{\mathrm{b}} \leq 3.2$,  there is only one  inflection point of $j_\mathrm{b}$ in $0\leq\phi\leq 1$ and this lies below~$\phi_\mathrm{c}$; see Figure~\ref{fig:evol_jb2}, which also shows that there may be a jump in the derivative of $j_\mathrm{b}$ at $\phi=\phi_\mathrm{c}$.
Since
\begin{align*}
j_\mathrm{b}'(\phi) &=
\begin{cases}
v_{\mathrm{term}}(1-\phi)^{n_{\mathrm{b}} -1}(1-(1+n_{\mathrm{b}})\phi) 
& \text{for $0<\phi<\phi_\mathrm{c}$,} \\
v_{\mathrm{term}}\dfrac{(1-\phi)^{2n_\mathrm{S}}(1-(2+2n_\mathrm{S})\phi)}{(1-\phi_\mathrm{c})^{2n_\mathrm{S}+1-n_{\mathrm{b}}}} 
&\text{for $\phi_\mathrm{c}<\phi< 1$,} 
\end{cases}
\end{align*}
we get that 
\begin{equation}\label{eq:ncond}
j_\mathrm{b}'(\phi_\mathrm{c}^-)\leq j_\mathrm{b}'(\phi_\mathrm{c}^+)
\quad\Leftrightarrow\quad
n_{\mathrm{b}}\geq1+2n_\mathrm{S}\approx 1.92.
\end{equation}
When this is satisfied, the exponent in the compatibility condition~\eqref{eq:compat} is nonnegative and the entire $j_\mathrm{b}$ has only one inflection point $\phi_{\mathrm{infl}}=\phi_{\mathrm{infl,l}}\in(0,\phi_\mathrm{c})$.

\subsection{Properties of the zone flux functions}\label{sec:jprop}

The zone flux functions~$j_k$, $f_k(\cdot,\phi)$, $k=1,2$, have an additional linear term due to the bulk velocity of the zone.
Let $j(\phi)=j_\mathrm{b}(\phi)+q\phi$ denote a general zone flux function, where we drop the $t$-variable when considering steady states.
We will sometimes write out the dependence on~$q$; $j(\phi;q)$.
The inflection point~$\phi_{\rm infl}$ of~$j$ is independent of~$q$, however, the local maximum~$\phi^\mathrm{M}=\phi^\mathrm{M}(q)<\phi_\mathrm{c}$ depends on~$q$.
To  provide  an explicit definition, we first define
\begin{equation*}
q_\mathrm{neg} :=-j_\mathrm{b}'(0),\qquad
\qbbar :=-j_\mathrm{b}'(\phi_\mathrm{infl}).
\end{equation*}
For~$q\leq q_\mathrm{neg}$, $j(\cdot,q)$ is decreasing and for $q\geq\qbbar$, $j(\cdot,q)$ is increasing.
For intermediate values of~$q$, the local maximum exists and satisfies $0=j'(\phi^\mathrm{M})=j_\mathrm{b}'(\phi^\mathrm{M})+q$.
Since the restriction $\smash{(j_\mathrm{b}|_{(0,\phi_\mathrm{infl})})'}$ is a strictly decreasing function, we can define
\begin{equation*}
\phi^{\mathrm{M}}=\phi^{\mathrm{M}}(q):=
\begin{cases}
0  & \text{if $q\leq q_\mathrm{neg}$,} \\
((\left.j_\mathrm{b}\right|_{(0,\phi_\mathrm{infl})})')^{-1}(-q) & \text{if $q_\mathrm{neg}<q<\qbbar$,} \\
\phi_\mathrm{infl} & \text{if $q\ge\qbbar$.}
\end{cases}
\end{equation*}
For $q_\mathrm{neg}<q<0$, there is a zero of $j(\cdot;q)$ which we denote by $\phi_\mathrm{Z}=\phi_\mathrm{Z}(q)\in(0,1)$.
For a specific zone flux functions $j_k$, we use the notation $\phi_k^{\mathrm{M}}=\phi^{\mathrm{M}}(q_k)$ and $\phi_{k\mathrm{Z}}=\phi_{\mathrm{M}}(q_k)$.

In a similar way, one can define the local minimum point, greater than the inflection point, for $0\leq q<\qbbar$. We denote it by $\phi_{k \mathrm{M}}=\phi_{\mathrm{M}}(q_k)$.
For $q\geq\qbbar$, we define $\phi_{k \mathrm{M}}(q):=\phi_\mathrm{infl}$. Furthermore, for a given~$\phi_{k \mathrm{M}}$, we define~$\phi_{k \mathrm{m}}$ as the unique value that satisfies
\begin{equation}\label{eq:varphi1m}
j_k(\phi_{k\mathrm{m}};q)=j_k(\phi_{k\mathrm{M}};q),\qquad
0\leq\phi_{k \mathrm{m}}\leq\phi_{\mathrm{infl}}.
\end{equation}
Analogous definitions can be made for the flux functions $f_k(\cdot,\phi,t)$, $k=1,2$.
\subsection{Construction of steady states}\label{sec:constructionSS}

We seek piecewise smooth and piecewise monotone steady-state solutions~$\phi=\phi(z)$ of~\eqref{eq:gov}.
Such solutions may contain jump discontinuities within or between the zones.
In the case $D\equiv 0$,  \cite{SDIMAflot2019} outlined the details on how to 
 construct unique steady-state solutions and we will not go through the entire machinery here.
The basic idea is to glue together solutions within each zone  in a unique way so that the conservation of mass holds across the zone borders.
Two such so-called Rankine-Hugoniot conditions (jump conditions) are~\eqref{eq:BCgasE} and \eqref{eq:BCsolidE}.
Since each such equation has two unknowns; for example, $\phi_\mathrm{E}^-$ and $\phi_\mathrm{E}^+$ in~\eqref{eq:BCgasE}, another so-called entropy condition in the theory of degenerate parabolic PDEs with spatially discontinuous coefficients is needed to establish a unique pair of boundary values~\citep{SDjhde1}.
Furthermore, as the values $\phi_\mathrm{E}^-$ and $\phi_\mathrm{E}^+$ are obtained, these are substituted into~\eqref{eq:BCsolidE} and a  similar procedure yields  $\varphi_\mathrm{E}^-$ and $\varphi_\mathrm{E}^+$.

The new ingredient due to the drainage is the term $\partial_z D(\phi)|_{z=z_\mathrm{E}^-}$ in~\eqref{eq:BCgasE} and \eqref{eq:BCsolidE}.
The property~\eqref{by3.2} implies the following (see \cite{Evje2000,SDjhde1} for further details).
A discontinuity of the solution~$\phi(\cdot,t)$, within or between zones, is possible only between two values in the interval $0\leq\phi\leq\phi_\mathrm{c}$.
Furthermore, since we are seeking piecewise smooth and piecewise monotone steady-state solutions, the fact $d(\phi)>0$ implies that if one of the values of the discontinuity is $\phi_\mathrm{c}$, this must be the larger value and located on the right; i.e., the left value of the jump $\phi^-<\phi_\mathrm{c}$.
Furthermore, $j_2(\phi)\geq j_2(\phi_\mathrm{c})$ for $\phi^-\leq\phi\leq\phi_\mathrm{c}$, and in a right neighbourhood of the jump, $\phi'(z)\geq 0$.

With these facts in mind, we  now construct steady-state solutions.
Let $H(z)$ denote the Heaviside function and assume that all volumetric flows and feed volume fractions are time independent.
A stationary solution~$\phi=\phi(z)$ of~\eqref{eq:gov} satisfies, in the weak sense,
\begin{align*}
\frac{\mathrm{d}}{\mathrm{d} z}   \left(
A(z)\left(J(\phi,z)
-\gamma(z) \frac{\mathrm{d}  D (\phi)}{\mathrm{d}   z }  
\right) - Q_\mathrm{F}\phi_\mathrm{F}H(z-z_\mathrm{F})
\right) =0,\quad z\in\mathbb{R}. 
\end{align*} 
Integrating this identity with respect to~$z$ yields 
\begin{align} 
A(z)\bigl(J(\phi,z)
-\gamma(z)d(\phi)  \phi'(z)  
\bigr) - Q_\mathrm{F}\phi_\mathrm{F}H(z-z_\mathrm{F})
=\mathcalold{M},
\quad z\in\mathbb{R},\label{eq:Phi}
\end{align}
where the constant mass flux~$\mathcalold{M}$ can be determined by  setting~$z$ to a value either less than~$z_\mathrm{U}$ or greater than $z_\mathrm{E}$; then one gets
\begin{align*}
\mathcalold{M}&=A_\mathrm{U}j_\mathrm{U}(\phi_\mathrm{U}) = -Q_\mathrm{U}\phi_\mathrm{U},\notag\\
\mathcalold{M}&=A_\mathrm{E}j_\mathrm{E}(\phi_\mathrm{E}) - Q_\mathrm{F}\phi_\mathrm{F}
=:\mathcalold{M}_\mathrm{E} - Q_\mathrm{F}\phi_\mathrm{F},
\end{align*}
where the effluent constant mass flux of aggregates $\mathcalold{M}_\mathrm{E}:=
A_\mathrm{E}j_\mathrm{E}(\phi_\mathrm{E})
=Q_\mathrm{E}\phi_\mathrm{E}$
is also the constant mass flux above the feed inlet.
For a desired steady state satisfying~\eqref{eq:desired}, we have
$\phi_\mathrm{U}=0$; hence, $\mathcalold{M}=0$ and the feed mass flux equals the effluent:
\begin{equation*}
\phi_\mathrm{U}=0\quad
\Leftrightarrow\quad
 Q_\mathrm{F}\phi_\mathrm{F}=\mathcalold{M}_\mathrm{E}.
\end{equation*}
It is convenient to define the feed mass flux per area unit by
\begin{equation}\label{eq:sF}
s_\mathrm{F}:=\frac{Q_\mathrm{F}\phi_\mathrm{F}}{A_\mathrm{E}}.
\end{equation}
With $z$ in zone~2, \eqref{eq:Phi} gives $\mathcalold{M}=A_\mathrm{E}(j_2(\phi)-D(\phi)') - Q_\mathrm{F}\phi_\mathrm{F}$, which with $\mathcalold{M}=0$ and \eqref{eq:sF} implies that the solution $\phi$ in zone~2 satisfies
\begin{equation}\label{eq:ODE1}
\begin{aligned}
j_2(\phi)-d(\phi) \phi'(z) 
&=s_\mathrm{F},
\quad z_\mathrm{F}<z<z_\mathrm{E},\\
s_\mathrm{F}&=q_\mathrm{E}\phi_\mathrm{E}.
\end{aligned}
\end{equation}

The boundary condition in~\eqref{eq:ODE1}  also implies that~$\phi_\mathrm{E}$ can be expressed in terms of given and control variables (recall that $Q_\mathrm{E}>0$):
\begin{equation}\label{eq:phiE}
\phi_\mathrm{E}
=\frac{A_\mathrm{E}s_\mathrm{F}}{Q_\mathrm{E}}
=\frac{Q_\mathrm{F}\phi_\mathrm{F}}{Q_\mathrm{W} +Q_\mathrm{F} -Q_\mathrm{U}}.
\end{equation}

Since we require that there be no aggregates in zone~1, the steady-state solution there is zero.
For any jump across $z=z_\mathrm{F}$ from this zero volume fraction to any larger value~$\bar{\phi}_2$, from which there should be a discontinuity in zone~2 at $z=z_\mathrm{fr}$, the bottom of the froth layer, the uniqueness condition~\citep{SDjhde1} implies that $\bar{\phi}_2$ has to lie on an increasing part of~$j_2(\cdot;q_2)$. 
This corresponds to cases (a) and (c) in~\cite[Section~3.2]{SDIMAflot2019}.
The latter case can only occur under special circumstances with a large $\bar{\phi}_2>\phi_{2\rm M}$ (see definition in Section~\ref{sec:jprop}). 
Any small disturbance in a volumetric flow will make the case impossible and we therefore ignore that case.
Consequently, we consider only $\bar{\phi}_2\in[0,\phi_2^\mathrm{M}]$.
Then $\bar{\phi}_2$ is the smallest positive solution of the jump condition equation at the feed level, namely  
\begin{equation}\label{eq:FJC}\tag{FJC}
s_\mathrm{F}=j_2(\phi;q_2),
\end{equation}
under the  conditions  \citep{SDIMAflot2019} 
\begin{align}
&s_\mathrm{F}\leq j_2\big(\phi_2^\mathrm{M};q_2\big),
\label{eq:FIa}\tag{FIa}\\
&\bar{\phi}_2\leq\phi_{1\rm Z},\label{eq:FIb}\tag{FIb}
\end{align}
where $\phi_2^\mathrm{M}$ and $\phi_{1\rm Z}$ are defined in Section~\ref{sec:jprop}.
By the properties of $j_2$, we have $\phi_2^\mathrm{M}\leq\phi_\mathrm{infl}<\phi_\mathrm{c}$.
Therefore, $\bar{\phi}_2<\phi_2^\mathrm{M}<\phi_\mathrm{c}$.
Then $d(\bar{\phi}_2)=0$, and the equation in \eqref{eq:ODE1} reduces to~\eqref{eq:FJC}.
Again with reference to~\cite{SDjhde1} and without going into details, we claim that the solution in zone~2 is
\begin{equation}\label{eq:phi2sol}
\phi_2(z)=
\begin{cases}
\bar{\phi}_2,\quad z_\mathrm{F}<z<z_\mathrm{E},&
\text{if $\phi_\mathrm{E}\leq\phi_\mathrm{c}$},\\
\begin{cases}
\bar{\phi}_2, & z_\mathrm{F}<z<z_\mathrm{fr},\\
\phi_{2\mathrm{par}}(z), & z_\mathrm{fr}<z\leq z_\mathrm{E},
\end{cases}
&\text{if $\phi_\mathrm{E}>\phi_\mathrm{c}$ and $z_\mathrm{fr}>z_\mathrm{F}$},
\end{cases}
\end{equation}
where $\phi_{2\mathrm{par}}(z)$ is the strictly increasing solution of the  ordinary 
 differential equation (see~\eqref{eq:ODE1}):
\begin{equation}\label{eq:ODE2}
\begin{aligned}
& \phi'(z) 
=\frac{j_2(\phi;q_2)-s_\mathrm{F}}{d(\phi)},\\
&\phi(z_\mathrm{fr})=\phi_\mathrm{c},
\quad
\phi(z_\mathrm{E})=\phi_\mathrm{E},
\end{aligned}
\end{equation}
and where $z_\mathrm{fr}$ is the unknown location of the pulp-froth interface $\phi= \phi_{\mathrm{c}}$, which depends on~$s_\mathrm{F}$ and $\phi_\mathrm{E}$.

\begin{figure}[t!]
\centering\includegraphics[width=0.9\textwidth]{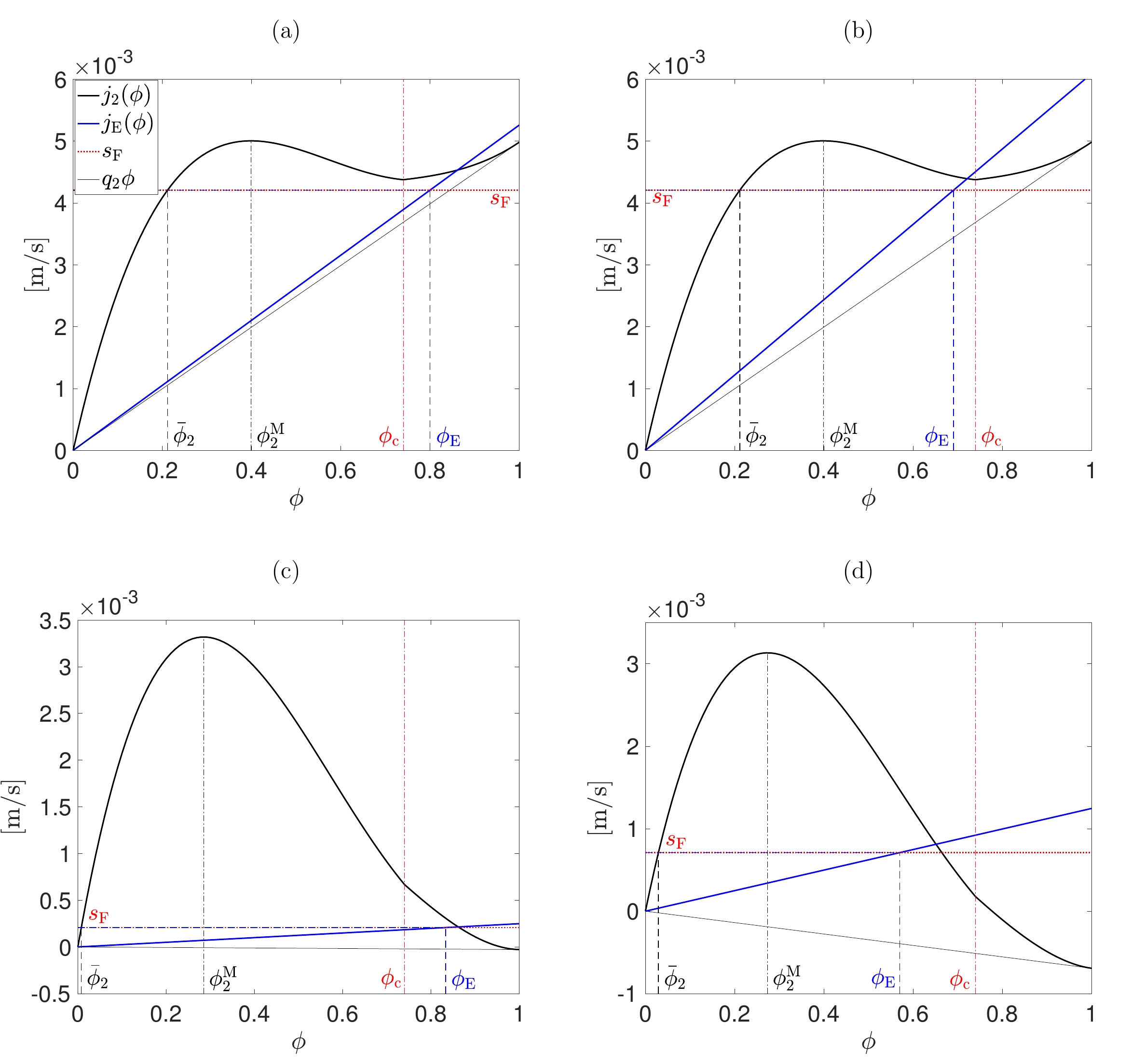}
\caption{Possible steady-state values for zone~2 with (a, b) $q_2>0$ and (c, d) $q_2<0$. 
The case $\phi_\mathrm{E}>\phi_\mathrm{c}$ is shown in (a) and (c), where there is a continuously increasing solution $\phi_\mathrm{par}(z)\in(\phi_\mathrm{c},\phi_\mathrm{E})$, while $\phi_\mathrm{E}\leq\phi_\mathrm{c}$ in (b) and (d), where the solution in zone~2 is the constant~$\bar{\phi}_2$.
For all the cases, we have $\phi_\mathrm{c}=0.74$, $n_\mathrm{b}=2.5$ and $n_\mathrm{S}=0.46$. 
For (a) and (b), $Q_2:=q_2A_\mathrm{E}=3.6 \times 10^{-5}\mathrm{m}^3/\mathrm{s}$, $s_\mathrm{F}=4.21 \times 10^{-3}\,\mathrm{m}/\mathrm{s}$ and (a) $Q_\mathrm{W}=2 \times 10^{-6}\,\mathrm{m}^3/\mathrm{s}$, (b) $Q_\mathrm{W}=8 \times 10^{-6}\,\mathrm{m}^3/\mathrm{s}$. 
For (c) we let $Q_2=-2 \times 10^{-6}\,\mathrm{m}^3/\mathrm{s}$, $Q_\mathrm{W}=2 \times 10^{-6}\,\mathrm{m}^3/\mathrm{s}$ and $s_\mathrm{F}=2.07 \times 10^{-4}\,\mathrm{m}/\mathrm{s}$, while for (d) we used $Q_2=-5 \times 10^{-6}\,\mathrm{m}^3/\mathrm{s}$, $Q_\mathrm{W}=10^{-5}\mathrm{m}^3/\mathrm{s}$ and $s_\mathrm{F}=7.1 \times 10^{-4}\,\mathrm{m}/\mathrm{s}$.
\label{fig:SS}}
\end{figure}%

See Figure~\ref{fig:SS} for illustrations of some steady-state solutions in zone~2.
In~\eqref{eq:phi2sol} lies the fact that if $\phi_\mathrm{E}>\phi_\mathrm{c}$, then there is no discontinuity at $z=z_\mathrm{E}$, so that $\phi_{2\mathrm{par}}(z_\mathrm{E}^-)=\phi_\mathrm{E}^-=\phi_\mathrm{E}^+=\phi_\mathrm{E}$.
The boundary value problem~\eqref{eq:ODE2} defines a function~$Z_\mathrm{fr}$ via
\begin{equation*}
z_\mathrm{fr} =Z_\mathrm{fr}(\phi_\mathrm{F},Q_\mathrm{F},Q_\mathrm{U},Q_\mathrm{W}).
\end{equation*}
In light of \eqref{eq:phiE} and \eqref{eq:phi2sol}, necessary conditions for a steady-state solution with a froth region are the inequalities
\begin{align}
\phi_\mathrm{c} &<\phi_\mathrm{E}\leq 1
\quad\Leftrightarrow\quad
Q_\mathrm{F}\left(1 - \frac{\phi_\mathrm{F}}{\phi_\mathrm{c}}\right)< Q_\mathrm{U} - Q_\mathrm{W}\leq Q_\mathrm{F}(1 - {\phi_\mathrm{F}}),\label{eq:Froth1}\tag{\rm Froth1}\\
z_\mathrm{F}&<Z_\mathrm{fr}(\phi_\mathrm{F},Q_\mathrm{F},Q_\mathrm{U},Q_\mathrm{W}) \label{eq:Froth2}\tag{\rm Froth2}
\end{align}
that should be satisfied for a steady-state solution with a froth-pulp interface in zone~2.
The requirement that $\phi_{2\mathrm{par}}(z)$ is strictly increasing from~$\phi_\mathrm{c}$ to $\phi_\mathrm{E}$ means that the left inequality of~\eqref{eq:Froth1} is equivalent to $Z_\mathrm{fr}(\phi_\mathrm{F},Q_\mathrm{F},Q_\mathrm{U},Q_\mathrm{W}) < z_\mathrm{E}$; hence, the latter inequality need not be invoked.

That $\phi_{2\mathrm{par}}(z)$ is strictly increasing, required by the entropy condition in~\cite{SDjhde1}, means that the right-hand side of~\eqref{eq:ODE2} is positive in the interval $[\phi_\mathrm{c},\phi_\mathrm{E})$.
Furthermore, a discontinuity at $z=z_\mathrm{fr}$ from $\bar{\phi}_2$ up to $\phi_\mathrm{c}$ can only occur (according to the entropy condition) if the graph of $j_2(\cdot;q_2)$ lies above $s_\mathrm{F}$ in the interval $(\bar{\phi}_2,\phi_\mathrm{c})$.
These conditions imply (FIa), which we can abandon.
The properties of $j_2(\cdot,;q_2)$ (see Section~\ref{sec:jprop}) imply that we can write these necessary conditions for a solution with a froth region:
\begin{align}
s_\mathrm{F} <j_2(\phi;q_2)\quad\text{for all $\phi\in(\phi_2^\mathrm{M},\phi_\mathrm{E})$} 
\quad\Leftrightarrow \quad  
s_\mathrm{F} \begin{cases}
<j_2(\phi_\mathrm{2M};q_2) &\text{if $\phi_\mathrm{2M}<\phi_\mathrm{E}$},\\
\leq j_2(\phi_\mathrm{E};q_2) &\text{if $\phi_\mathrm{2M}\geq\phi_\mathrm{E}$},
\end{cases}\label{eq:Froth3}\tag{Froth3}
\end{align}
where equality holds if and only if $\phi_\mathrm{2M}=\phi_\mathrm{E}$.

\subsection{Desired steady states and operating charts}\label{sec:desiredSS}

\begin{figure}[t] 
\centering\includegraphics[width=\textwidth]{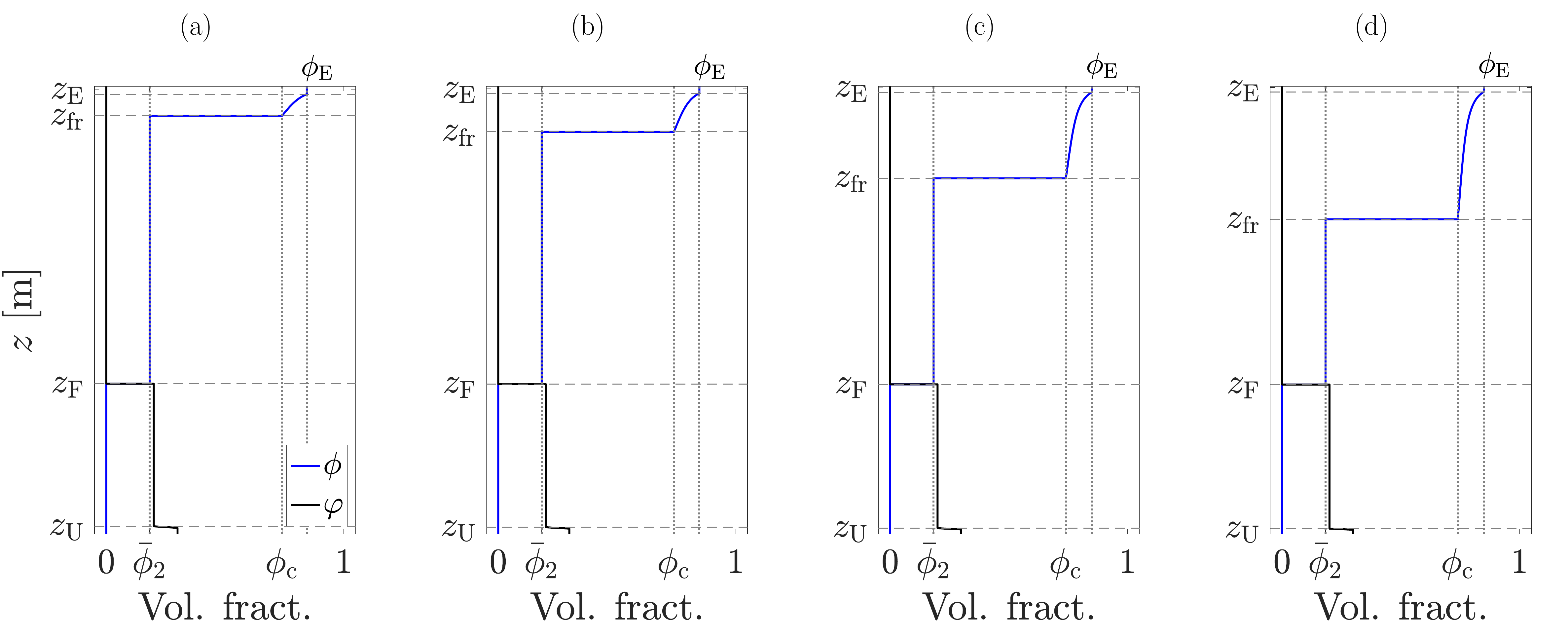}
\caption{Examples of desired steady states given by \eqref{eq:phiSS} and \eqref{eq:varphiSS}. We use fixed values of $\phi_\mathrm{F}=0.3$, $\psi_\mathrm{F}=0.2$, $Q_\mathrm{F}=8.9927 \times 10^{-5} \,\mathrm{m}^3/\mathrm{s}$ and $Q_\mathrm{W}=2 \times 10^{-6} \,\mathrm{m}^3/\mathrm{s}$ and vary $Q_\mathrm{U}$, choosing: (a) $Q_\mathrm{U}=5.9972 \times 10^{-5} \,\mathrm{m}^3/\mathrm{s}$, (b) $Q_\mathrm{U}=6.0083 \times 10^{-5} \,\mathrm{m}^3/\mathrm{s}$, (c) $Q_\mathrm{U}=6.0155 \times 10^{-5} \,\mathrm{m}^3/\mathrm{s}$ and (d) $Q_\mathrm{U}=6.0171 \times 10^{-5} \,\mathrm{m}^3/\mathrm{s}$. Once the values of $\phi_\mathrm{F}$, $Q_\mathrm{U}$, $Q_\mathrm{F}$ and $Q_\mathrm{W}$ are chosen, the values of the effluent concentration $\phi_\mathrm{E}$ are given by \eqref{eq:phiE} and used as input in the ODE \eqref{eq:ODE2} to calculate the value of $z_\mathrm{fr}$. In particular, we get (a) $\phi_\mathrm{E}=0.8443$, (b) $\phi_\mathrm{E}=0.8472$, (c) $\phi_\mathrm{E}=0.8491$ and (d) $\phi_\mathrm{E}=0.8495$. The values of $\phi_\mathrm{F}$, $\psi_\mathrm{F}$, $Q_\mathrm{U}$, $Q_\mathrm{F}$ and $Q_\mathrm{W}$ chosen here are used in Example~1 in Section~\ref{sec:numexp} to recover these profiles using the numerical method proposed in Section~\ref{sec:num_method}.\label{fig:desiredSS}}
\end{figure}

From the derivation above concerning the aggregates and from the  treatment in \cite{SDIMAflot2019} concerning the solids, we here summarize the \emph{desired steady states} that satisfy~\eqref{eq:desired}:
\begin{align}
\phi_\mathrm{SS}(z)&=
\begin{cases}
0 & \text{for $z_\mathrm{U}<z<z_\mathrm{F}$,} \\
\bar{\phi}_2 & \text{for $z_\mathrm{F}<z<z_\mathrm{fr}$,} \\
\phi_{2\mathrm{par}}(z)  & \text{for $z_\mathrm{fr}<z<z_\mathrm{E}$,} \label{eq:phiSS}\\
\phi_\mathrm{E} & \text{for $z>z_\mathrm{E}$,} 
\end{cases}\\
\varphi_\mathrm{SS}(z)&=\begin{cases}
0 & \text{for $z>z_\mathrm{F}$,} \\
\varphi_1\in[0,\varphi_{1\mathrm{m}}] & \text{for $z_\mathrm{U}<z<z_\mathrm{F}$,} \label{eq:varphiSS}\\
\varphi_\mathrm{U} = \varphi_1 + {A_\mathrm{U}}f_\mathrm{b}(\varphi_1)/{Q_\mathrm{U}}& \text{for $z<z_\mathrm{U}$.}  
\end{cases}
\end{align}
Here, $\phi_{2\mathrm{par}}(z)$ is the solution of the ODE problem~\eqref{eq:ODE2}, $\phi_\mathrm{E}$ is given by~\eqref{eq:phiE}, $\varphi_{1\mathrm{m}}$ is given by~\eqref{eq:varphi1m} and $\varphi_1>0$ satisfies the jump condition at the feed level $z=z_\mathrm{F}$ ($\varphi_1$ is unique if condition~\eqref{eq:FIas} below holds; see \cite{SDIMAflot2019})
\begin{equation}\label{eq:FJCs}
Q_\mathrm{F}\psi_\mathrm{F}=A_\mathrm{U}f_1(\varphi_1,0;q_1).
\tag{FJCs}
\end{equation}
In Figure~\ref{fig:desiredSS}, we have represented some examples of desired steady states with different values of $z_\mathrm{fr}$ obtained by fixing the values of the parameters $\phi_\mathrm{F}$, $\psi_\mathrm{F}$, $Q_\mathrm{F}$ and $Q_\mathrm{W}$ and choosing different values for $Q_\mathrm{U}$. As it can be seen, the location of $z_\mathrm{fr}$ is very sensitive to the choice of $Q_\mathrm{U}$. For instance, it changes from $z_\mathrm{fr}=0.8027\,\mathrm{m}$ in (c) to $z_\mathrm{fr}=0.7081\,\mathrm{m}$ in (d) with a small variation in $Q_\mathrm{U}$ of $-1.6 \times 10^{-8} \,\mathrm{m}^3/\mathrm{s}$.
We now collect the conditions for obtaining a desired steady state in terms of the input and control variables.

\begin{theorem}\label{theorem1}
The desired steady-state solution~\eqref{eq:phiSS} and \eqref{eq:varphiSS} of the PDE system~\eqref{eq:gov} is possible only if the following inequalities are satisfied:
\begin{align}
&\bar{\phi}_2\leq\phi_{\rm Z}(-Q_\mathrm{U}/A_\mathrm{U}),\tag{FIb}\label{eq:FIb2}\\
&A_\mathrm{U}f_1\big(\varphi_\mathrm{M}(-Q_\mathrm{U}/A_\mathrm{U}),0;-Q_\mathrm{U}/A_\mathrm{U}\big) \geq Q_\mathrm{F}\psi_\mathrm{F}.\tag{FIas}\label{eq:FIas}\\
&Q_\mathrm{F}\left(1 - \frac{\phi_\mathrm{F}}{\phi_\mathrm{c}}\right)< Q_\mathrm{U} - Q_\mathrm{W}\leq Q_\mathrm{F}(1 - {\phi_\mathrm{F}}),\tag{Froth1}\\
&z_\mathrm{F}<Z_\mathrm{fr}(\phi_\mathrm{F},Q_\mathrm{F},Q_\mathrm{U},Q_\mathrm{W}),\tag{Froth2}\\
&s_\mathrm{F}\begin{cases}
<j_2(\phi_\mathrm{2M};q_2) &\text{if $\phi_\mathrm{2M}<\phi_\mathrm{E}$,} \\
\leq j_2(\phi_\mathrm{E};q_2) &\text{if $\phi_\mathrm{2M}\geq\phi_\mathrm{E}$,} 
\end{cases}\tag{Froth3}
\end{align}
where we recall the definitions of $q_2$~\eqref{eq:q} and $\phi_\mathrm{E}$~\eqref{eq:phiE}:
\begin{equation*}
q_2=\frac{-Q_\mathrm{U}+Q_\mathrm{F}}{A_\mathrm{E}},
\qquad
\phi_\mathrm{E}
=\frac{Q_\mathrm{F}\phi_\mathrm{F}}{Q_\mathrm{W} +Q_\mathrm{F} -Q_\mathrm{U}}.
\end{equation*}
\end{theorem}

\begin{figure}[t!] 
\centering\includegraphics[width=0.85\textwidth]{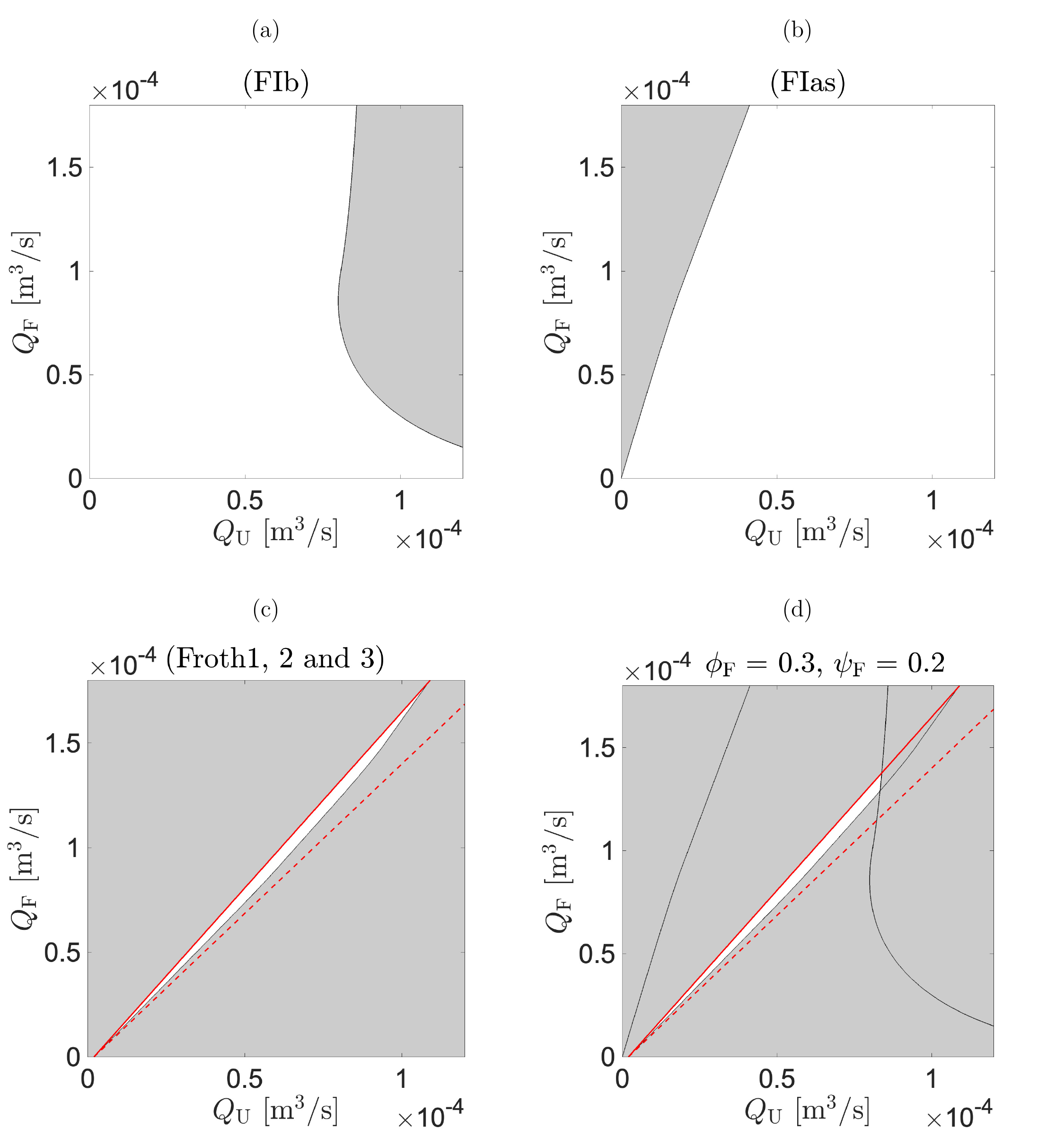}
\caption{(a--c) Visualization of the conditions of Theorem~\ref{theorem1} for $Q_\mathrm{W}=2 \times 10^{-6} \,\mathrm{m}^3/\mathrm{s}$, $\phi_\mathrm{F}=0.3$ and $\psi_\mathrm{F}=0.2$. (d)~Operating chart showing the intersection of all the conditions, which are true in the white region.\label{fig:opchart_cond}}
\end{figure}%

Inequalities~\eqref{eq:FIb} and \eqref{eq:FIas} can also be found in~\citep{SDIMAflot2019}.
We visualize them  together with \eqref{eq:Froth1}, \eqref{eq:Froth2} and \eqref{eq:Froth3} in the $(Q_\mathrm{U},Q_\mathrm{F})$-plane for fixed values of $Q_\mathrm{W}$, $\phi_\mathrm{F}$, $\psi_\mathrm{F}$ and $z_\mathrm{F}$; see Figure~\ref{fig:opchart_cond} for the choices $Q_\mathrm{W}=2 \times 10^{-6} \,\mathrm{m}^3/\mathrm{s}$, $\phi_\mathrm{F}=0.3$, $\psi_\mathrm{F}=0.2$,  and $z_\mathrm{F}=0.33 \,\mathrm{m}$.
All the conditions are shown together in Figure~\ref{fig:opchart_cond}~(d), which we call an {\em operating chart}.
For any chosen point $(Q_\mathrm{U},Q_\mathrm{F})$ in the white region, where all conditions in Theorem~\ref{theorem1} are satisfied, a desired steady-state solution given by~\eqref{eq:phiSS} and \eqref{eq:varphiSS} can be reached.

\begin{figure}[t!] 
\centering\includegraphics[width=0.85\textwidth]{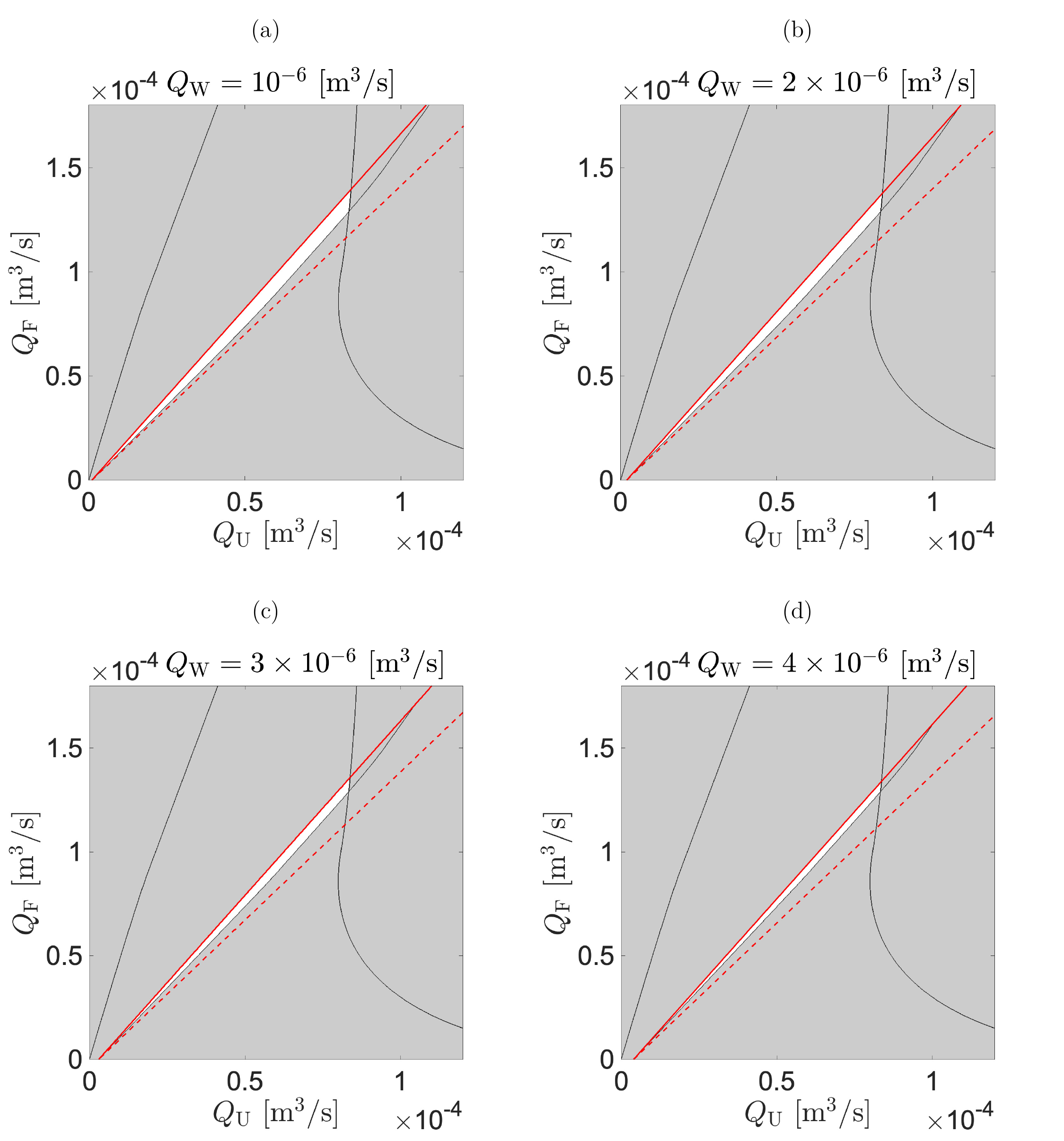}
\caption{Dependence of the operating chart on the wash water flow~$Q_\mathrm{W}$ for $\phi_\mathrm{F}=0.3$ and $\psi_\mathrm{F}=0.2$. \label{fig:opchart1}}
\end{figure}%

The two inequalities in~\eqref{eq:Froth1} give rise to a wedge-shaped region with vertex at $(Q_\mathrm{U},Q_\mathrm{F})=(Q_\mathrm{W},0)$; see Figure~\ref{fig:opchart_cond}~(c).
Thus, each wedge displayed in Figure~\ref{fig:opchart1} corresponds to a fixed value of~$Q_\mathrm{W}$, which can be read off at its vertex on the $Q_\mathrm{U}$-axis. 
The strict inequality of~\eqref{eq:Froth1} corresponds to the lower dashed line of a wedge, and its slope is positive or negative depending on whether~$\phi_\mathrm{F}$ is greater or less than $\phi_\mathrm{c}$.
The difference in slope of the two lines is $\phi_\mathrm{F}(1/\phi_\mathrm{c}-1)$, so the angle of the wedge increases with $\phi_\mathrm{F}$ and decreases with $\phi_\mathrm{c}$. 
The lower part of the wedge is, however, cut off by conditions~\eqref{eq:Froth2} and \eqref{eq:Froth3}; see Figure~\ref{fig:opchart_cond}~(c).
Figure~\ref{fig:opchart1} shows also that the white region of the operating chart thins and will eventually disappear as $Q_\mathrm{W}$ increases.

\begin{figure}[t!] 
\centering
\includegraphics[width=0.9\textwidth]{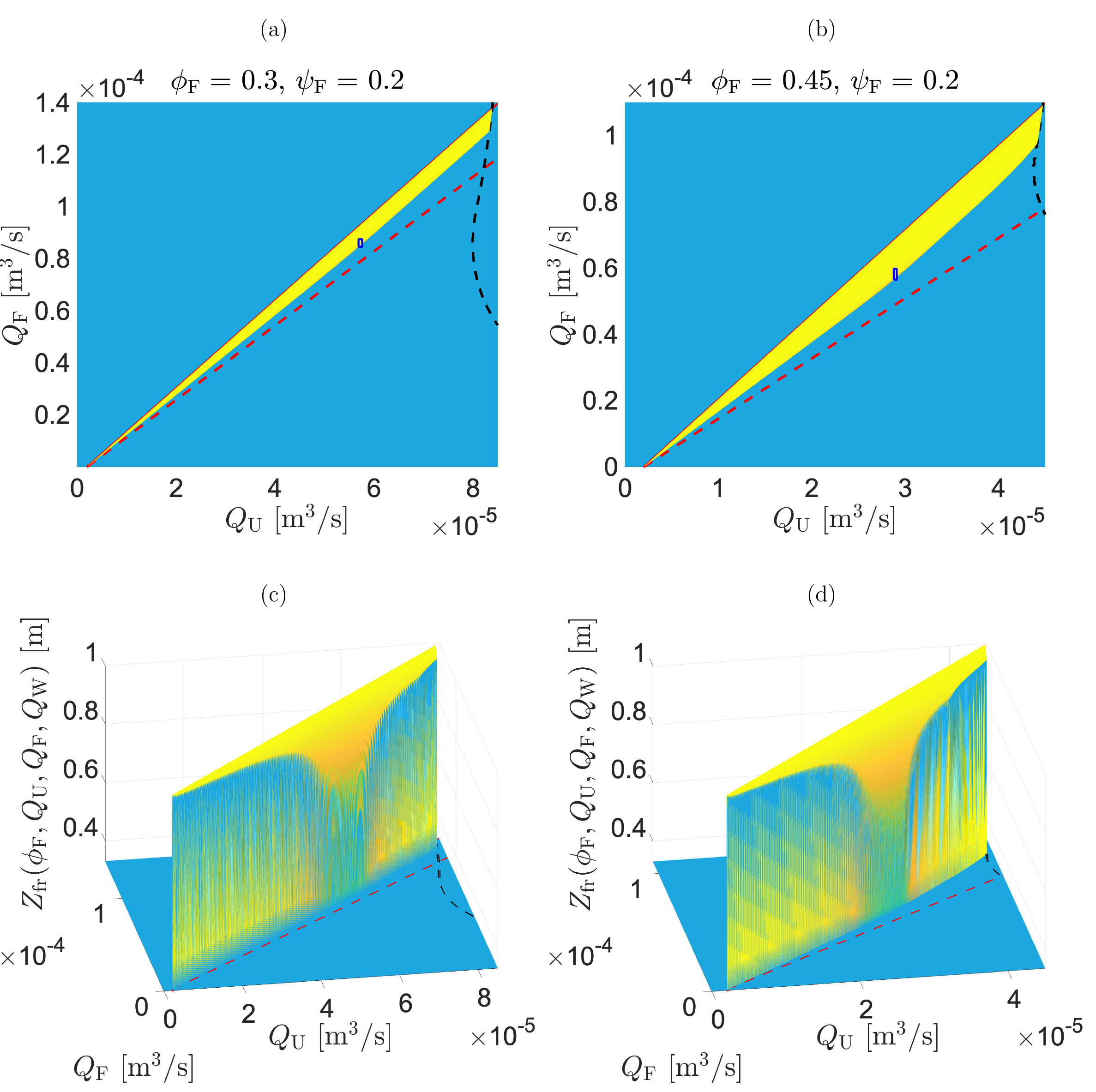}
\caption{Operating charts for $Q_\mathrm{W}=2 \times 10^{-6} \,\mathrm{m}^3/\mathrm{s}$ and $\psi_\mathrm{F}=0.2$ with (a, c) $\phi_\mathrm{F}=0.3$ and (b, d) $\phi_\mathrm{F}=0.45$, showing  
the graphs of $(Q_\mathrm{F},Q_\mathrm{U})\mapsto Z_\mathrm{fr}(\phi_\mathrm{F},Q_\mathrm{F},Q_\mathrm{U}, Q_\mathrm{W})$ obtained by~\eqref{eq:ODE2}.
The small rectangles in (a, b) are enlarged in Figure~\ref{fig:opchart2b}. 
\label{fig:opchart2}}
\end{figure}%

\begin{figure}[t] 
\centering
\includegraphics[width=0.8\textwidth]{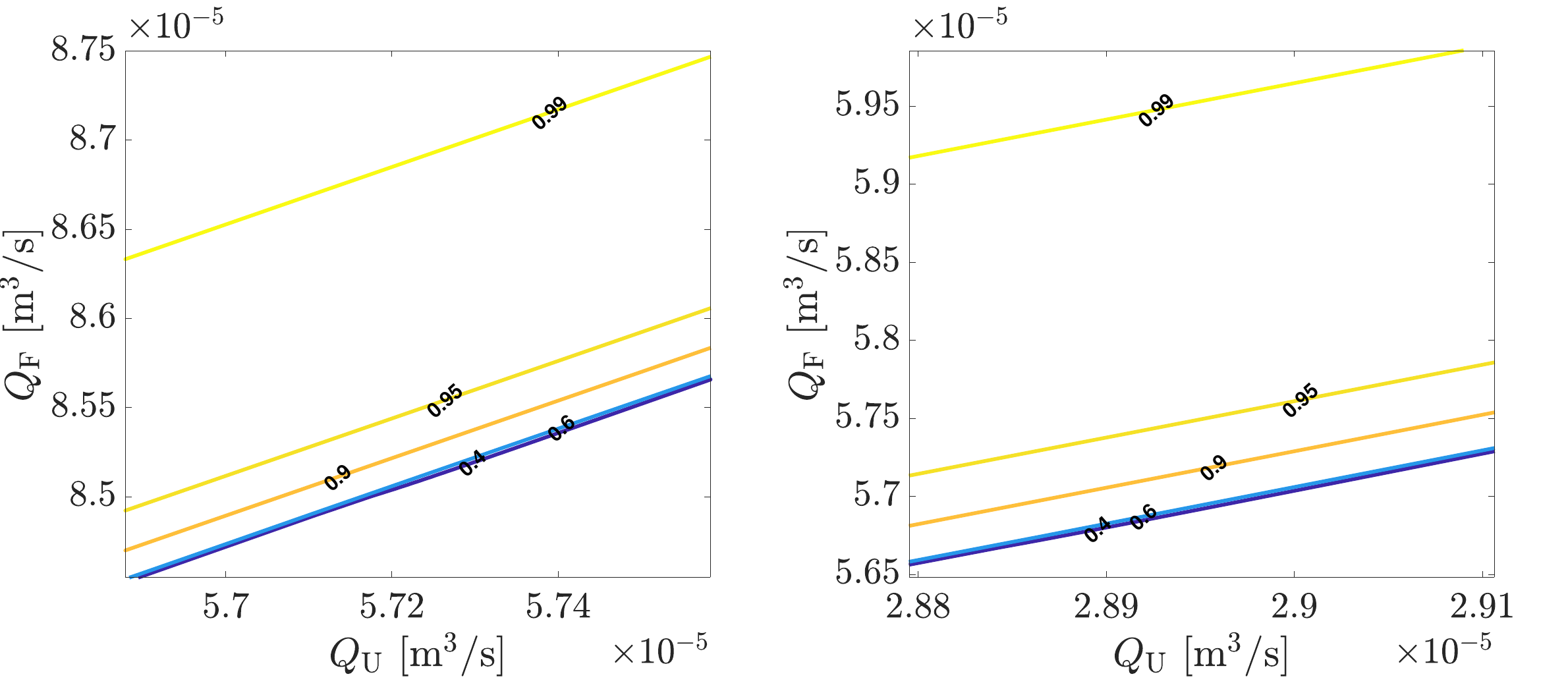}
\caption{Enlarged views of the small rectangles marked in Figure~\ref{fig:opchart2}~(a) and~(b), respectively, showing contours of the function $(Q_\mathrm{F},Q_\mathrm{U})\mapsto Z_\mathrm{fr}(\phi_\mathrm{F},Q_\mathrm{F},Q_\mathrm{U}, Q_\mathrm{W})$.
\label{fig:opchart2b}}
\end{figure}%

Inequality~\eqref{eq:Froth2} is more involved than the others.
For every given set of input and control values, one has to integrate the ODE of~\eqref{eq:ODE2} backwards from $z=z_\mathrm{E}$ for given $\phi_\mathrm{E}$ towards lower $z$-values until $\phi_\mathrm{c}$ is reached; the corresponding location defines $z=z_\mathrm{fr}$. In Figure~\ref{fig:opchart2}, the surface $z=Z_\mathrm{fr}(\phi_\mathrm{F},Q_\mathrm{F},Q_\mathrm{U},Q_\mathrm{W})$ has been computed for two different values of $\phi_\mathrm{F}$ and fixed values of $Q_\mathrm{W}$ and $\psi_\mathrm{F}$.
The same red and black curves as in Figures~\ref{fig:opchart_cond} and~\ref{fig:opchart1} of conditions \eqref{eq:Froth1} and \eqref{eq:FIb2}, respectively, limit the white region in the operating charts.
The study of the surface $Z_\mathrm{fr}(\phi_\mathrm{F},Q_\mathrm{F},Q_\mathrm{U},Q_\mathrm{W})$ is crucial for the 
 choice of  values of $Q_\mathrm{F}$, $Q_\mathrm{U}$ and $Q_\mathrm{W}$ for which we obtain a desired steady state with a froth region in $z=z_\mathrm{fr}$, with $z_\mathrm{fr}$ a given value, as we will see in the numerical results in Section \ref{sec:numexp}.

\section{Numerical method}\label{sec:num_method}

\subsection{Discretization and CFL condition}\label{sec:discret}

\begin{figure}[t] 
\centering\includegraphics[height=0.45\textheight]{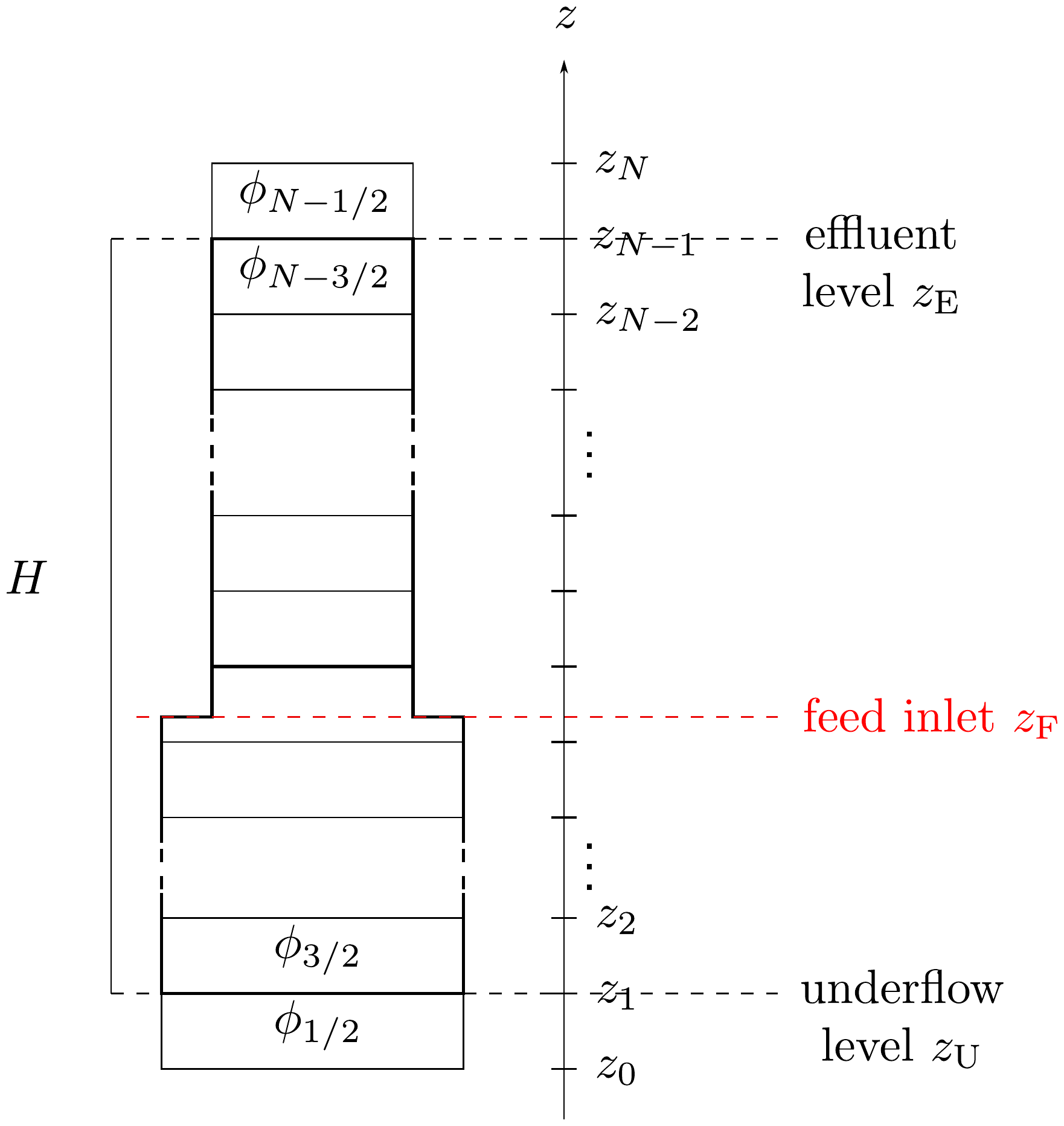} 
\caption{Grid  covering the flotation column for the discretization of $\phi$ and $\psi$.
The outlets $z_\mathrm{U}$ and $z_\mathrm{E}$ are each fixed on the boundaries between two cells and the feed inlet~$z_\mathrm{F}$ is then located in a cell.} 
\label{fig:mesh}
\end{figure}%

We define a computational domain of $N$ cells by covering the vessel with $N-2$ cells and placing one cell each below and above for the calculation of the outlet volume fractions see Figure~\ref{fig:mesh}.
Given the column height~$H$, we define $\Delta z:=H/(N-2)$ and the cell 
boundaries $z_{i}:=i\Delta z$, $i=0,1,\ldots,N$.
Furthermore, we define the cell intervals $I_{i-1/2}:=[z_{i-1},z_{i})$  and $I_{i}:=[z_{i-1/2},z_{i+1/2})$.
We place the column between $z_\mathrm{U}:=\Delta z=z_{1}$ and $z_\mathrm{E}:=z_\mathrm{U}+H=(N-1)\Delta z=z_{N-1}$.
The injection point $z_{\mathrm{F}}$ is assumed to belong to one cell $I_{i-1/2}$ and we define the dimensionless function
\begin{align*}
\delta_{\mathrm{F},i-1/2}:= \int_{I_{i-1/2}}^{}\delta_{z_{\mathrm{F}}}(z)\,\mathrm{d}z := \begin{cases}
1 &\text{if $z_{\mathrm{F}}\in I_{i-1/2}$,} \\
0 &\text{otherwise}.
\end{cases}
\end{align*}

The cross-sectional area $A= A(z)$ is allowed to have a finite number of discontinuities and it is discretized by
\begin{align*}
A_{i}:=\frac{1}{\Delta z}\int_{I_i}^{}A(z)\, \mathrm{d}z, \qquad
A_{i+1/2}:=\frac{1}{\Delta z}\int_{I_{i+1/2}}A(z)\, \mathrm{d}z.
\end{align*}

We simulate $N_T$ time steps up to the final time $T:=N_T\Delta t$, with the fixed time step~$\Delta t$ satisfying the Courant-Friedrichs-Lewy (CFL) condition 
\begin{equation}\label{eq:CFL}\tag{CFL}
{\Delta t}\left(\frac{2\|Q\|_{\infty,T}}{A_\mathrm{min}} 
+ M_1\|\tilde{v}'\|_\infty
+ \max\{\beta_1,\beta_2\}\right) \leq\Delta z,
\end{equation}
where 
\begin{equation*}
\begin{split} 
& \beta_1:=M_1 \|\tilde{v}\|_\infty
+ M_2\dfrac{\|d\|_\infty}{\Delta z},
\qquad
\beta_2:= M_1\max\left\{v_\mathrm{hs}(0),\|v_\mathrm{hs}'\|_\infty\right\}
+ M_2(1-\phi_\mathrm{c})\dfrac{\|d\|_\infty}{\Delta z},
\\
& M_1 :=\max\limits_{i=1,2,\ldots,N}
\left\{\frac{A_{i-1}}{A_{i-1/2}},\frac{A_{i}}{A_{i-1/2}}\right\}, \quad
M_2 :=\max\limits_{i=1,2,\ldots,N}
\left\{\frac{A_{i-1}+A_{i}}{A_{i-1/2}}\right\},\\
&A_\mathrm{min} :=  \min\limits_{k=0,\frac12,1,\frac32,\ldots,N}A_{k},\qquad
\|Q\|_{\infty,T} :=\max\limits_{0\leq t\leq 
T}(Q_{\mathrm{F}}(t) + Q_{\mathrm{W}}(t)),
\qquad \|d\|_\infty :=\max_{0\leq\phi\leq 1}|d(\phi)|.
\end{split} 
\end{equation*}
Finally, we set $t^n:=n\Delta t$ for $n=0,1,\ldots,N_T$. 

The time-dependent feed functions are discretized as
\begin{align*}
Q_{\mathrm{F}}^n:=\frac{1}{\Delta t}\int_{t^n}^{t^{n+1}}Q_{\mathrm{F}}(t)\,\mathrm{d}t,\qquad \phi_{\mathrm{F}}^n:=\frac{1}{\Delta t}\int_{t^n}^{t^{n+1}}\phi_{\mathrm{F}}(t)\,\mathrm{d}t,
\end{align*}
and the same is made for $\psi_{\mathrm{F}}$.

\subsection{Update of $\phi$}\label{sec:updatephi}

The first equation of~\eqref{eq:gov} depends only on $\phi$ and is discretized by a simple scheme on the cells $I_{i-1/2}$.
The initial data are discretized by
\begin{equation*} 
\phi_{i-1/2}^{0}:=\frac{1}{A_{i-1/2}\Delta z}\int_{I_{i-1/2}}^{}\phi(z,0)A(z)\,\mathrm{d}z.
\end{equation*}
To advance from $t^n$ to $t^{n+1}$, we assume that 
$\phi_{i-1/2}^{n}$, $i=1,\ldots, N$, are given.
With the notation 
\begin{align*} 
a^+:=\max\{ a,0 \}, \quad a^-:=\min \{ a,0\}, \quad \gamma_{i}:=\gamma(z_{i}), 
 \quad \text{and}\quad  q_{i}^{n+}:=q(z_{i},t^n)^+,
 \end{align*} 
we define the numerical total flux at $z=z_{i}$ at time $t=t^n$ by
\begin{equation}\label{eq:Jnumflux}
\Phi_{i}^n
:= \begin{cases}
\phi_{1/2}^nq_{0}^{n-} & \text{for $i=0$,} \\
\phi_{i-1/2}^nq_{i}^{n+} + \phi_{i+1/2}^nq_{i}^{n-} +\gamma_{i}\phi_{i-1/2}^n \tilde{v}(\phi_{i+1/2}^{n})
- \gamma_{i}\dfrac{D(\phi_{i+1/2}^n) - D(\phi_{i-1/2}^n)}{\Delta z}&
 \text{for $i=1,\ldots,N-1$,} \\
\phi_{N-1/2}^nq_{N}^{n+}
& \text{for $i=N$,}
\end{cases}
\end{equation}
where $\tilde{v}(\phi)$ is defined by~\eqref{eq:vtilde}.
Since the bulk fluxes above and below the tank are directed away from it, the following terms that appear in~\eqref{eq:Jnumflux} are zero:
\begin{align*} 
\phi_{-1/2}^nq_{0}^{n+}=0 \quad \text{and} \quad \phi_{N+1/2}^nq_{N}^{n-}=0
 \quad \text{for any values of $\smash{\phi_{-1/2}^n} \text{ and } \smash{\phi_{N+1/2}^n}$.}
 \end{align*} 
To simplify the presentation, we  use the middle line of~\eqref{eq:Jnumflux} as the definition of~$\Phi_{i}^n$, $i=0,\ldots,N$, together with $\smash{\phi_{-1/2}^n:=0}$ and $\smash{\phi_{N+1/2}^n:=0}$.
With the notation $\smash{\lambda:=\Delta t/\Delta z}$ and $\smash{Q_{i}^{n+}:=A_iq_{i}^{n+}}$ etc., the conservation law on the interval~$I_{i-1/2}$ implies the update formula
\begin{align}\label{eq:phi_updateR} 
\begin{split} 
\phi_{i-1/2}^{n+1}& =\phi_{i-1/2}^{n}+\frac{\lambda}{A_{i-1/2}}
\big(A_{i-1}\Phi_{i-1}^{n} -A_{i}\Phi_{i}^{n} + Q_{\mathrm{F}}^n \phi_{\mathrm{F}}^n\delta_{\mathrm{F},i-1/2}\big)\\
& =\phi_{i-1/2}^{n}+\frac{\lambda}{A_{i-1/2}}
\bigg( \phi_{i-3/2}^nQ_{i-1}^{n+} + \phi_{i-1/2}^nQ_{i-1}^{n-} +(A\gamma)_{i-1}\phi_{i-3/2}^n \tilde{v}(\phi_{i-1/2}^{n})\\
& \quad - \dfrac{(A\gamma)_{i-1}}{\Delta z}(D(\phi_{i-1/2}^n) - D(\phi_{i-3/2}^n)) -\phi_{i-1/2}^nQ_{i}^{n+} - \phi_{i+1/2}^nQ_{i}^{n-} -(A\gamma)_{i}\phi_{i-1/2}^n \tilde{v}(\phi_{i+1/2}^{n}) \\
& \quad +\dfrac{(A\gamma)_{i}}{\Delta z}(D(\phi_{i+1/2}^n) - D(\phi_{i-1/2}^n))  + Q_{\mathrm{F}}^n \phi_{\mathrm{F}}^n\delta_{\mathrm{F},i-1/2}\biggr), 
\quad i=1,\ldots, N. 
\end{split} 
\end{align}

\begin{theorem}\label{thm:bounded_phiR} 
If the CFL condition~\eqref{eq:CFL} is satisfied and the initial data satisfy $0\leq\phi(z,0)\leq 1$, then the update formula 
for~$\phi$, \eqref{eq:phi_updateR}, is monotone and produces approximate solutions that satisfy 
\begin{align}  \label{phibound} 
 0\leq\phi_{i-1/2}^n \leq 1 \quad \text{for $i=1,\ldots,N$ and $n =1, \dots, N_T$}.
 \end{align}  
\end{theorem}

The proof is outlined  in Appendix~A.

\subsection{Update of $\psi$}\label{sec:updatepsi}

We discretize the initial data by
\begin{align*}
\psi_{i-1/2}^{0}:=\frac{1}{A_{i-1/2}\Delta z}\int_{I_{i-1/2}}^{}\psi(z,0)A(z)\,\mathrm{d}z.
\end{align*}
A consistent numerical flux corresponding to~\eqref{eq:Ftotalflux} is, for $i=0,\ldots,N$,
\begin{align*}
\Psi_{i}^n
 & :=\psi_{i-1/2}^n q_{i}^{n+} + \psi_{i+1/2}^n q_{i}^{n-} \\
& \quad  -\gamma_{i} \bigg( G_{i}^{n}\big(\psi_{i-1/2}^{n},\psi_{i+1/2}^{n}\big) 
+\dfrac{\psi_{i+1/2}^n}{1-\phi_{i+1/2}^{n}}
\bigg(\phi_{i-1/2}^n\tilde{v}(\phi_{i+1/2}^{n} ) 
- \dfrac{\Delta D_i^{n-}}{\Delta z}\bigg)
- \dfrac{\psi_{i-1/2}^n}{1-\phi_{i-1/2}^{n}}\dfrac{\Delta D_i^{n+}}{\Delta z}
\bigg),
\end{align*}
where $\Delta D_i^n:=D(\phi_{i+1/2}^n) - D(\phi_{i-1/2}^n)$, $\Delta D_i^{n-}:=(\Delta D_i^n)^-$, and we set 
\begin{align*} 
\psi_{-1/2}^n:=0 \quad \text{and} \quad  \psi_{N+1/2}^n:=0
\end{align*} 
 with the same motivation as for $\phi$ above (these values are irrelevant). Here 
$\smash{G_{i}^{n} (\psi_{i-1/2}^{n},\psi_{i+1/2}^{n} )}$  
is the Engquist-Osher numerical flux \citep{EO81}  associated with the function
\begin{align} \label{fbi1} 
f_{\mathrm{b},i}^{n}(\psi):= \psi\tilde{v}_\mathrm{hs} 
\biggl(\frac{\psi}{\psi_{\mathrm{max},i}^n}\biggr), \qquad 
\tilde{v}_\mathrm{hs}(u)
:=\begin{cases}
v_\mathrm{hs}(u)  &\text{for $u<1$},\\
0 & \text{for $u\geq 1$},
\end{cases}
\end{align}
where we recall that~$v_\mathrm{hs}$ is given by~\eqref{eq:vhs}, and we define 
\begin{align} \label{fbi2} 
\psi_{\mathrm{max},i}^n:= \min \bigl\{1-\phi_{i-1/2}^n,1-\phi_{i+1/2}^n \bigr\} = 1-\max
 \bigl\{\phi_{i-1/2}^n,\phi_{i+1/2}^n \bigr\}.
\end{align}
If $\smash{\hat{\psi}_{i}^n}$ is the maximum point of~$f_{\mathrm{b},i}^{n}$, then the   Engquist-Osher numerical flux is given by 
\begin{align*}  
\begin{split} 
& G_{i}^{n}  \bigl(\psi_{i-1/2}^{n},\psi_{i+1/2}^{n} \bigr)   \\ & = 
\begin{cases} 
f_{\mathrm{b},i}^n ( \psi_{i+1/2}^n ) & \text{if $\psi_{i-1/2}^n, \psi_{i+1/2}^n \leq \hat{\psi}_i^n$,} \\
 f_{\mathrm{b},i}^n  (\hat{\psi}_i^n)  & 
 \text{if \phantom{$\psi_{i-1/2}^n,\,$}$\psi_{i-1/2}^n \leq \hat{\psi}_i^n<\psi_{i+1/2}^n$,} \\
-f_{\mathrm{b},i}^n  (\hat{\psi}_i^n)   + f_{\mathrm{b},i}^n ( \psi_{i-1/2}^n ) +  f_{\mathrm{b},i}^n (\psi_{i+1/2}^n)
 &  \text{if \phantom{$\psi_{i-1/2}^n,\,$}$\psi_{i+1/2}^n \leq  \hat{\psi}_i^n<\psi_{i-1/2}^n$,} \\
 f_{\mathrm{b},i}^n  (\psi_{i-1/2}^n) &  \text{if \phantom{$\psi_{i-1/2}^n, \psi_{i+1/2}^n\leq$} $\hat{\psi}_i^n<\psi_{i-1/2}^n, \psi_{i+1/2}^n$.}
\end{cases} \end{split} 
\end{align*} 
The marching formula is (for $i=1,\ldots, N$)
\begin{align} \label{eq:psi_updateR} 
\begin{split} 
& \psi_{i-1/2}^{n+1} \\ & =
\psi_{i-1/2}^{n} + \dfrac{\lambda}{A_{i-1/2}}\big( A_{i-1}\Psi_{i-1}^{n} - A_{i}\Psi_{i}^{n}
+  \displaystyle Q_{\mathrm{F}}^n\psi_{\mathrm{F}}^n\delta_{\mathrm{F},i-1/2} \big)\\
&=\psi_{i-1/2}^{n} + \dfrac{\lambda}{A_{i-1/2}}\bigg\{
\psi_{i-3/2}^{n}Q_{i-1}^{n+} + \psi_{i-1/2}^n Q_{i-1}^{n-}
-\psi_{i-1/2}^{n}Q_{i}^{n+} - \psi_{i+1/2}^n Q_{i}^{n-} +Q_{\mathrm{F}}^n\psi_{\mathrm{F}}^n\delta_{\mathrm{F},i-1/2}\\
& \quad -(A\gamma)_{i-1} 
\biggl( G_{i-1}^{n}\big(\psi_{i-3/2}^{n},\psi_{i-1/2}^{n}\big)
+\dfrac{\psi_{i-1/2}^n}{1-\phi_{i-1/2}^{n}}\bigg(\phi_{i-3/2}^n\tilde{v}(\phi_{i-1/2}^{n} ) 
- \dfrac{\Delta D_{i-1}^{n-}}{\Delta z}\bigg)
- \dfrac{\psi_{i-3/2}^n}{1-\phi_{i-3/2}^{n}}\dfrac{\Delta D_{i-1}^{n+}}{\Delta z}
\biggr) \\
& \quad +(A\gamma)_{i} 
\biggl( G_{i}^{n}\big(\psi_{i-1/2}^{n},\psi_{i+1/2}^{n}\big)
+\dfrac{\psi_{i+1/2}^n}{1-\phi_{i+1/2}^{n}}
\bigg(\phi_{i-1/2}^n\tilde{v}(\phi_{i+1/2}^{n} ) 
- \dfrac{\Delta D_i^{n-}}{\Delta z}\bigg)
-\dfrac{\psi_{i-1/2}^n}{1-\phi_{i-1/2}^{n}}\dfrac{\Delta D_i^{n+}}{\Delta z}
\biggr) \bigg\}.
\end{split} 
\end{align}

\begin{theorem}\label{thm:bounded_varphiR} 
Assume that the assumptions of Theorem~\ref{thm:bounded_phiR} are in effect. 
If the initial data satisfy $0\leq\psi(z,0)\leq 1-\phi(z,0)$ and the feed volume fraction $\psi_{\mathrm{F}}(t)\leq 1-\phi_{\mathrm{F}}(t)$, then the update formula~\eqref{eq:psi_updateR} is monotone and together with~\eqref{eq:phi_updateR} it produces approximate solutions that satisfy 
\begin{align*} 
0\leq\psi_{i-1/2}^n\leq 1-\phi_{i-1/2}^n \quad \text{for all~$i$ and~$n$.}
\end{align*}  
\end{theorem}

The proof is sketched  in Appendix~A.

\section{Numerical simulations} \label{sec:numexp}

We simulate the flotation process in the column in Figure~\ref{fig:Column} with the specific measures $A_\mathrm{E}=7.225\times 10^{-3} \,\mathrm{m}^2\leq A_\mathrm{U}=8.365\times 10^{-3} \,\mathrm{m}^2$, $z_\mathrm{U}=0\,$m, $z_\mathrm{F}=0.33\,$m,  $z_\mathrm{E}=1\,$m and $H=1\,$m.
For all the examples, we use the parameters given in \cite[Table~1]{Brito2012advantages} to define $v_\mathrm{drain}$ and $d_\mathrm{cap}$ in \eqref{eq:vtilde}--\eqref{eq:D}: $\rho_{\mathrm{f}} = 10^3 \, \mathrm{kg}/\mathrm{m}^3$, $\mu = 10^{-3} \, \mathrm{Pa} \,\mathrm{s}$,  $r_{\mathrm{b}} = 4.13 \times 10^{-4} \, \mathrm{m}$, $C_{\mathrm{PB}} = 50$,  $ \gamma_{\mathrm{w}}  = 3.5 \times 10^{-2} \, \mathrm{N} / \mathrm{m}$, $g = 9.81 \, \mathrm{m}/\mathrm{s}^2$, and by  \cite{Stevenson2004}, $n_\mathrm{S}=0.46$ and $m=1.28$, from which we obtain $d_\mathrm{cap}=3.1045 \times 10^{-3}$\,m. 
For the velocity functions $\tilde{v}$ and $v_{\mathrm{hs}}$, given by \eqref{eq:vtilde} and \eqref{eq:vhs}, respectively, we use $n_{\mathrm{b}} =2.5$, $v_{\mathrm{term}}=2.7\times 10^{-2}\, \mathrm{m}/\mathrm{s}$, $n_{\mathrm{RZ}}=1.5$ and $v_{\infty}=5.0\times 10^{-3}\, \mathrm{m}/\mathrm{s}$. The critical volume fraction is $\phi_{\mathrm{c}} = 0.74$ according to~\cite[Eq.~(21)]{Neethling2003}.

\subsection{Example 1}

\begin{figure}[t!] 
\centering 
\begin{tabular}{cc} 
(a) & (b) \\
\includegraphics[width=0.45\textwidth]{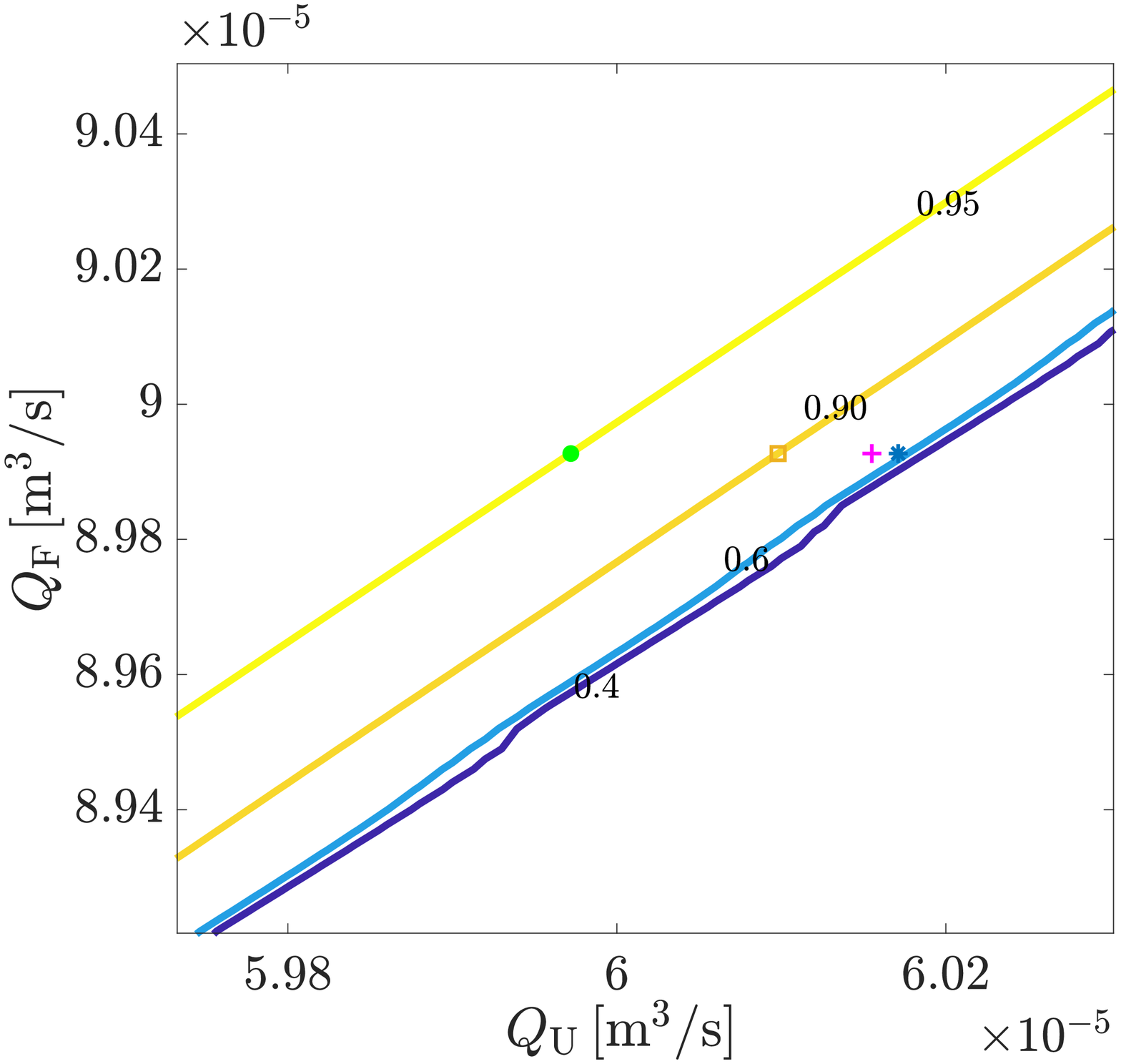} &  \includegraphics[width=0.45\textwidth]{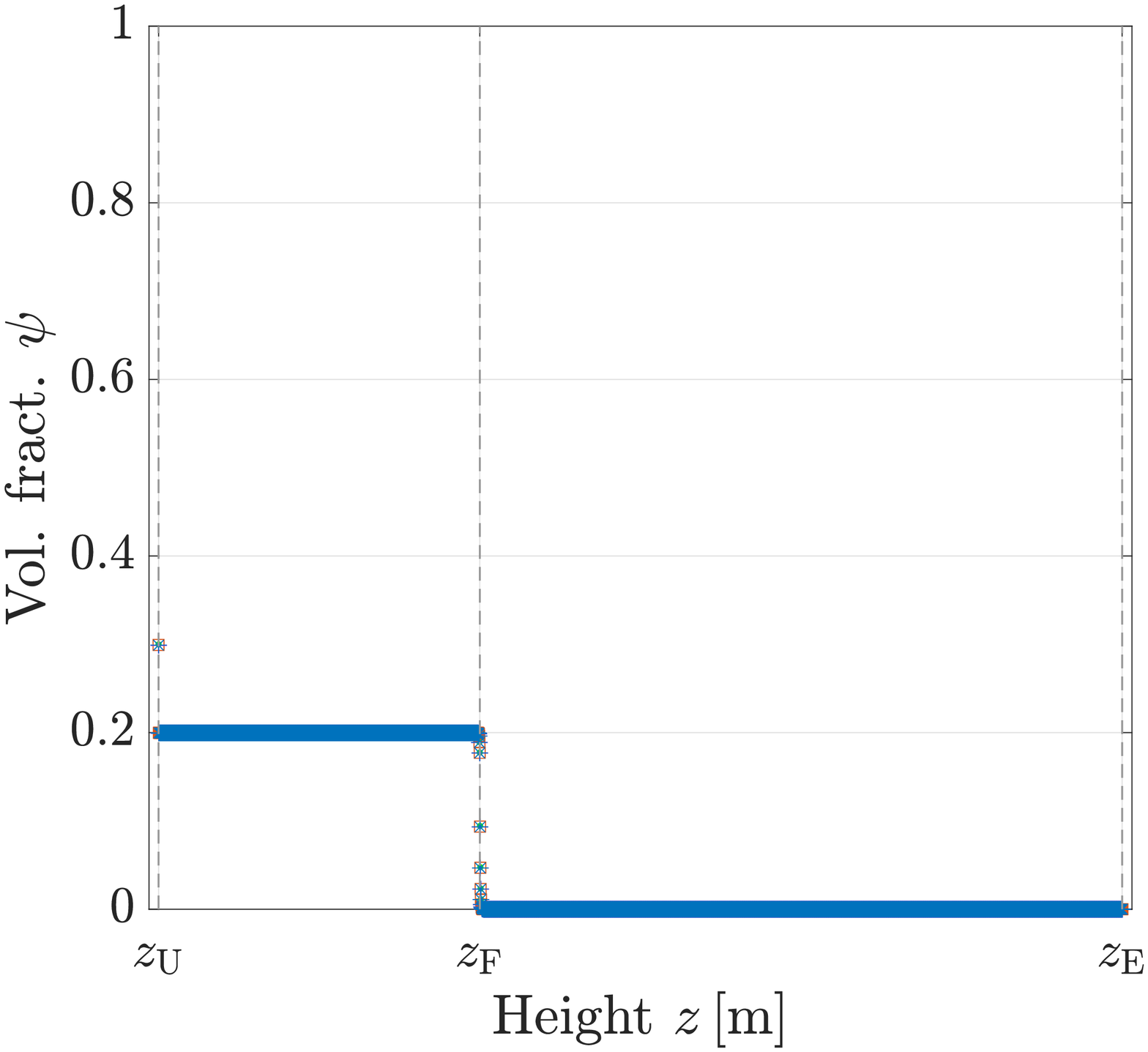}  \\
(c) & (d) \\  
\includegraphics[width=0.45\textwidth]{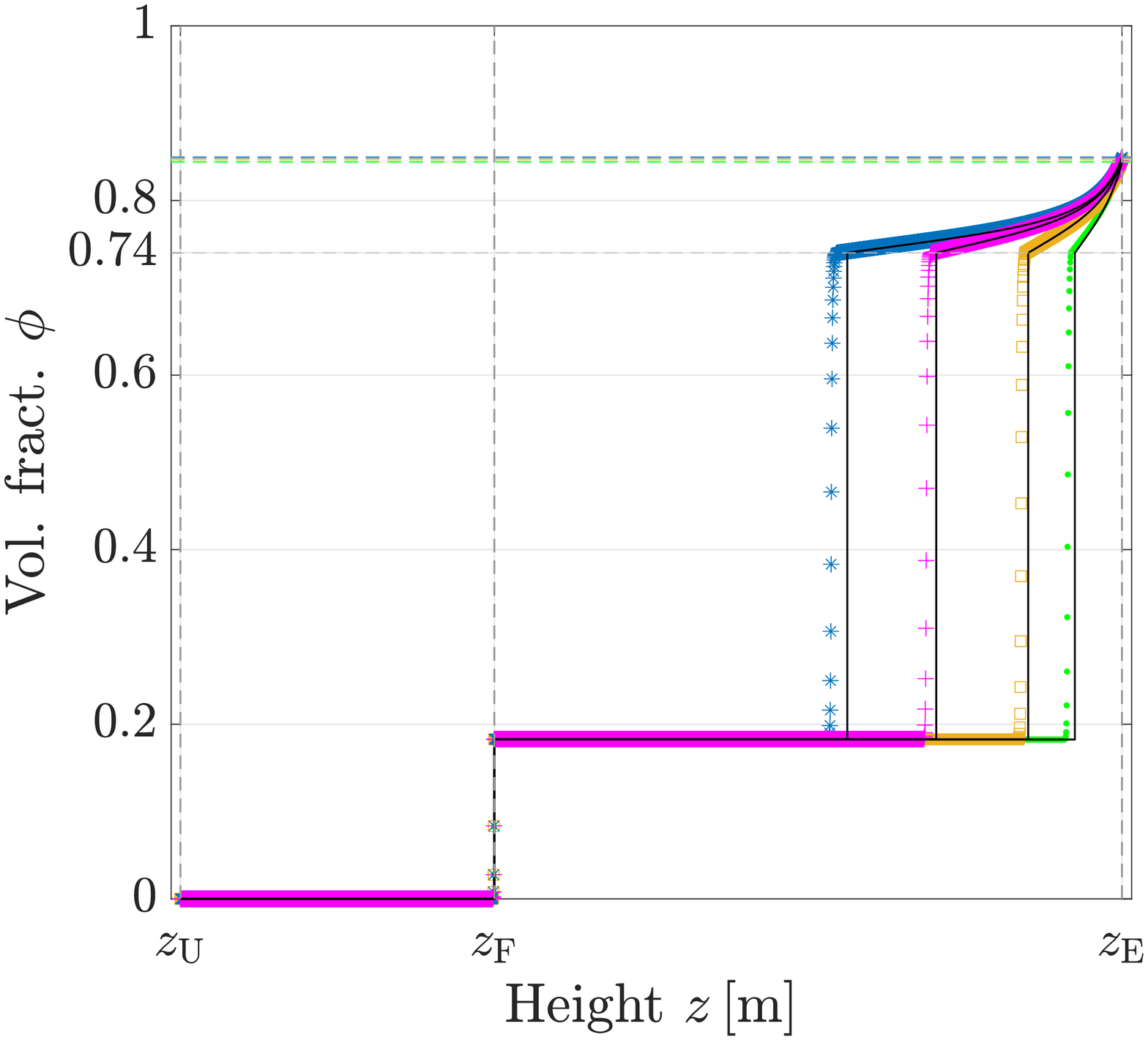} &  \includegraphics[width=0.45\textwidth]{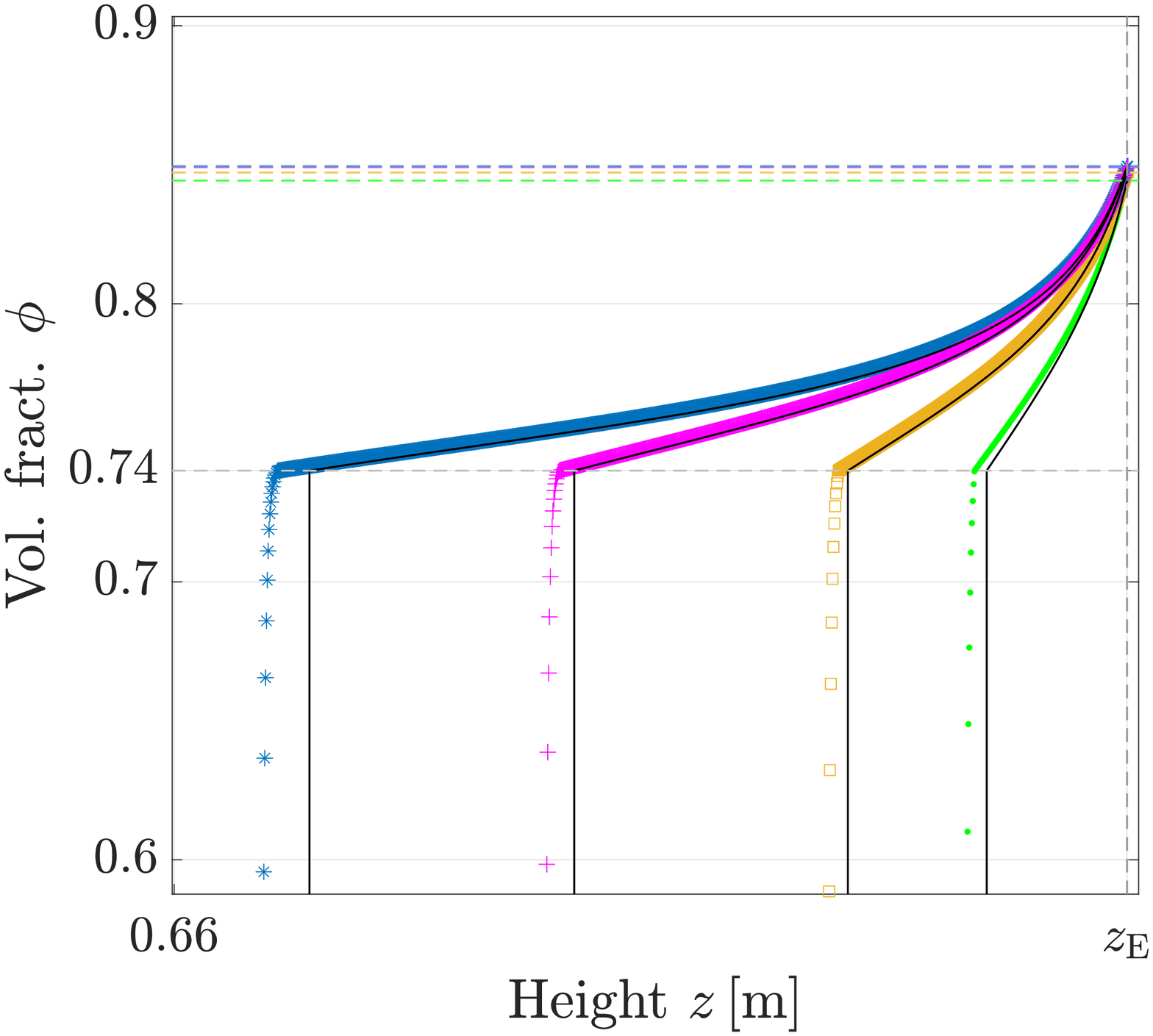} 
\end{tabular} 
\caption{Example 1: (a) Contour lines of $(Q_\mathrm{F},Q_\mathrm{U})\mapsto Z_\mathrm{fr}(\phi_\mathrm{F},Q_\mathrm{F},Q_\mathrm{U}, Q_\mathrm{W})$ for $Q_\mathrm{W}=2 \times  10^{-6} \,\mathrm{m}^3/\mathrm{s}$, $\psi_\mathrm{F}=0.2$ and $\phi_\mathrm{F}=0.3$. (b)~Approximate volume fraction of solids~$\psi$ computed with $N=3200$. 
(c) Approximate solution (dots) versus exact solution (solid lines) of volume fraction of aggregates~$\phi$ corresponding to the four point in plot~(a) computed with $N=3200$. (d)~Enlarged view of (c). \label{fig:ContourEx1}}
\end{figure}%

\begin{figure}[t] 
\centering 
\begin{tabular}{cc}  
(a) & (b) \\ 
\includegraphics[width=0.45\textwidth]{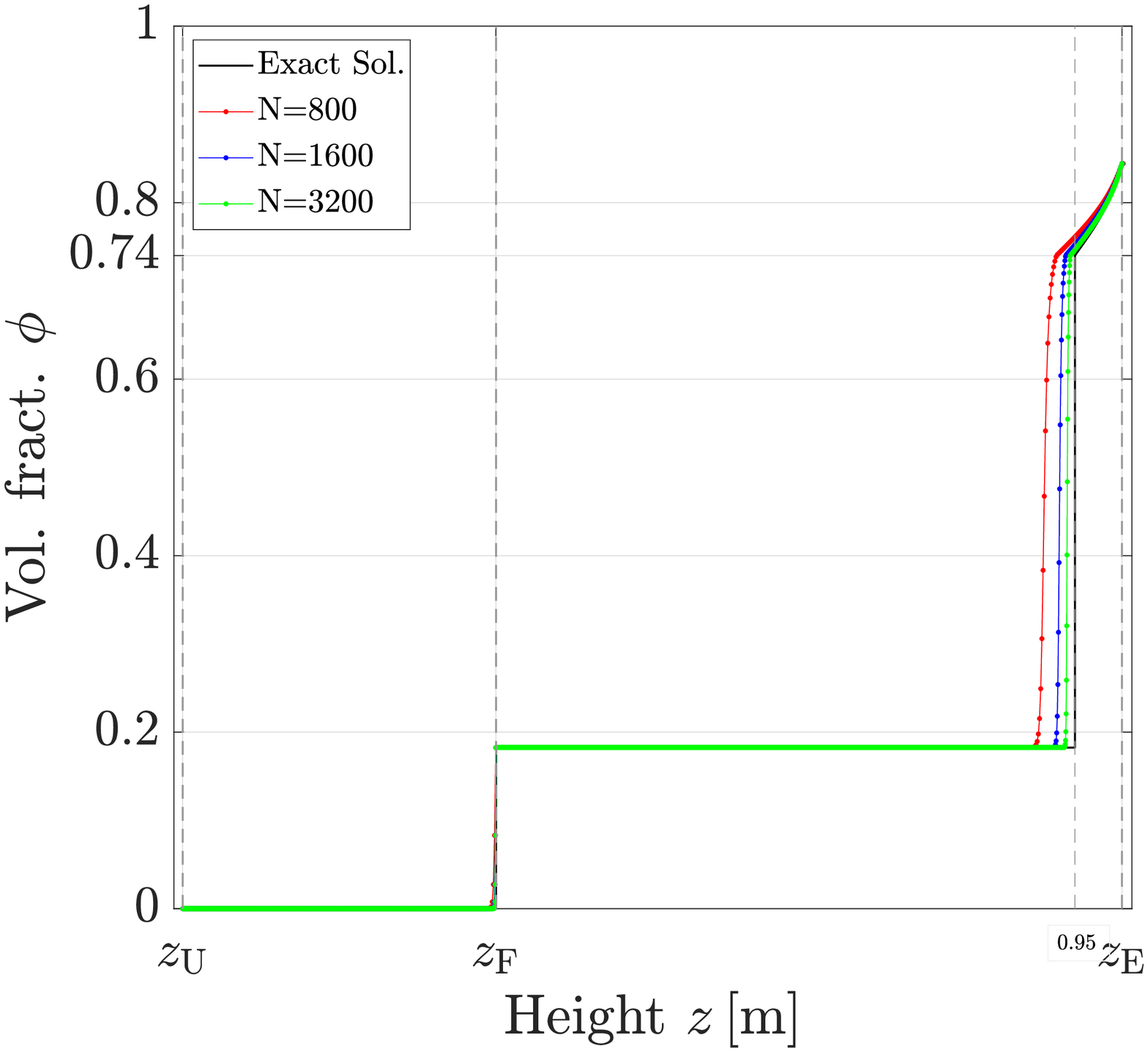} &  \includegraphics[width=0.45\textwidth]{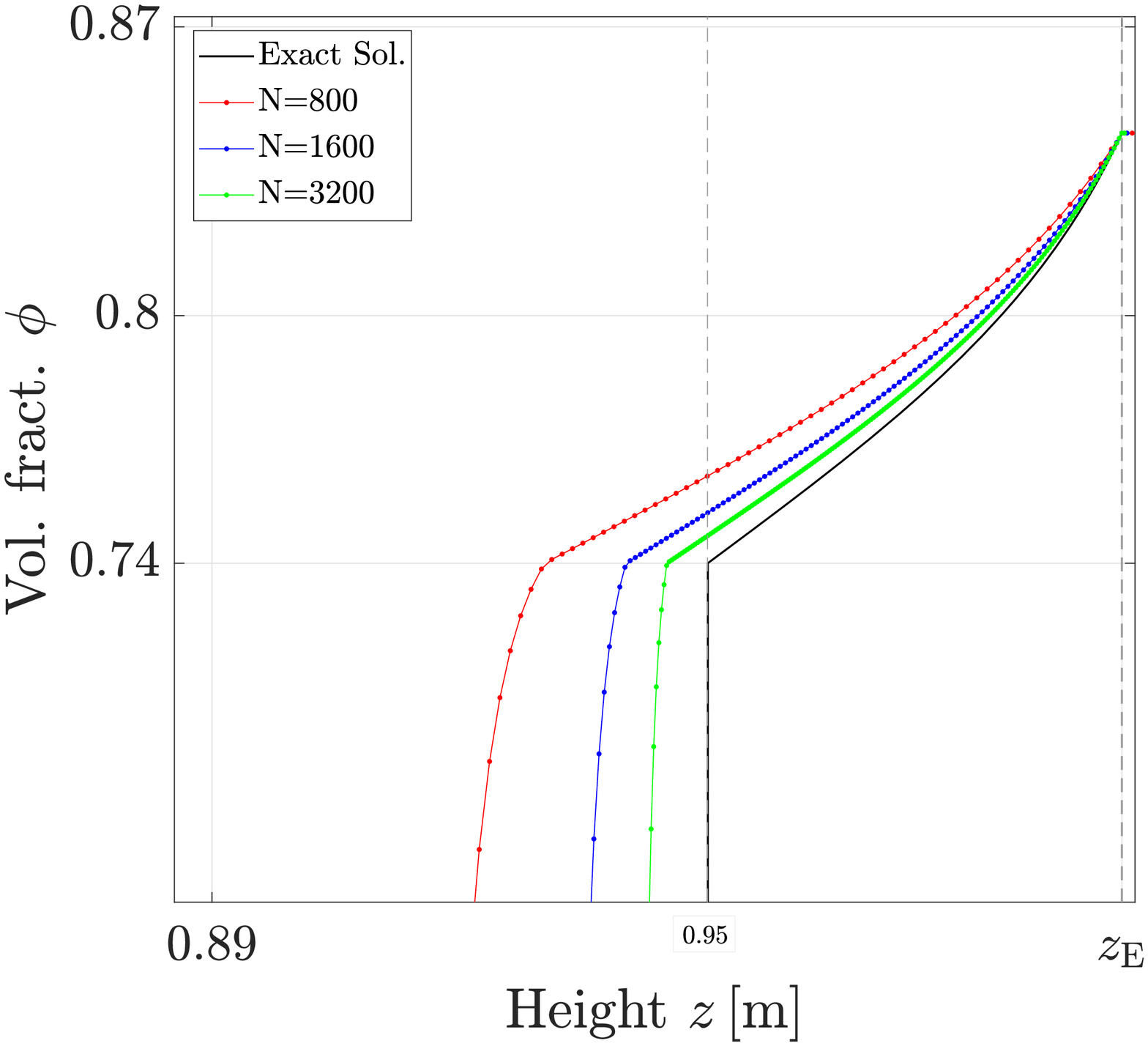} 
\end{tabular} 
\caption{Example 1: (a) Approximate solution for the point represented by a dot in Figure~\ref{fig:ContourEx1}(a) with various values of $N$.
 (b)~Enlarged view of~(a).  \label{fig:SolEx1Approx}}
\end{figure}%

We show steady-state solutions for fixed $Q_\mathrm{F}=8.9927\times10^{-5}\,\mathrm{m}^3/\mathrm{s}$ and $Q_\mathrm{W}=2.0\times10^{-6}\,\mathrm{m}^3/\mathrm{s}$ for various  values of~$Q_\mathrm{U}$; see Figure~\ref{fig:ContourEx1}(a).
For these values and with the feed volume fractions $\phi_\mathrm{F}=0.3$ and $\psi_\mathrm{F}=0.2$, we solve the ODE~\eqref{eq:ODE2} to obtain `exact'   solutions (i.e.,  the ODE is solved numerically), and the value of $z_\mathrm{fr}$ for each point; see the solid lines in Figure~\ref{fig:ContourEx1}(c) and (d).
The dots in the same plots show the numerical solutions, which are obtained by simulating a long time from any initial data.
All solutions have the same volume fraction~$\phi_\mathrm{E}$ at the top, since this is given by the explicit formula~\eqref{eq:phiE}.
Figure~\ref{fig:ContourEx1}(b) shows the steady state for the solids with particles only below the feed level.

A clear difference between the two types of solution of~$\phi$ can be seen near the discontinuity.
This is an inaccuracy of the numerical solution, which seems to converge to the exact one as $N\to\infty$; see
Figure~\ref{fig:SolEx1Approx}, which shows the steady-state solution for the solid point in Figure~\ref{fig:ContourEx1}(a) for  various  values of~$N$.

\begin{figure}[t] 
\centering 
\includegraphics[width=0.45\textwidth]{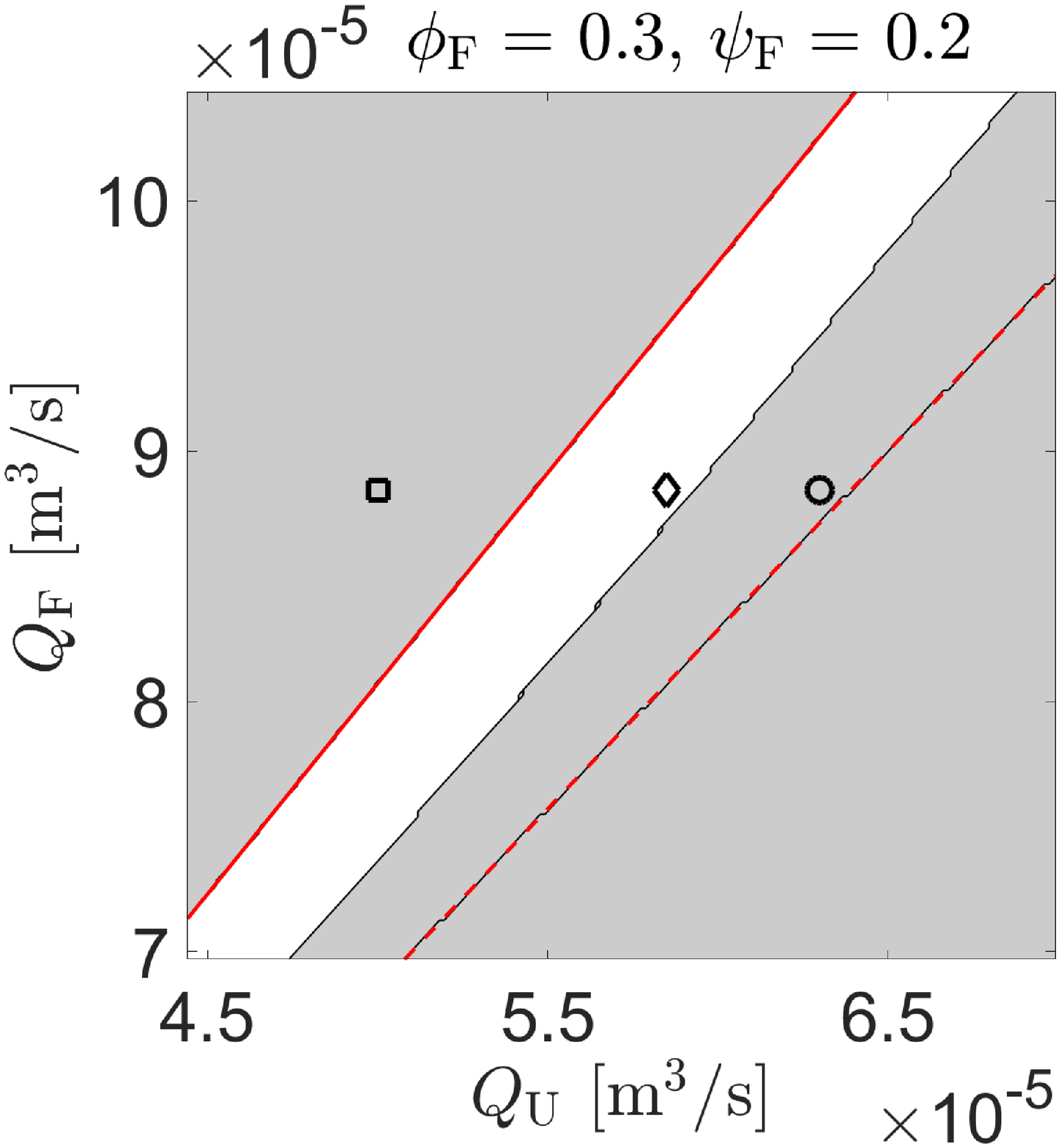} 
\caption{Example 2: An operating charts for $\phi_\mathrm{F}=0.3$ and $\psi_\mathrm{F}=0.2$. The point $(Q_\mathrm{U},Q_\mathrm{F})=(5.85,8.846)\times10^{-5}\, \mathrm{m}^3/\mathrm{s}$ marked with a diamond in the white region results in a desired steady state with a froth layer at the top of the column.
The points marked with a square $(Q_\mathrm{U},Q_\mathrm{F})=(5.0,8.846)\times10^{-5}\, \mathrm{m}^3/\mathrm{s}$ and a circle $(Q_\mathrm{U},Q_\mathrm{F})=(6.3,8.84)\times10^{-5}\, \mathrm{m}^3/\mathrm{s}$ result in no froth (Figure~\ref{fig:SolEx2a}) or a tank full of froth (Figure~\ref{fig:SolEx2b}), respectively. 
(The plot is a zoom of Figure~\ref{fig:opchart1} (b) and the black curves are smoother than they here appear due to numerical resolution.)
 \label{fig:OpChartEx2}}
\end{figure}%

\begin{figure}[t!] 
\centering 
\begin{tabular}{cc}  
(a) & (b) \\ 
\includegraphics[width=0.45\textwidth]{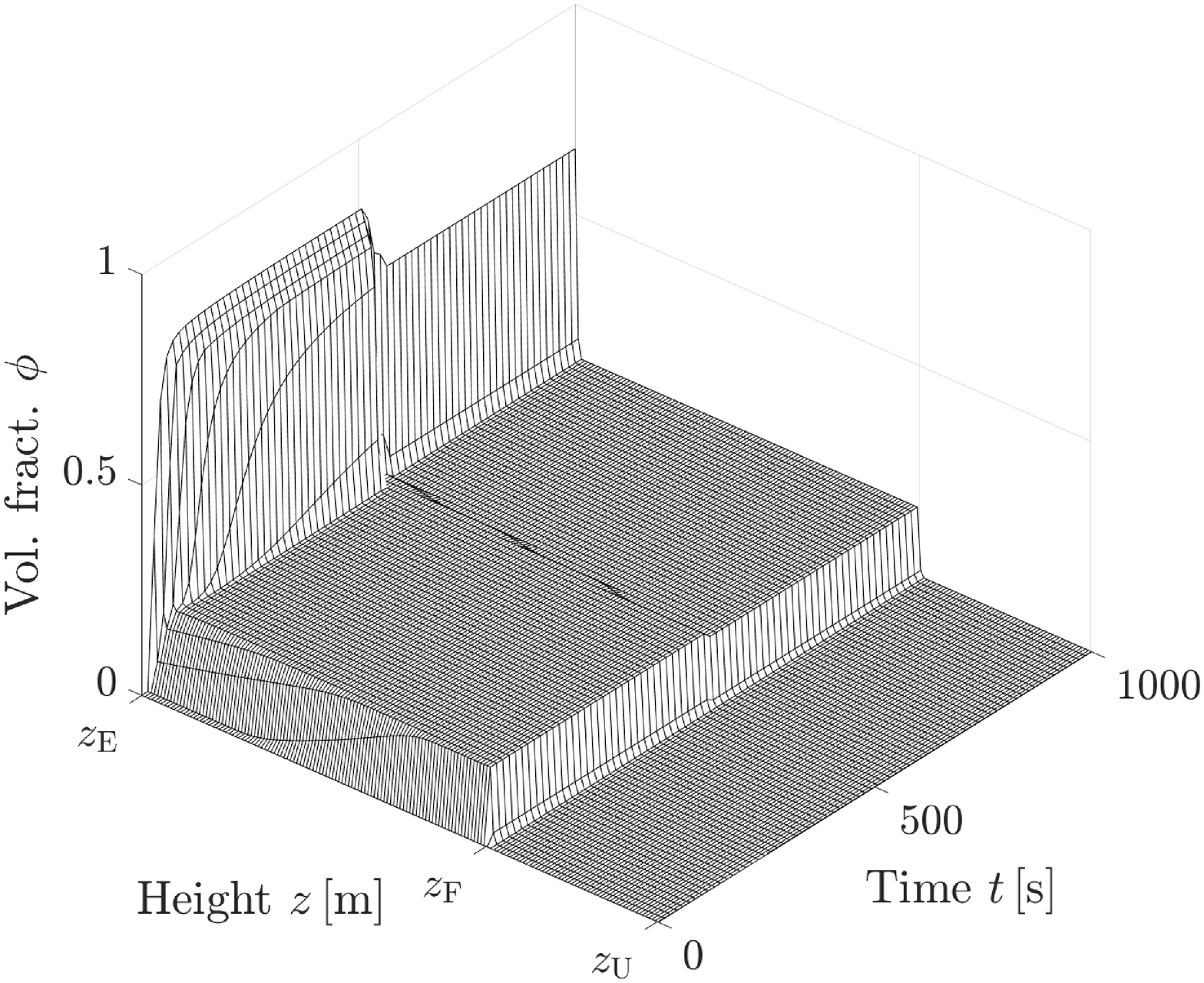} &  \includegraphics[width=0.45\textwidth]{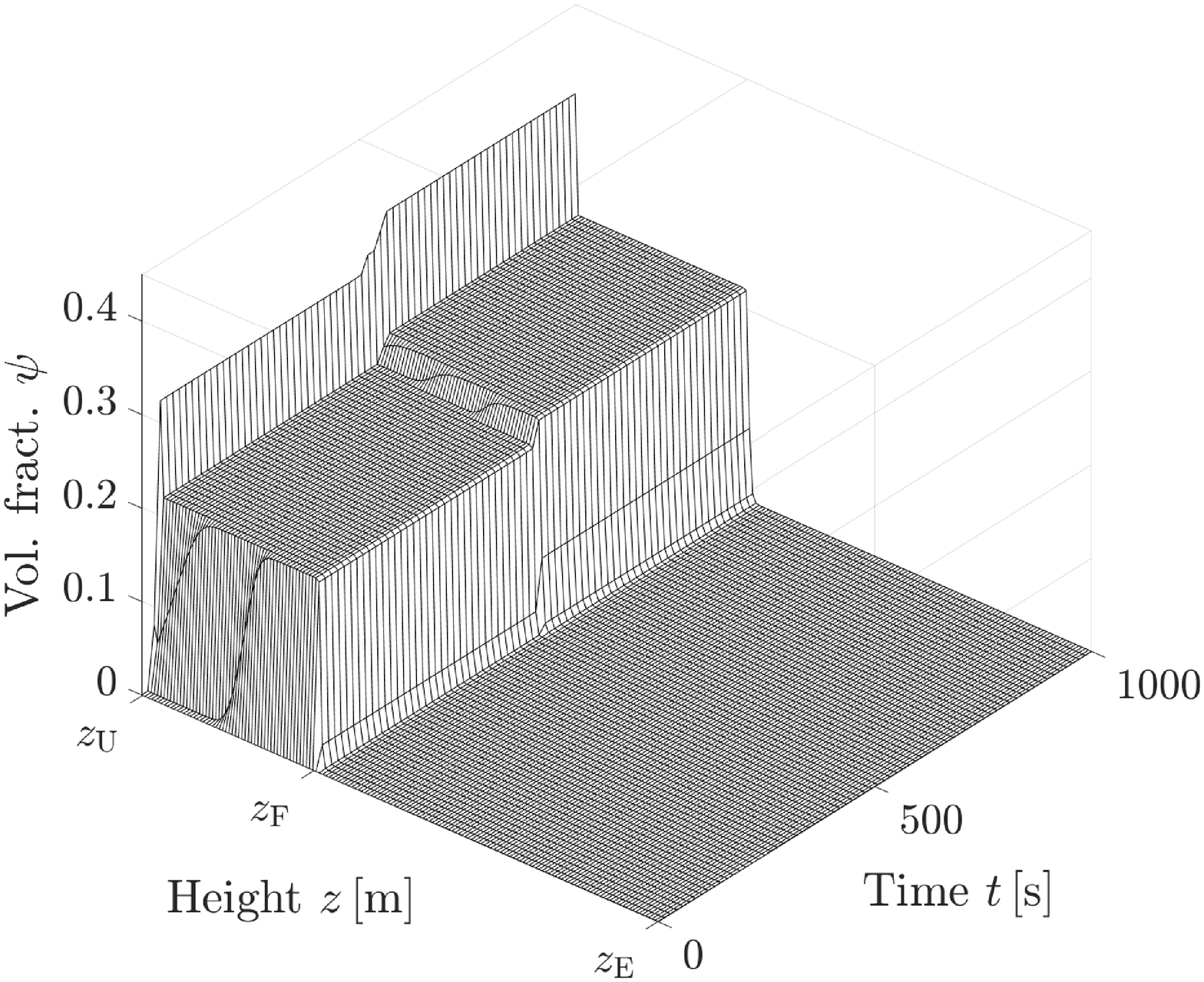} 
\end{tabular} 
\caption{Example 2: Simulation with $N=1600$ of the volume fractions of  (a) aggregates $\phi$  and (b) 
 solids $\psi$   from a tank filled of only water.
The initial operating point $(Q_\mathrm{U},Q_\mathrm{F})=(5.85,8.846)\times10^{-5}\, \mathrm{m}^3/\mathrm{s}$ (diamond in Figure~\ref{fig:OpChartEx2}) is at $t=500\,$s changed to $(5.0,8.846)\times10^{-5}\, \mathrm{m}^3/\mathrm{s}$ (square in Figure~\ref{fig:OpChartEx2}).\label{fig:SolEx2a}}
\end{figure}%

\begin{figure}[t] 
\centering 
\begin{tabular}{cc}  
(a) & (b) \\ 
\includegraphics[width=0.45\textwidth]{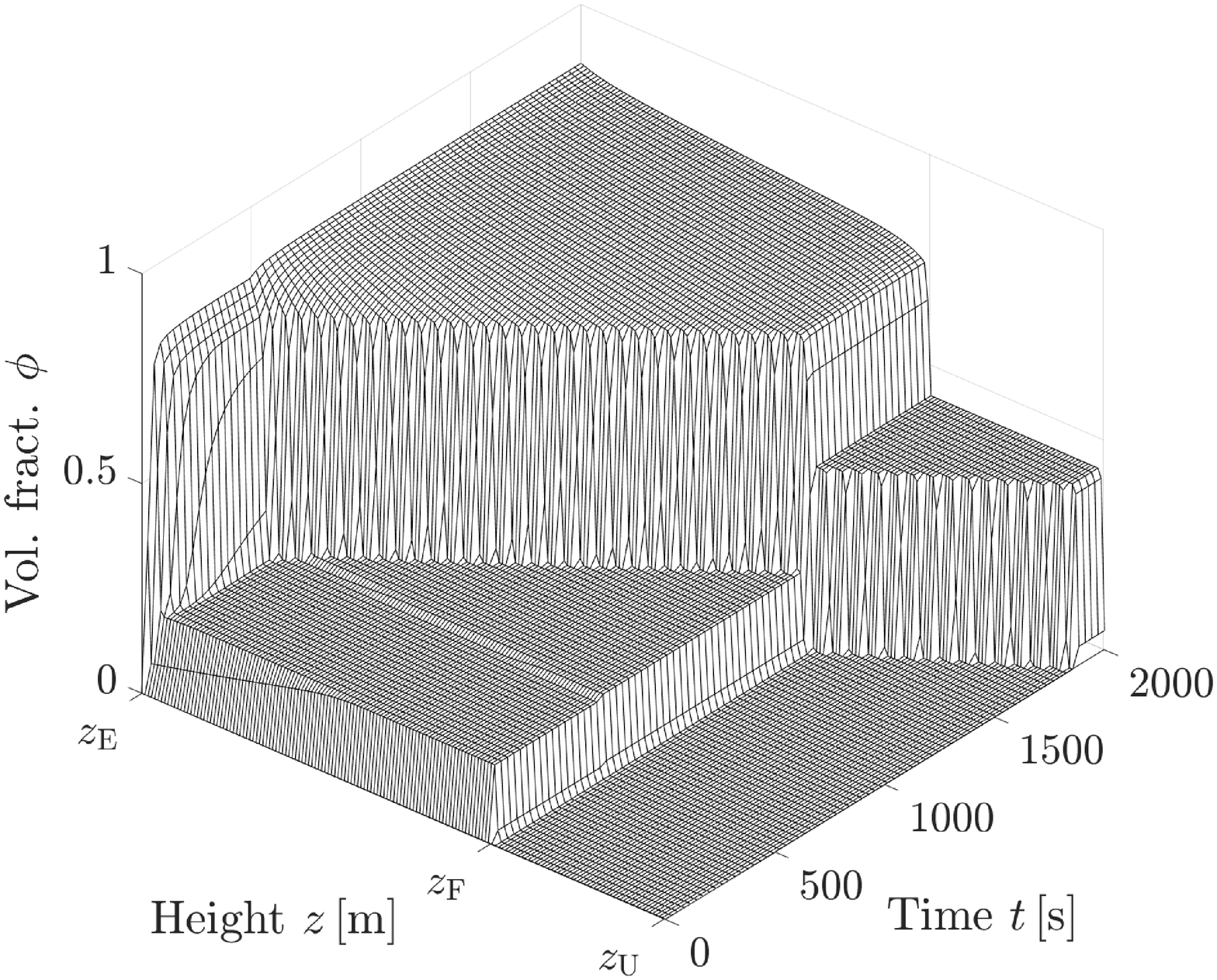} & 
  \includegraphics[width=0.45\textwidth]{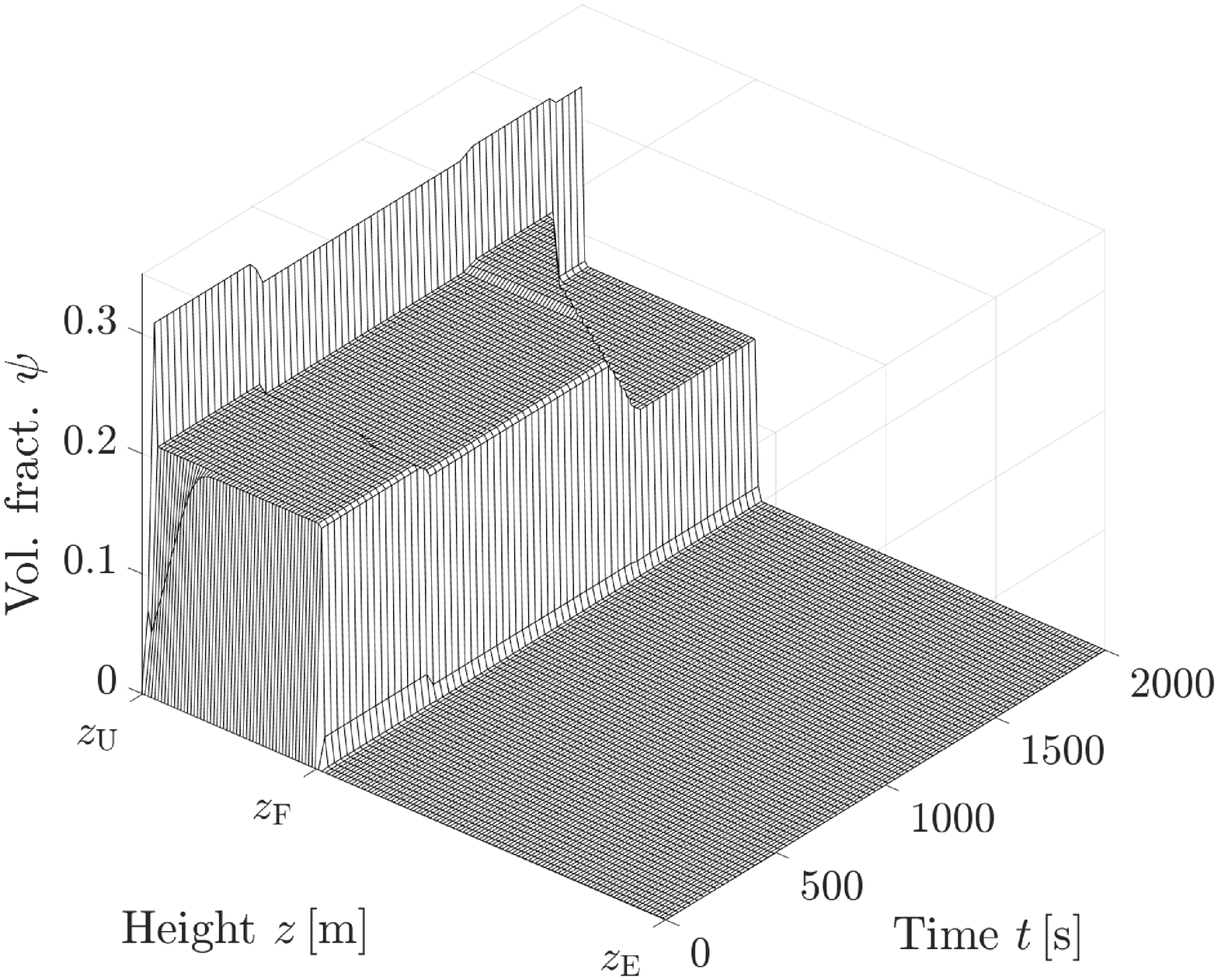} 
\end{tabular} 
\caption{Example 2: Simulation with $N=1600$ of the volume fractions of  (a) aggregates $\phi$  and (b) solids $\psi$   from a tank filled of only water.
The initial operating point $(Q_\mathrm{U},Q_\mathrm{F})=(5.85,8.846)\times10^{-5}\, \mathrm{m}^3/\mathrm{s}$ (diamond in Figure~\ref{fig:OpChartEx2}) is at $t=500\,$s changed to  $(6.3,8.846)\times10^{-5}\, \mathrm{m}^3/\mathrm{s}$ (circle in Figure~\ref{fig:OpChartEx2}).\label{fig:SolEx2b}}
\end{figure}%

\subsection{Example 2}

We start from a tank filled with only water at time $t=0 \,\mathrm{s}$, i.e., $\phi(z,0)=\psi(z,0)=0$ for all $z$, when we start pumping aggregates, solids, fluid and wash water with $\phi_\mathrm{F}=0.3$ and $\psi_\mathrm{F}=0.2$.
In the white region of the operating chart in Figure~\ref{fig:OpChartEx2}, we choose the point (diamond symbol) $(Q_\mathrm{U},Q_\mathrm{F})=(5.85,8.846)\times10^{-5}\, \mathrm{m}^3/\mathrm{s}$.
The wash water volumetric flow is $Q_\mathrm{W}=2.0\times 10^{-6}\, \mathrm{m}^3/\mathrm{s}$.
Then $Q_\mathrm{E}=1.4496\times 10^{-4}\, \mathrm{m}^3/\mathrm{s}$ and one obtains a desired steady state with a thin layer of froth at the top and solids only below the feed level after about $500\,\mathrm{s}$; see Figures~\ref{fig:SolEx2a}~(a) and~\ref{fig:SolEx2b}~(a).

Once the system is in steady state at $t=500\,\mathrm{s}$, we perform two different changes corresponding to the points marked with a square (left) and a circle (right) in the operating chart in Figure~\ref{fig:OpChartEx2} with the corresponding responses seen in Figures~\ref{fig:SolEx2a} and \ref{fig:SolEx2b}, respectively. The jump from the middle point (diamond) to the left point (square) means a jump from $Q_\mathrm{U}=5.85\times10^{-5} \,\mathrm{m}^3/\mathrm{s}$ to the smaller value $5.0\times10^{-5} \,\mathrm{m}^3/\mathrm{s}$ and
produces  the solution in Figure~\ref{fig:SolEx2a}.
After $t=1000\,\mathrm{s}$, there is no froth in zone~2 and the solids volume fraction is slightly higher in the new steady state.
 
If the jump from the middle point (diamond) instead goes to the right point (circle), i.e., the new value at $t=500\,\mathrm{s}$ is the larger $Q_\mathrm{U}=6.3\times10^{-5} \,\mathrm{m}^3/\mathrm{s}$, Figure~\ref{fig:SolEx2b} shows the reaction of the system until $t=2000\,\mathrm{s}$.
The aggregates fill the entire column while the solids volume fraction has a lower value in the new steady state.
We have demonstrated that operating points outside the white region lead to non-desired steady states.

\subsection{Example 3} 
\begin{figure}[t!] 
\centering 
\begin{tabular}{cc}  
(a) & (b) \\
\includegraphics[width=0.43\textwidth]{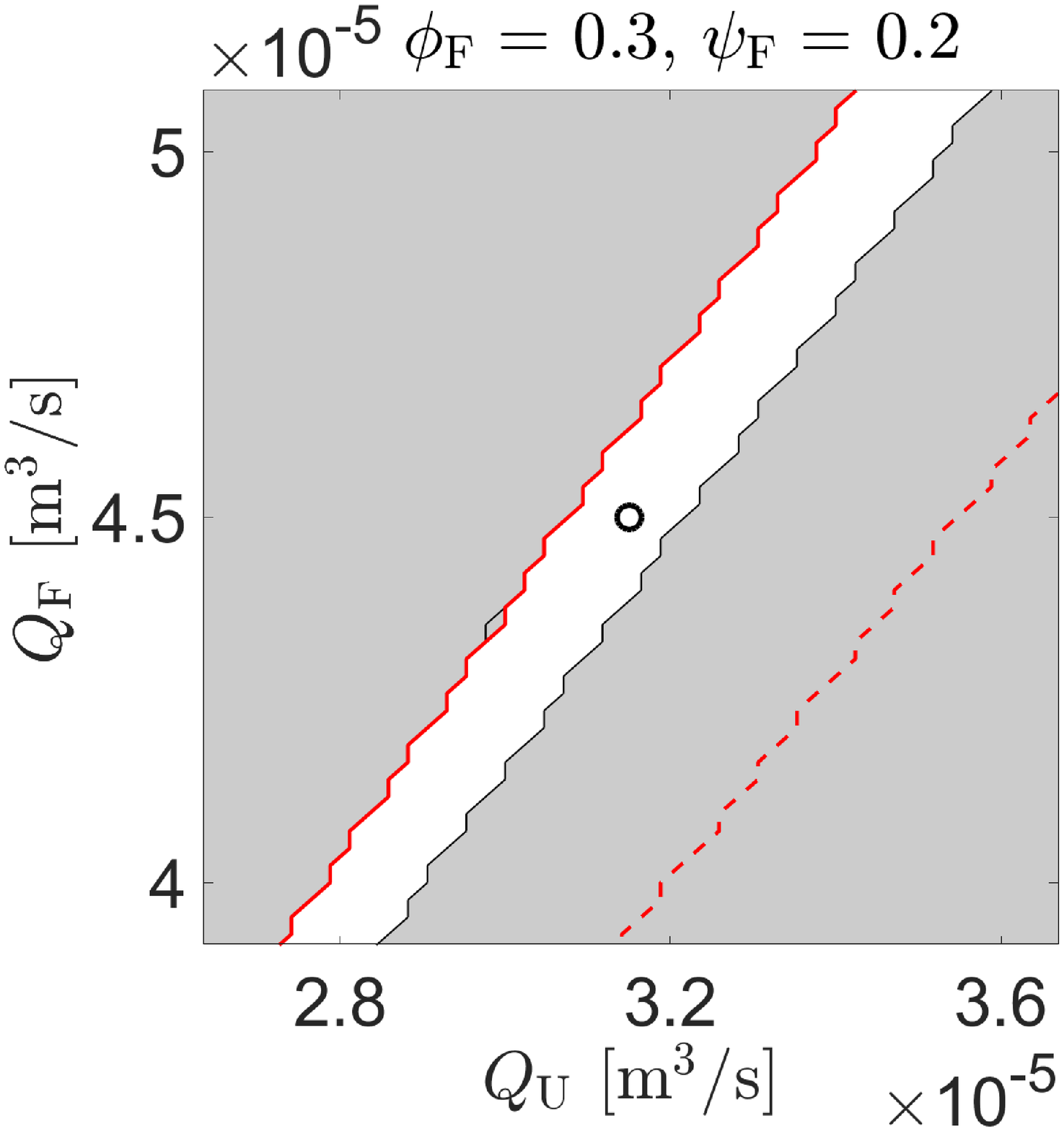}  & \includegraphics[width=0.43\textwidth]{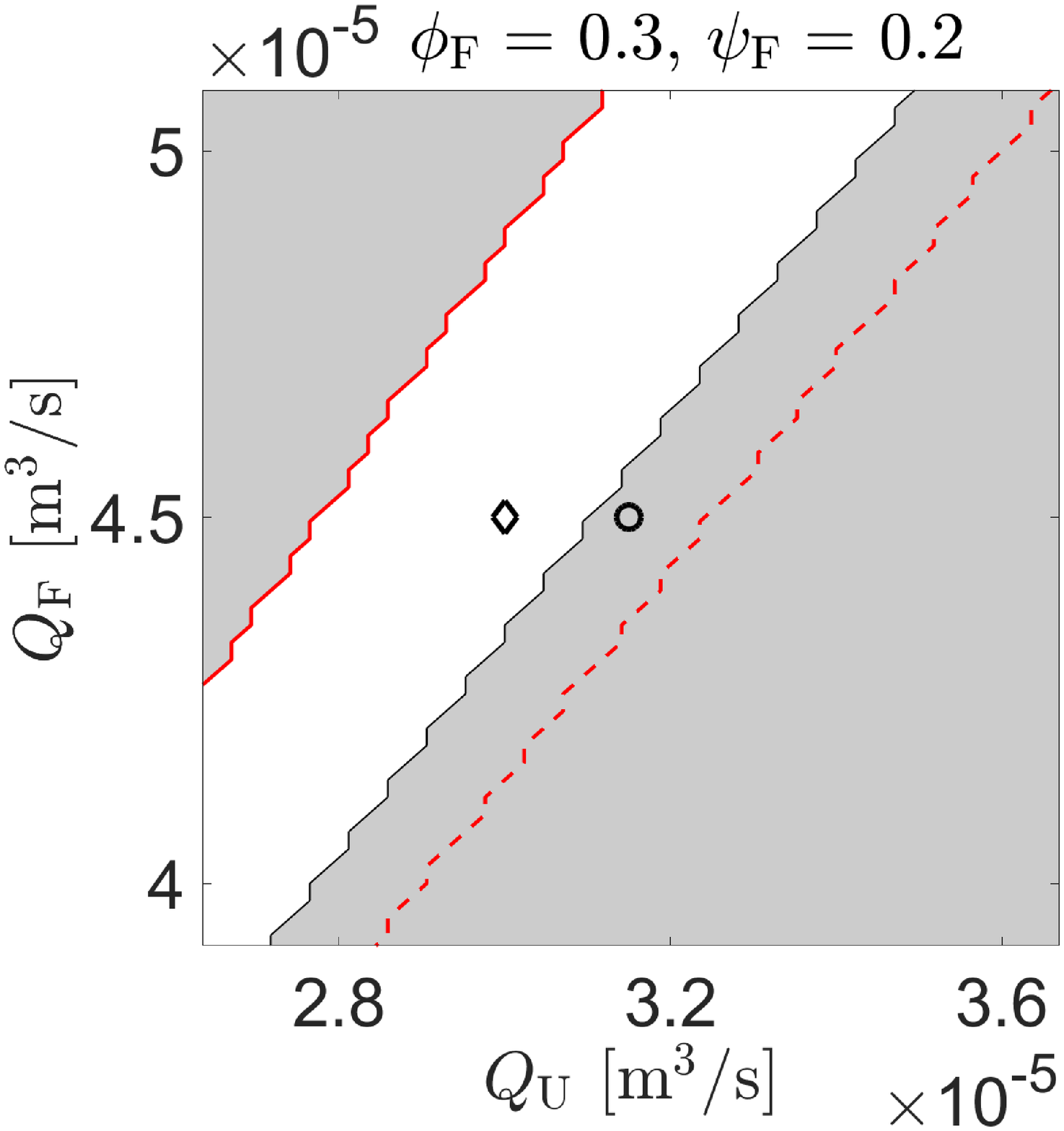} 
\end{tabular} 
\caption{Example 3. Operating charts for $\phi_\mathrm{F}=0.3$ and $\psi_\mathrm{F}=0.2$ with (a) $Q_\mathrm{W}=3.15\times10^{-5}\,\mathrm{m}^3/\mathrm{s}$, (b)   $Q_\mathrm{W}=3.0 \times10^{-5}\,\mathrm{m}^3/\mathrm{s}$.
The initial point $(Q_\mathrm{U},Q_\mathrm{F})=(3.15,4.5)\times10^{-5}\, \mathrm{m}^3/\mathrm{s}$ is marked with a circle and the one after the control action $(Q_\mathrm{U},Q_\mathrm{F})=(3.0,4.5)\times10^{-5}\, \mathrm{m}^3/\mathrm{s}$ with a diamod.
(The curves are smoother than they appear here due to numerical resolution.)
\label{fig:OpChartEx3}}
\end{figure}

\begin{figure}[t!] 
\centering 
\begin{tabular}{cc}  (a) & (b) \\ 
\includegraphics[width=0.45\textwidth]{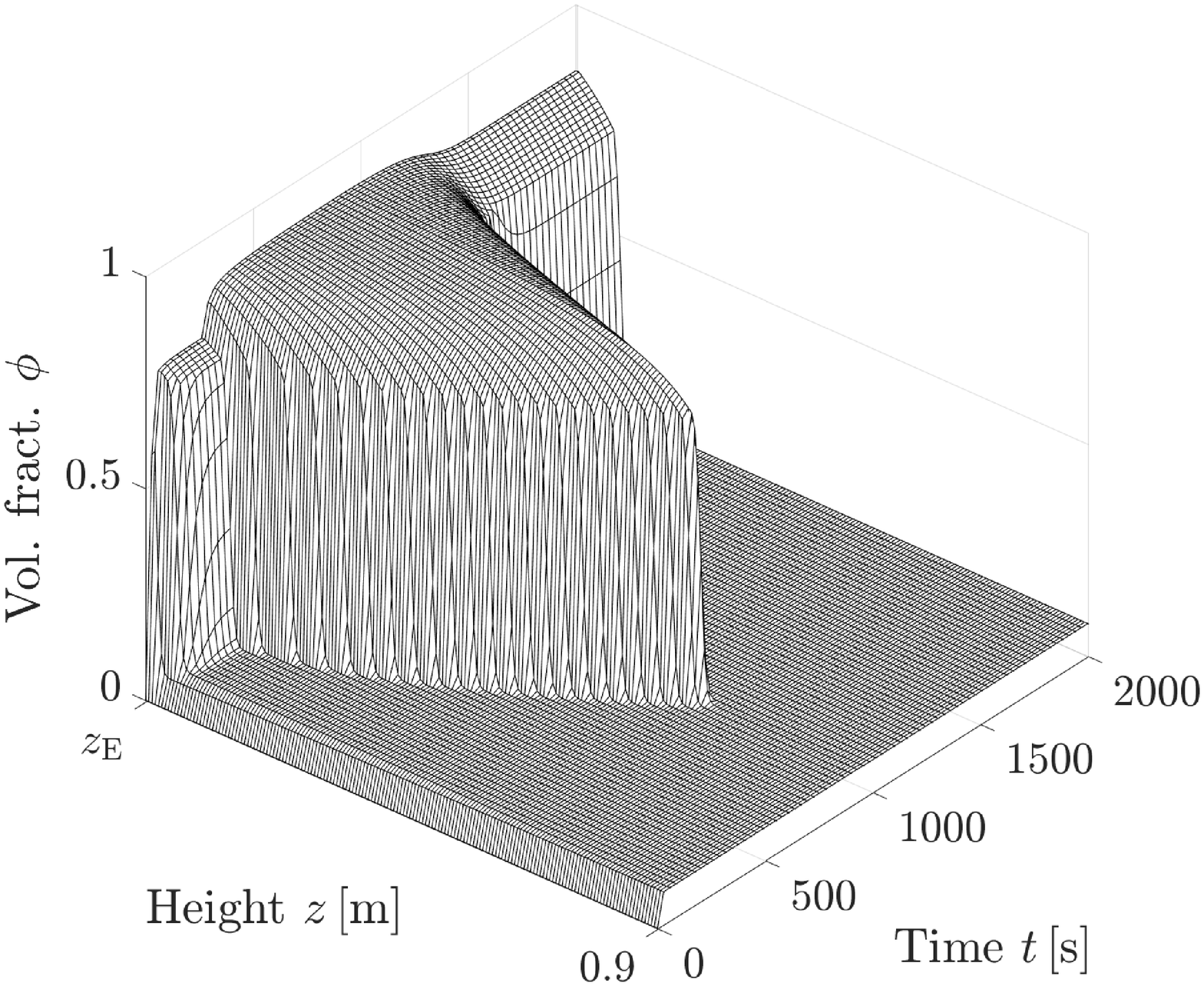} & \includegraphics[width=0.45\textwidth]{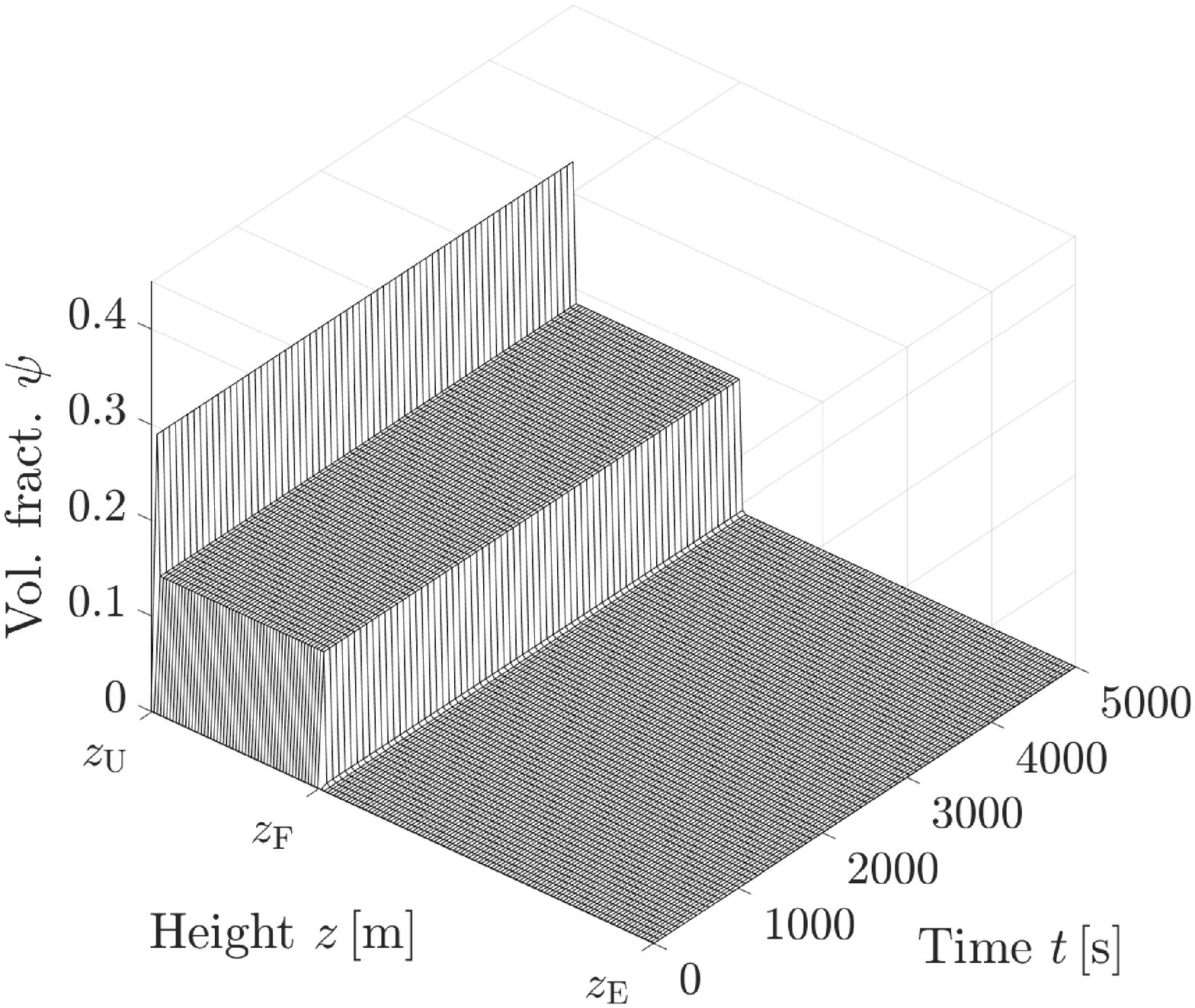} \\
(c) & (d) \\ 
\includegraphics[width=0.45\textwidth]{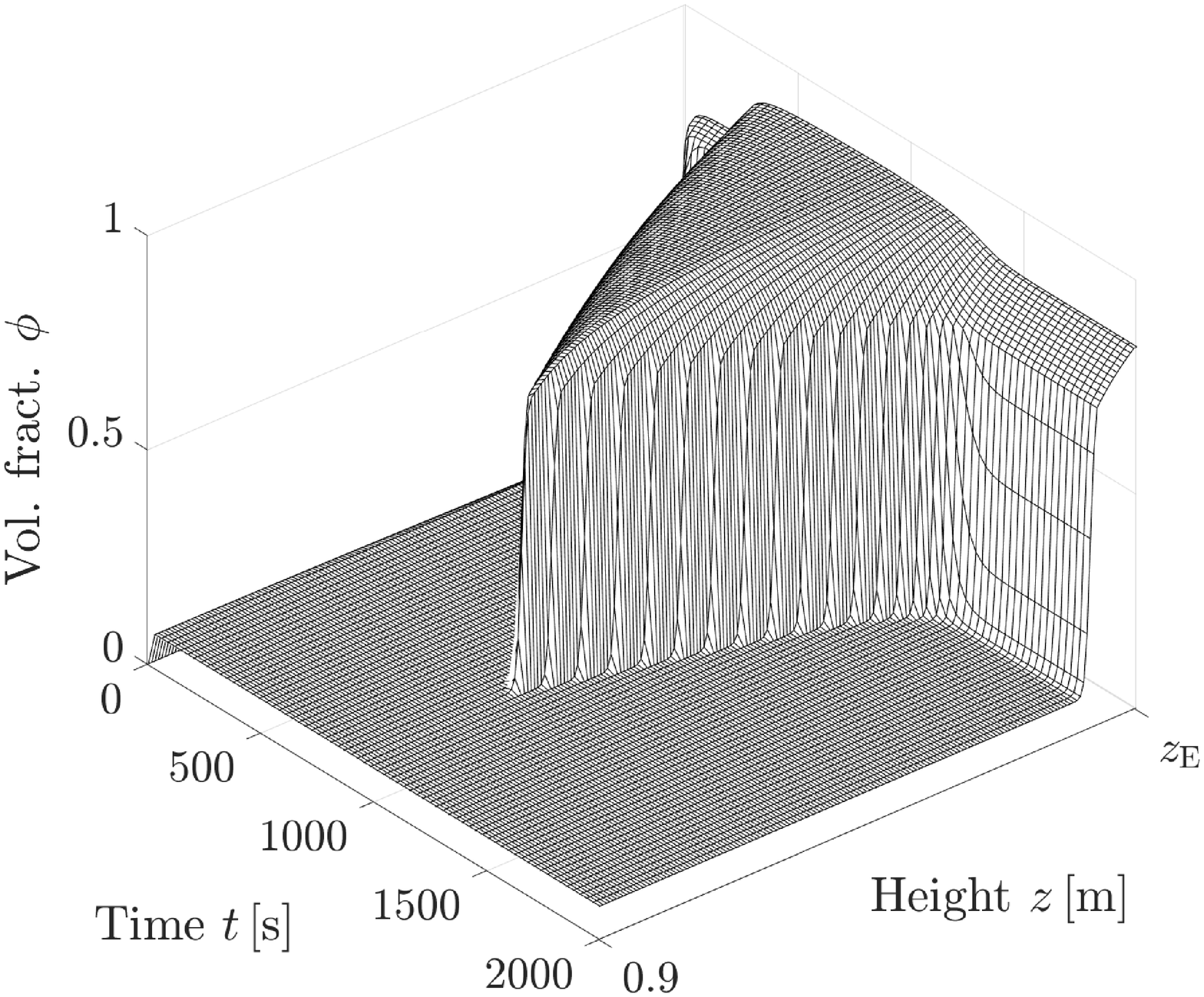} & \includegraphics[width=0.45\textwidth]{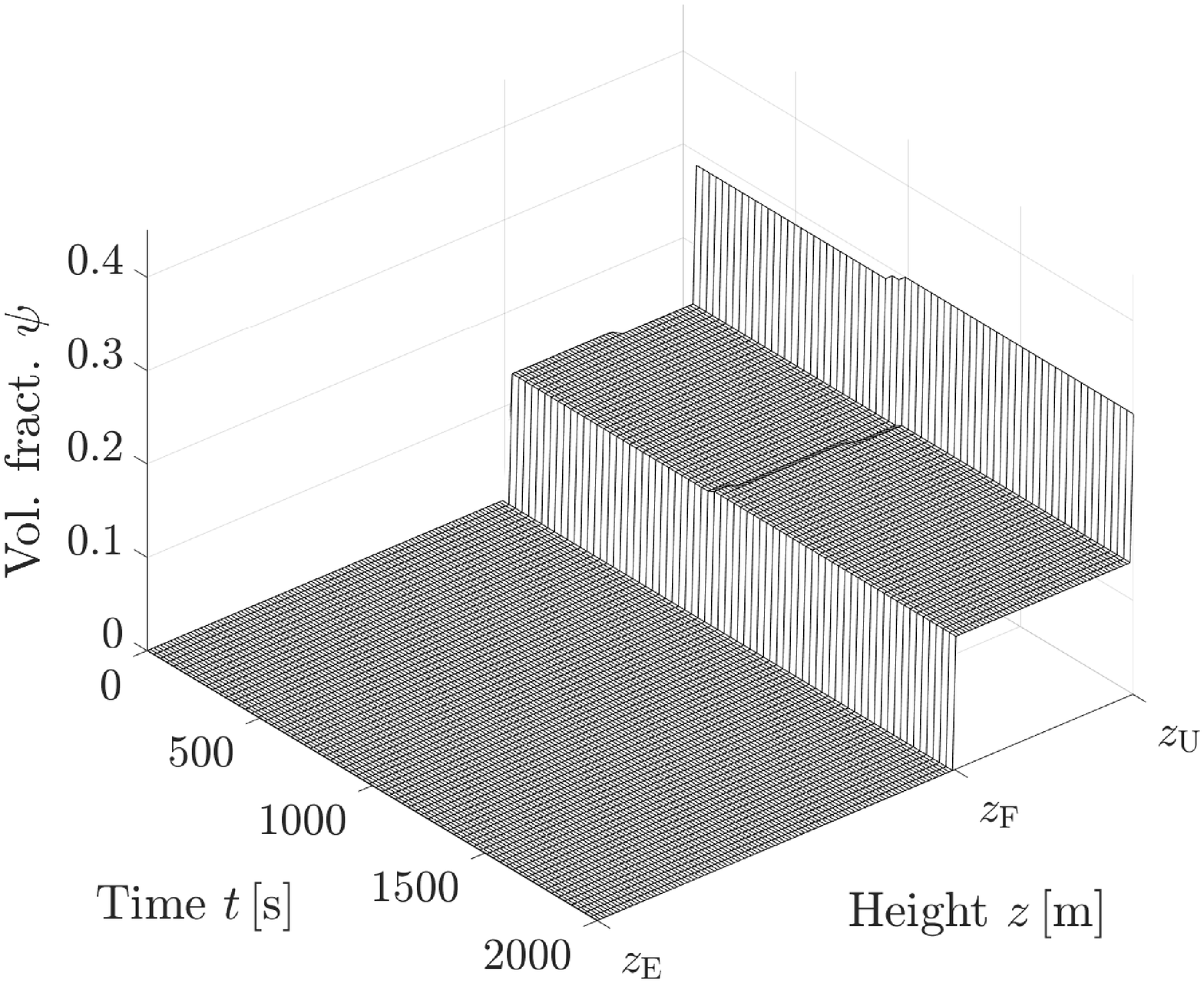}
\end{tabular} 
\caption{Example 3: Time evolution of the volume fraction (a, c) of aggregates $\phi$ and (b,d) solids $\psi$  computed with $N=1600$ and seen from two different angles.
A step change down in $Q_\mathrm{W}$ occurs at $t=300\,\mathrm{s}$ and a control action by decreasing $Q_\mathrm{U}$ is made at $t=1000\,\mathrm{s}$. \label{fig:SolEx3}}
\end{figure}

Again, the tank is filled with only water at time $t=0 \,\mathrm{s}$ when we start feeding it with $\phi_\mathrm{F}=0.3$ and $\psi_\mathrm{F}=0.2$. 
The wash water flow is $Q_\mathrm{W}=4.0\times10^{-6}\,\mathrm{m}^3/\mathrm{s}$ and hence the effluent volumetric flow is $Q_\mathrm{E}=1.75\times10^{-5}\,\mathrm{m}^3/\mathrm{s}$. 
From the corresponding operating chart in Figure~\ref{fig:OpChartEx3}~(a), we choose the point of volumetric flows $(Q_\mathrm{U},Q_\mathrm{F})=(3.15,4.5)\times10^{-5}\, \mathrm{m}^3/\mathrm{s}$ lying in the white region.
Then a desired steady state builds up quickly and at $t=250 \,\mathrm{s}$ there is a thin froth layer at the top of in zone~2 and with solids only in zone~1; see Figure~\ref{fig:SolEx3}.

Once the system is in steady state, we change at $t=300 \,\mathrm{s}$ the volumetric flow of the wash water from $Q_\mathrm{W}=4.0\times10^{-6}\,\mathrm{m}^3/\mathrm{s}$ to $1.0\times10^{-6}\,\mathrm{m}^3/\mathrm{s}$ and simulate the reaction of the system.
In the corresponding operating chart for this new set of variables, the point $(Q_\mathrm{U}, Q_\mathrm{F}) =(3.15,4.5)\times10^{-5}\,\mathrm{m}^3/\mathrm{s}$ is no longer in the white region; see Figure~\ref{fig:OpChartEx3} (b, circle point), and no desired steady state is feasible. 
As it can be seen in Figure~\ref{fig:SolEx3} (a), with less flow of wash water flushing the aggregates out at the top, the froth layer increases downwards. 
At time $t=1000 \,\mathrm{s}$, we make a control action and change the volumetric flow from  $Q_\mathrm{U}=3.15\times10^{-5}\,\mathrm{m}^3/\mathrm{s}$ to $3.0 \times10^{-5}\,\mathrm{m}^3/\mathrm{s}$ so that the new point lies inside the white region of the corresponding operating chart in Figure~\ref{fig:OpChartEx3}~(b, diamond point).
Figures~\ref{fig:SolEx3}~(a) and~(c) show that a second desired steady state is reached after $t=1500 \,\mathrm{s}$.
Figures~\ref{fig:SolEx3}~(b) and~(d) show that the solids settle in any case.

\section{Conclusions}\label{sec:conclusions}

Our previous one-dimensional model of a flotation column, where the movement of rising aggregates and settling solids follow the drift- and solids-flux theories, is a triangular hyperbolic $2\times 2$ system of nonlinear PDEs of the first order.
Here, we propose an extended model where the drainage of liquid in the froth layer due to capillarity is included.
The traditional derivation of the drainage PDE, valid only within the froth, is combined with further experimental findings from the literature to end up in a constitutive relationship between the relative velocity of aggregates to fluid (or suspension of hydrophilic solids), which in the governing equations yields a second-order-derivative degenerate nonlinear term.

An analysis of the possible steady states with a froth layer at the top of the column (desired steady states) leads to several inequalities involving the feed input variables and other control volumetric flows; see Theorem~\ref{theorem1}.
Those inequalities are visualized in operating charts; see Figure~\ref{fig:opchart1}, in which the white region shows the necessary location of an operating point $(Q_\mathrm{U},Q_\mathrm{F})$ for having a desired steady state after a time of transient behaviour.

With parameters extracted from the literature, the white region of an operating chart is quite small, meaning that the existence of a froth layer is very sensitive to small changes in any of the control variables $Q_\mathrm{U}$ and $Q_\mathrm{W}$.
Different operating points $(Q_\mathrm{U},Q_\mathrm{F})$ in the white region give rise to different thicknesses of the froth layer.
Unfortunately, our model anticipates a very sensitive dependence of the pulp-froth interface location~$z_\mathrm{fr}$ on the operating point; see Figure~\ref{fig:opchart2}, where the yellow surface shows that the most common values of~$z_\mathrm{fr}$ is close to one, meaning a thin froth layer.
The surfaces seen in plots~(c) and (d) indicate a very large gradient from $z_\mathrm{fr}$ just below one down to $z_\mathrm{fr}=0.33=z_\mathrm{F}$.
(Even finer resolutions indicate that the graph is continuous.)
The numerical scheme suggested resolves discontinuities well and  numerical results
 (e.g., those of Figure~\ref{fig:SolEx2b}) show that the volume fractions, and their sum, stay between zero and one, as is proven in Theorems~\ref{thm:bounded_phiR} and~\ref{thm:bounded_varphiR}.

Overall, the steady-state analysis, boundedness properties of the numerical solutions, and simulation results indicate that 
 the model is useful for the simulation of flotation columns and could be used, for example, to simulate the effect of various 
  alternative control actions. In light of this practical interest it would be desirable to obtain a well-posedness 
   (existence and uniqueness) result for the underlying model. The first equation of \eqref{eq:gov} 
    (the one for~$\phi$) as a scalar strongly degenerate parabolic equation with discontinuous flux that 
     is independent of~$\psi$ can be handled by known arguments (cf., e.g., \cite{Karlsen&R&T2002,Karlsen&R&T2003,Burger&K&T2005a}). The corresponding model for $D \equiv 0$ 
      (without capillarity) is a triangular system of conservation laws for which 
       convergence results of monotone schemes to a weak solution are available, at least  for 
        the case of fluxes without spatial discontinuity \citep{Karlsen2008_conv_triangular,CMR10}. It is 
         not clear at the moment whether the corresponding arguments, based on the compensated compactness 
          method, can also be applied to the system \eqref{eq:gov}. Therefore, well-posedness of the 
           model is at the moment left as an open problem.

The model of a flotation column with drainage can certainly be extended to include additional processes. 
 For instance, 
 one could incorporate the possibility of bursting bubbles at the top by assuming that the flotation column is constructed in a way such that a portion~$\alpha$ of the froth overflows with unbursting aggregates, whereas for  the portion~$1-\alpha$, the aggregates burst.  
The latter means that the gas `disappears' (i.e., is released into the surrounding air), whereas the suspension and the hydrophobic particles attached to the bursting bubbles follow the effluent stream. There is practical interest in quantifying this effect \citep{Neethling2018}.
Hence, the effluent volume fraction of aggregates would then be $\phi_\mathrm{E}:=\alpha \phi_\mathrm{E}^+$, and the solids $\varphi_\mathrm{E}:=\varphi_\mathrm{E}^+$.
Under the present assumptions, the factor $\alpha$ does not influence the solution inside the column.
It could however depend on the wash water flow, and this could be an extension of the present model.
Another thinkable extension  could consist in the explicit description of the aggregation process itself. 
 A common    variant of the column drawn in Figure~\ref{fig:Column} has gas feed and pulp feed at different levels, 
  which  
  form a so-called collection zone where the attachment of hydrophobic particles takes place.  Such a description 
   would compel a distinction between hydrophobic and hydrophilic particles as distinct solid phases.

\section*{Acknowledgements}

R.B.\ acknowledges  support from  ANID (Chile) through Fondecyt project 1210610; Anillo project ANID/PIA/210030; Centro de Modelamiento Matem\'{a}tico (CMM), projects 
  ACE210010 and FB210005 of BASAL funds for Centers of Excellence; and CRHIAM, project ANID/FONDAP/15130015. 
 S.D.\ acknowledges support from the Swedish Research Council (Vetenskapsr\aa det, 2019-04601). M.C.M.\ is supported by grant MTM2017-83942 funded by Spanish MINECO and by grant PID2020-117211GB-I00 funded by MCIN/AEI/10.13039/501100011033. Y.V.\ is supported by SENACYT (Panama).

\bibliographystyle{apalike}
\bibliography{ref_RB}

\appendix

\section{Proofs of boundedness of numerical solutions}

We outline the proofs of Theorems~\ref{thm:bounded_phiR} and~\ref{thm:bounded_varphiR}, which are both based on 
 monotonicity arguments. The calculations are straightforward in both cases; details for the case $D \equiv 0$ 
  are provided by~\cite{bdmv_nhm22}. 

\subsection*{Outline of the proof of Theorem~\ref{thm:bounded_phiR}}

 Assume that $\smash{\phi_{i-1/2}^n}$ 
 and   $\smash{\tilde{\phi}_{i-1/2}^n}$, $i \in \mathbb{Z}$, $n=0,1,2, \dots$ are two numerical solutions 
   produced by the  numerical scheme~\eqref{eq:phi_updateR}. Then monotonicity means that if 
       $\smash{\phi_{i-1/2}^n \leq \tilde{\phi}_{i-1/2}^n}$ for all~$i$,  
     then  $\smash{\phi_{i-1/2}^{n+1} \leq \tilde{\phi}_{i-1/2}^{n+1}}$ for all~$i$, for all $n=0,1,2, \dots$. 
      For the case of a three-point scheme such as~\eqref{eq:phi_updateR} this property 
       can be verified 
  by showing that $\smash{\partial\phi_{i-1/2}^{n+1}/\partial\phi_{k-1/2}^{n}\geq 0}$ for all~$i$ and $k=i-1,i,i+1$.
 In fact, we have 
\begin{align*} 
\frac{\partial\phi_{i-1/2}^{n+1}}{\partial\phi_{i-3/2}^{n}}&=
\frac{\lambda}{A_{i-1/2}}\left(Q_{i-1}^{n+}  + (A\gamma)_{i-1}\big(\tilde{v}(\phi_{i-1/2}^{n}) + {d(\phi_{i-3/2}^{n})}/{\Delta z}\big)\right)\geq 0,\\
\frac{\partial\phi_{i-1/2}^{n+1}}{\partial\phi_{i+1/2}^{n}}&=
\frac{\lambda}{A_{i-1/2}}\left(-Q_{i}^{n-} +(A\gamma)_{i}\big(
-\phi_{i-1/2}^n \tilde{v}'(\phi_{i+1/2}^{n}) + {d(\phi_{i+1/2}^{n})}/{\Delta z}
\big)\right)\geq 0,\\
\frac{\partial\phi_{i-1/2}^{n+1}}{\partial\phi_{i-1/2}^{n}}&=
1 + 
\frac{\lambda}{A_{i-1/2}}\Big(
Q_{i-1}^{n-} + (A\gamma)_{i-1}\phi_{i-3/2}^n \tilde{v}'(\phi_{i-1/2}^{n})
- Q_{i}^{n+} 
- (A\gamma)_{i}\tilde{v}(\phi_{i+1/2}^{n})\\
&\quad
-\big((A\gamma)_{i-1}+(A\gamma)_{i}\big)d(\phi_{i-1/2})/\Delta z
)\Big)\\ 
&\geq 1-\lambda\left(\frac{2\|Q\|_{\infty,T}}{A_\mathrm{min}}  + M_1\big(\|\tilde{v}'\|_\infty+\|\tilde{v}\|_\infty\big)
+M_2\frac{\|d\|_\infty}{\Delta z}\right)\geq 0, 
\end{align*}
where we have used the CFL condition~\eqref{eq:CFL}.
The rest of the proof, the boundedness $0\leq\phi_{i-1/2}^n\leq 1$, follows by standard arguments, namely one verifies 
 that if $\smash{\phi_{i-1/2}^n=0}$ for all~$i$, then $\smash{\phi_{i-1/2}^{n+1}=0}$ for all~$i$  
  and likewise that if $\smash{\phi_{i-1/2}^n=1}$ for all~$i$, then $\smash{\phi_{i-1/2}^{n+1}=1}$ for all~$i$. 
   Thus, appealing to the monotonicity of the scheme, one deduces  that if 
     $\smash{0 \leq \phi_{i-1/2}^n \leq 1}$ for all~$i$, then $\smash{0 \leq \phi_{i-1/2}^{n+1} \leq 1}$ for all~$i$, 
      which proves \eqref{phibound}.

\subsection*{Proof of Theorem~\ref{thm:bounded_varphiR}}

The proof is similar to that  of Theorem~\ref{thm:bounded_phiR}.  
 We note that~\eqref{eq:psi_updateR} is again  a three-point scheme, and 
  show that   $\smash{\partial\psi_{i-1/2}^{n+1}/\partial\psi_{k-1/2}^{n}\geq 0}$ for all $i=1,\ldots, N$ and $k=i-1,i,i+1$. 
   The contributions of the terms that contain~$D$ to $\smash{\partial\psi_{i-1/2}^n/\partial\psi_{i-3/2}^n}$ and 
   $\smash{\partial\psi_{i-1/2}^n/\partial\psi_{i-3/2}^n}$ are 
\begin{equation*}
\frac{\lambda(A\gamma)_{i-1}}{A_{i-1/2}(1-\phi_{i-3/2}^n)} \dfrac{\Delta D_{i-1}^{n+}}{\Delta z}\geq 0 \quad \text{and} 
 \quad 
-\frac{\lambda(A\gamma)_{i}}{A_{i+1/2}(1-\phi_{i+1/2}^n)} \dfrac{D_{i-1}^{n-}}{\Delta z}\geq 0, 
\end{equation*}
respectively. 
Now we utilize the estimations similar to those in the previous proof (see \cite{bdmv_nhm22}) and add the terms with~$D$ to obtain
\begin{align*}
\frac{\partial\mathcalold{K}_{i-1/2}^n}{\partial\psi_{i-1/2}^{n}}& \geq
1-{\lambda}\biggl( \frac{2\|Q\|_{\infty,T}}{A_\mathrm{min}} 
+ M_1\big(\max\left\{v_\mathrm{hs}(0),\|v_\mathrm{hs}'\|_\infty\right\} 
+ \|\tilde{v}'\|_\infty \big)\\
&\quad+
\frac{1}{A_{i-1/2}(1-\phi_{i-1/2}^n)}
\biggl(
(A\gamma)_{i-1}
\frac{\Delta D_{i-1}^{n-}}{\Delta z}
-(A\gamma)_{i}
\frac{\Delta D_{i}^{n+}}{\Delta z}
\biggr)
\biggr). 
\end{align*}
It is easy to estimate the integrated terms
\begin{equation*}
\Delta D_i^n=D(\phi_{i+1/2}^n)-D(\phi_{i-1/2}^n) 
\leq\|d\|_{\infty}(1-\phi_\mathrm{c}).
\end{equation*}
Hence, we obtain with~(CFL)
\begin{align*}
\frac{\partial\mathcalold{K}_{i-1/2}^n}{\partial\psi_{i-1/2}^{n}}& \geq
1-{\lambda}\biggl( \frac{2\|Q\|_{\infty,T}}{A_\mathrm{min}} 
+ M_1\big(\max\left\{v_\mathrm{hs}(0),\|v_\mathrm{hs}'\|_\infty\right\} 
+ \|\tilde{v}'\|_\infty \big)
+ M_2(1-\phi_\mathrm{c})\frac{\|d\|_{\infty}}{\Delta z}
\biggr) \geq 0.
\end{align*}
The inequalities proven imply that $\smash{\psi_{i-1/2}^{n+1}}$ is a non-decreasing of 
 each of  $\smash{\psi_{k-1/2}^{n}}$ for $k=i-1,i,i+1$, and therefore the scheme is monotone. 
 Writing the scheme as 
 \begin{align*} 
  \psi_{i-1/}^{n+1} = \mathcalold{K}_{i-1/2} \bigl( \psi_{i-3/2}^n,\psi_{i-1/2}^n,\psi_{i+1/2}^n \bigr), 
  \end{align*} 
 assuming that $0\leq\psi_{i-1/2}^n\leq 1-\phi_{i-1/2}^n$ for all $i$ and using  that 
  \eqref{fbi1} and \eqref{fbi2} ensure that $G_i^n(1-\phi_{i-1/2}^n,1-\phi_{i+1/2}^n)=0$ (see~\cite{bdmv_nhm22}), we get 
\begin{align*}
0&\leq \frac{\lambda}{A_{i-1/2}} Q_{\mathrm{F}}^{n} \psi_{\mathrm{F}}^{n}\delta_{\mathrm{F},{i-1/2}} = \mathcalold{K}_{i-1/2}(0,0,0)\leq \psi_{i-1/2}^{n+1}\\ &=\mathcalold{K}_{i-1/2}(\psi_{i-3/2}^n,\psi_{i-1/2}^n,\psi_{i+1/2}^n)
\leq \mathcalold{K}_{i-1/2}(1-\phi_{i-3/2}^n,1-\phi_{i-1/2}^n,1-\phi_{i+1/2}^n)\\
&=1-\phi_{i-1/2}^n + \frac{\lambda}{A_{i-1/2}}\biggl( (1-\phi_{i-3/2}^n)Q_{i-1}^{n+} + (1-\phi_{i-1/2}^n)Q_{i-1}^{n-}\\ 
&\quad - (A\gamma)_{i-1}\bigg(\phi_{i-3/2}^n 
\tilde{v}(\phi_{i-1/2}^n)
-\frac{\Delta D_{i-1}^{n-}}{\Delta z}
-\frac{\Delta D_{i-1}^{n+}}{\Delta z}
\bigg) 
-(1-\phi_{i-1/2}^n)Q_{i}^{n+} - (1-\phi_{i+1/2}^n)Q_{i}^{n-}\\ 
&\quad 
+ (A\gamma)_{i}\bigg(\phi_{i-1/2}^n 
\tilde{v}(\phi_{i+1/2}^n)
-\frac{\Delta D_{i}^{n-}}{\Delta z}
-\frac{\Delta D_{i}^{n+}}{\Delta z}
\bigg)
+ {Q_{\mathrm{F}}^n\psi_{\mathrm{F}}^n\delta_{\mathrm{F},i-1/2}} \biggr).
\end{align*}
Now we use that $\Delta D_{i}^{n-}+\Delta D_{i}^{n+}=\Delta D_{i}^{n}=D(\phi_{i+1/2})-D(\phi_{i-1/2})$, that $\psi_{\mathrm{F},k}^n\leq 1-\phi_{\mathrm{F},k}^n$ and the update formula for $\phi$~\eqref{eq:phi_updateR} to obtain
\begin{align*}
\psi_{i-1/2}^{n+1}&\leq 1-\phi_{i-1/2}^{n+1}+ \frac{\lambda}{A_{i-1/2}}\big( Q_{i-1}^{n+} + Q_{i-1}^{n-} - Q_{i}^{n+} - Q_{i}^{n-} 
+ Q_{\mathrm{F}}^n\delta_{\mathrm{F},i-1/2} \big)\\
&= 1-\phi_{i-1/2}^{n+1}+ \frac{\lambda}{A_{i-1/2}}\big( Q_{i-1}^{n} - Q_{i}^{n} + Q_{\mathrm{F}}^n\delta_{\mathrm{F},i-1/2}
\big)
 =1-\phi_{i-1/2}^{n+1},
\end{align*}
since the latter parenthesis is zero irrespective of whether there is a source in the cell;  $Q_{i-1}^{n} - Q_{i}^{n} + Q_{\mathrm{F}}^n=0$, or not; $Q_{i-1}^{n} - Q_{i}^{n}=0$.

\end{document}